\documentclass[11pt, a4paper,twoside]{article}
\usepackage{amssymb,amsfonts,amsmath,amsthm, dsfont,mathrsfs}
\usepackage{graphicx,subfig}
\usepackage{multicol, multirow}
\usepackage[top = 2cm, bottom = 2cm, outer = 2.2cm, inner = 2.2cm, headsep=0.5cm]{geometry}
\usepackage{natbib}
\usepackage{enumerate}
\usepackage{color}
\usepackage[pdftex]{hyperref}
\hypersetup{
  unicode=true,
  colorlinks = true,
  linkcolor = blue,
  anchorcolor = blue,
  citecolor = blue,
  filecolor = blue,
  urlcolor = blue
}

\newtheorem{prop}{Proposition}[section]
\newtheorem{thm}{Theorem}[section]
\newtheorem{lemma}{Lemma}[section]
\newtheorem{cor}{Corollary}[section]
\theoremstyle{definition}
\newtheorem{defn}{Definition}[section]
\newtheorem{remark}{Remark}[section]
\newtheorem{example}{Example}[section]
\renewcommand{\qed}{\hfill \mbox{\raggedright \rule{.07in}{.1in}}}

\DeclareMathSizes{10.95}{10}{6}{5}

\renewcommand{\theequation}{\thesection.\arabic{equation}}

\numberwithin{equation}{section}
\numberwithin{table}{section}
\numberwithin{figure}{section}

\newcommand{\B}{\mathcal{B}}
\newcommand{\cov}{\mathrm{Cov}}
\newcommand{\E}{\mathds {E}}
\newcommand{\erf}{\operatorname{erf}}
\newcommand{\im}{\ensuremath{\mathfrak{i}{}}}
\newcommand{\Ind}[2]{\ensuremath{\mathbb {I}_{#1}{\mbox{\footnotesize $\left(#2\right)$}}}}
\newcommand{\N}{\mathds {N}}
\newcommand{\p}{\ensuremath{\mathfrak{p}}}
\renewcommand{\r}{\ensuremath{{\sf r}}}
\newcommand{\m}{\ensuremath{{\sf m}}}

\newcommand{\R}{\mathds{R}}
\newcommand{\var}{\mathrm{Var}}
\newcommand{\X}{\mbox{\tiny$X$}}
\newcommand{\Z}{\ensuremath{\mathds {Z}}}

\title{\bf  Seasonal FIEGARCH Processes \\}

\author{{\bf S\'ilvia R.C. Lopes and Taiane S. Prass\footnote{Corresponding
      author. E-mail: taianeprass@gmail.com}}  \vspace{0.2cm}\\
  Mathematics Institute \\
  Federal University of Rio Grande do Sul\\
  Porto Alegre - RS - Brazil}

\date{May 25, 2013}

\begin{document}
\allowdisplaybreaks
\maketitle

\begin{abstract}
  Here we develop the theory of seasonal FIEGARCH processes, denoted by SFIEGARCH, establishing conditions for
  the existence, the invertibility, the stationarity and the ergodicity of these processes. We analyze their
  asymptotic dependence structure by means of the autocovariance and autocorrelation functions. We also present
  some properties regarding their spectral representation.  All properties are illustrated through graphical
  examples and an application of SFIEGARCH models to describe the volatility of the S\&P500 US stock index
  log-return time series in the period from December 13, 2004 to October 10, 2009 is provided.
\end{abstract}



%
\maketitle
\section*{Introduction}

Introduced by \cite{BM1996}, FIEGARCH processes are one of the main models used to describe the volatility in
financial time series.  This class of models has not only the capability of capturing the asymmetry in the
log-returns, as in the EGARCH models, but also it takes into account the characteristic of long memory in the
volatility, as in the FIGARCH models, with the advantage of been weakly stationary.  \cite{LP2013} present a
study on the theoretical properties of these processes, including results on the volatility forecast. The
authors also analyze the finite sample performance of the quasi-likelihood estimator for four different
FIEGARCH models and present the analysis of an observed time series. The simulated study presented by
\cite{LP2013} considers the same parameters values as the ones in the models adjusted to the observed time
series considered in \cite{PL2012,PL2013}.

More recently, economists have noticed that FIEGARCH models are not fully satisfactory, specially when
modelling volatility of intra-daily financial returns.  The main discovery is that volatility of high
frequency financial time series shows long-range dependence merged with periodic behavior. According to
\cite{BEA2007}, these patterns, in the case of exchange rate returns, are generally attributed to different
openings of European, Asian and North American markets superimposed each other. Similar patterns are found in
stock markets, mainly due to the so-called time-of-day phenomena, such as market opening, closing operations,
lunch-hour and overlapping effects.  Once again, the focus is on the squared, log-squared and absolute
returns. Periodic components are represented as marked peaks at some frequencies of the time series
periodogram function and it can also be identified through a persistent cyclical behavior on the
autocorrelation function with oscillations decaying very slowly.  From the theoretical point of view,
modelling and prediction of the volatility dynamics may be seriously affected if this empirical evidence is
neglected.

\cite{BEA2007, BEA2009} introduced new GARCH-type models characterized by long memory behavior of periodic
type.  The generalized long memory GARCH (G-GARCH) introduces generalized periodic long-memory filters, based
on Gegenbauer polynomials, into the equation describing the time-varying volatility of standard GARCH models.
The periodic long-memory GARCH (PLM-GARCH) process represents a natural extension of the FIGARCH model
proposed for modelling the volatility long-range persistence.  Although periodic long memory versions of
EGARCH (PLM-EGARCH) models were also considered in \cite{BEA2009} , we feel that there are several theoretical
results related to these processes that were not yet explored. For instance, conditions for the existence,
stationarity and ergodicity are yet to be established.  Moreover, the autocovariance structure and the
spectral representation of these processes are of extreme importance in both theoretical and practical point
of view and hence, their study is an important matter.

Here we develop the theory of seasonal FIEGARCH processes, denoted by SFIEGARCH$(p,d,q)_s$, where $p,d$ and
$q$ have the same meaning as in the so-called FIEGARCH$(p,d,q)$ process and $s$ is the length of the periodic
component.  This model is similar to the PLM-EGARCH process introduced by \cite{BEA2009} but, in the
definition considered here, for any SFIEGARCH process $\{X_t\}_{t \in \mathds{Z}}$, the process
$\{\ln(\sigma_t^2)\}_{t\in\Z}$ is a SARFIMA one, where $\sigma_t^2$ is the conditional variance of $X_t$, for
all $t\in \Z$.  In particular, if $s = 1$, it is an ARFIMA$(q,d,p)$ process \citep[see][]{L2008}.  This result
is useful for establishing whether the process $\{X_t\}_{t \in \mathds{Z}}$ is well defined.

Results regarding the process $\{\ln(\sigma_t^2)\}_{t\in\Z}$ are already known in the literature and can be
found in \cite{BL2009} and references therein. Moreover, for an SFIEGARCH process the sequence of random
variables $\{\ln(\sigma_t^2)\}_{t\in\Z}$ is not directly observable and we study its characteristics only to
obtain the properties of the processes $\{X_t\}_{t \in \mathds{Z}}$ and $\{\ln(X_t^2)\}_{t\in\Z}$, which are
the observable ones.  In this work we extend the range of the parameter $d$ for the invertibility and we
present an alternative asymptotic expression for the autocovariance function
$\gamma_{\ln(\sigma_t^2)}(\cdot)$.  These results are useful to derive the exact and the asymptotic
expressions for the autocovariance and spectral density functions of the process $\{\ln(X_t^2)\}_{t\in\Z}$.

The paper is organized as follows: in Section \ref{sec:sfiegarch} we present the SFIEGARCH$(p,d,q)_s$
processes and we discuss the existence of a power series representation for the function $\lambda(z)
=\frac{\alpha(z)}{\beta(z)}(1-z^s)^{-d}$ and the asymptotic behavior of the coefficients in this
representation.  A recurrence formula to calculate those coefficients is also provided.  In Section
\ref{sec:sfiegarch} we also analyze the existence of the process $\{\ln(\sigma_t^2)\}_{t\in\Z}$ and its
invertibility property. This analysis is important to guarantee the existence of the process $\{X_t\}_{t \in
  \mathds{Z}}$ itself. Section \ref{sec:estac_erg} is devoted to study the asymptotic dependence structure of
both $\{\ln(\sigma^2_t)\}_{t \in \Z}$ and $\{\ln(X^2_t)\}_{t \in \Z}$ processes, where $\{X_t\}_{t \in
  \mathds{Z}}$ is an SFIEGARCH process. Section \ref{sec:spectralrep} presents the spectral representation of
both processes $\{\ln(\sigma^2_t)\}_{t \in \Z}$ and $\{\ln(X^2_t)\}_{t \in \Z}$. Section
\ref{sec:analysisObserved} shows an application of SFIEGARCH models to describe the volatility of the S\&P500
US stock index log-return time series in the period from December 13, 2004 to October 02, 2009.  Section
\ref{sec:conclusions} presents the final conclusions.  All proofs are presented in Appendix A.

\section{SFIEGARCH Process}\label{sec:sfiegarch}

In this section we define the \emph{Seasonal} FIEGARCH (SFIEGARCH) process which describes the volatility
varying in time, volatility clusters (known as ARCH/GARCH effects), volatility periodic long-memory and
asymmetry.  Since the existence of a solution $\{X_t\}_{t \in \mathds{Z}}$ for expression \eqref{eq:xt}
depends on the existence of the stochastic process $\{\ln(\sigma_t^2)\}_{t\in\Z}$ satisfying expression
\eqref{eq:sigmat}, we show that the random variable $\ln(\sigma_t^2) - \omega$ is finite with probability one,
for all $t\in\Z$, if and only if $d<0.5$. We show that $\{\ln(\sigma_t^2) - \omega\}_{t\in\Z}$ is an
invertible process, with respect to $\{g(Z_t)\}_{t \in \mathds{Z}}$, if and only if, $d \in (-1,0.5)$,
extending the range given in \cite{BL2009}. We also discuss the similarities between this model and the
PLM-EGARCH model, introduced by \cite{BEA2009}.

Hereafter, $\lfloor\cdot\rfloor$ and $\lceil\cdot\rceil$ denote, respectively, the floor and ceiling functions
and $\Ind{A}{\cdot}$ is the indicator function defined as $\Ind{A}{z} = 1$, if $z\in A$, and 0, otherwise.
 Whenever $T = \N$ or $T = \Z$, we define $T^*:= T\backslash\{0\}$.  Throughout the paper, given two
real/complex valued functions $f(\cdot)$ and $g(\cdot)$, $f(x) = O(g(x))$, means that $|f(x)|\leq c|g(x)|$,
for some $c>0$, as $x \rightarrow \infty$; $f(x) = o(g(x))$ means that $f(x)/g(x) \rightarrow 0$, as
$x\rightarrow \infty$; $f(x) \sim g(x)$ means that $f(x)/g(x) \rightarrow 1$, as $x\rightarrow \infty$. We
also say that $f(x) \approx g(x)$, as $x \rightarrow \infty$, if for any $\varepsilon > 0$, there exists $x_0
\in \R$ such that $|f(x) - g(x)| < \varepsilon$, for all $x \geq x_0$.  Similar definitions can be obtained
upon replacing the functions $f(\cdot)$ and $g(\cdot)$ by sequences of real numbers $\{a_k\}_{k\in\N}$ and
$\{b_k\}_{k\in\N}$ or if one considers any constant $a$ or $-\infty$ instead of $\infty$.

\begin{defn}\label{defn:sfiegarch}
  Let $\{X_t\}_{t \in \mathds{Z}}$ be the stochastic process defined by the expressions
  \begin{align}
    X_t &= \sigma_tZ_t,\label{eq:xt}\\
    \ln(\sigma_t^2) &= \omega +
    \frac{\alpha(\B)}{\beta(\B)}(1-\B^s)^{-d}g(Z_{t-1}), \quad \mbox{ for all
      $t\in\Z$,}\label{eq:sigmat}
  \end{align}
  \noindent where $\{Z_t\}_{t \in \mathds{Z}}$ is a sequence of i.i.d.\! random variables, with zero mean and
  variance equal to one, $g(\cdot)$ is defined by
  \begin{equation}
    g(Z_{t}) = \theta Z_t +
    \gamma\big[|Z_t|-\E\big(|Z_t|\big)\big], \quad \mbox{with } \theta, \gamma \in
    \R, \quad \mbox{for all } t\in\Z,\label{eq:gzt}
  \end{equation}
  \noindent $\omega \in \R$, $\B$ is the backward shift operator defined by $\B^{sk}(X_t) = X_{t-sk}$, for all
  $s,k\in\N$, $\alpha(\cdot)$ and $\beta(\cdot)$ are, respectively, polynomials of order $p$ and $q$, with no
  common roots defined by
  \begin{equation}
    \alpha(z) = \sum_{i=0}^{p}(-\alpha_i)z^i     \quad \mbox{and} \quad
    \beta(z) = \sum_{j=0}^{q}(-\beta_j)z^j,
    \label{eq:alphabeta}
  \end{equation}
  \noindent with $\alpha_0 = -1 = \beta_0$, and $\beta(z)\neq 0$ in the closed disk $\{z: |z|\leq 1\}$,
  $d\in\R$ is the differencing parameter, $s\in\N ^*$ is the length of the periodic component, $(1-\B^s)^{-d}$
  is the seasonal difference operator, defined by its Maclaurin series expansion, namely,
  \begin{equation}\label{eq:pik}
    (1-\B^s)^{-d} =
    \sum_{k=0}^{\infty}\frac{\Gamma(k+d)}{\Gamma(k+1)\Gamma(d)}(\B^s)^k
    :=\sum_{k=0}^{\infty} \delta_{-d,k}\B^{sk} := \sum_{k = 0}^
    \infty\pi_{d,k}\B^k ,
  \end{equation}
  \noindent where $\Gamma(\cdot)$ is the Gamma function, $\pi_{d,k} := 0$, if $k/s \notin \N$, and $\pi_{d,sj}
  = \delta_{-d,j}:= \frac{\Gamma(k+d)}{\Gamma(k+1)\Gamma(d)}$, for all $j\in\N$.  Then, $\{X_t\}_{t \in
    \mathds{Z}}$ is a seasonal FIEGARCH process, with seasonal period $s$ and differencing parameter $d$, denoted
  by SFIE\-GARCH$(p,d,q)_s$.
\end{defn}

\begin{remark}
  The assumptions that $\beta(z)\neq 0$ in the closed disk $\{z: |z|\leq 1\}$ and that $\alpha(\cdot)$ and
  $\beta(\cdot)$ have no common roots guarantee that the operator $\frac{\alpha(\B)}{\beta(\B)}$ is well
  defined.
\end{remark}

\begin{figure}[!htb]
  \centering
  \subfloat[]{\includegraphics[width = 0.45\textwidth]{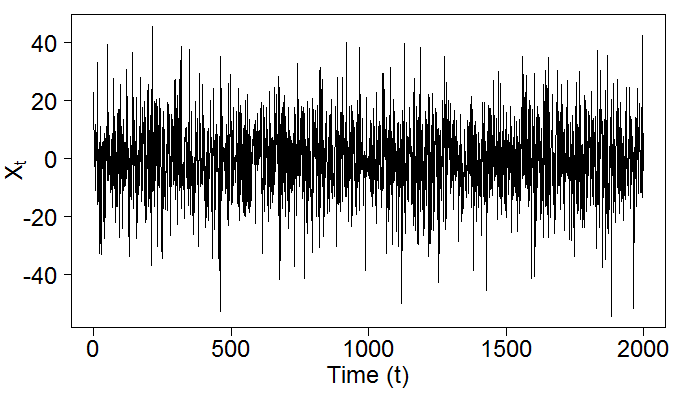}}\hspace{0.5cm}
  \subfloat[]{\includegraphics[width = 0.45\textwidth]{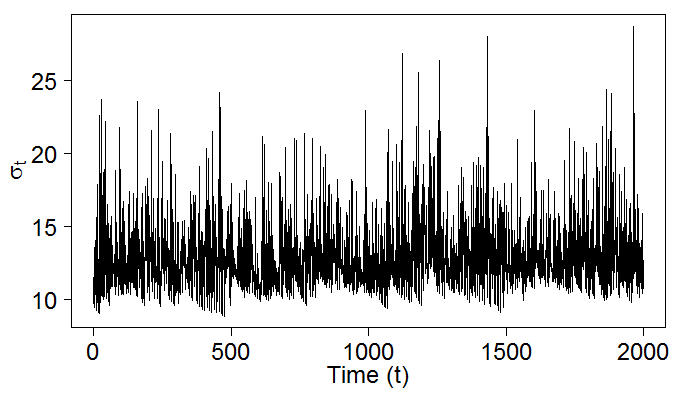}}
  \caption{ Samples from an SFIEGARCH$(0,d,0)_s$ processes, with $n = 2000$ observations, considering $Z_0 \sim
    \mathcal{N}(0,1)$, $\omega = 5.0$, $\theta = -0.25$, $\gamma = 0.24$, $d = 0.35$ and $s = 6$.  Panel (a) shows
    the time series $\{X_t\}_{t=1}^{n}$. Panel  (b) presents the time series $\{\sigma_t\}_{t=1}^{n}$, where $\sigma_t$
    is the conditional standard deviation of $X_t$, for all $t\in \{1,\cdots, n\}$. }\label{fig:sfiegsim}
\end{figure}

\begin{example}
  Figure \ref{fig:sfiegsim} presents a simulated SFIEGARCH$(0,d,0)_s$ time series $\{X_t\}_{t=1}^n$ and its
  conditional standard deviation $\{\sigma_t\}_{t=1}^n$, defined by expressions \eqref{eq:xt} and
  \eqref{eq:sigmat}.  For these graphs,   $Z_0 \sim \mathcal{N}(0,1)$, $\omega = 5.0$,   $\theta= -0.25$, $\gamma =
  0.24$, $d = 0.35$ and $s = 6$.
\end{example}

\begin{remark}
  In this work we consider the case where the conditional variance $\sigma_t^2$ is defined through expression
  \eqref{eq:sigmat}, with $\E(Z_0) = 0$, $\var(Z_0) = \E(Z_0^2) = 1$ and the function $g(\cdot)$ defined by
  expression \eqref{eq:gzt}. However the results presented here can be easily extended if one considers
  $\var(Z_0) = \sigma^2 \neq 1$ and replaces $g(\cdot)$ by any measurable function
  satisfying $\var(g(Z_0)) < \infty$.
\end{remark}

Observe that the series expansion of the operator $(1-\B^s)^d$ is obtained upon replacing $-d$ by $d$ in
expression \eqref{eq:pik}. Moreover, when $d\in\N$, $(1-\B^s)^d$ is merely the seasonal difference operator
$1-B^s$ iterated $d$ times. Thus, one can easily see that an equivalent definition for SFIEGARCH process is
given if one replaces expression \eqref{eq:sigmat} by
\begin{equation}\label{eq:sigmat01}
  \beta(\B)(1-\B^s)^{d}(\ln(\sigma_t^2) - \omega) =
  \alpha(\B)g(Z_{t-1}), \quad \mbox{for all } t\in\Z.
\end{equation}
\noindent This expression is similar to the one in the definition of the PLM-EGARCH process, presented by
\cite{BEA2009}.  For a PLM-EGARCH$(p,m,d,q,s)$, the conditional variance $\sigma_t^2$ of $X_t$ is defined
through the equation
\begin{equation}\label{eq:plmegarch}
  (1-\B^s)^d\phi(\B)\big(\ln(\sigma_t^2)-\omega)= a(\B)Z_t +
  c(\B)\big(|Z_t|-\E(|Z_t|)\big), \quad \mbox{for all } t\in\Z,
\end{equation}
\noindent where $a(z) = \sum_{k=1}^pa_kz^k$ and $c(z) = \sum_{l=1}^m c_lz^l$ are polynomials of order $p$ and
$m$, respectively, $\phi(z) = \sum_{j=0}^{q-s} \phi_jz^j$ is a polynomial of order $q-s$, which satisfies
$(1-\B^s)^d\phi(\B) = 1-b(\B)$, where $b(z) = \sum_{i=1}^qb_iz^i$ is a polynomial of order $q$.

Notice that in the PLM-EGARCH, the polynomials $a(\cdot)$ and $c(\cdot)$ do not necessarily have the same
order.  Also, it is easy to see that, by setting
\begin{equation}\label{eq:relation}
  \begin{split}
    \phi_i & := -\beta_i, \quad \mbox{for all } i=0,\cdots, q, \\
    a_{j+1} & := -\theta\alpha_j \quad \mbox{and} \quad c_{j+1}:=
    -\gamma\alpha_j, \quad \mbox{for all } j=0,\cdots, p,
  \end{split}
\end{equation}
\noindent one can rewrite the right hand side of expression \eqref{eq:sigmat01} as the right hand side of
\eqref{eq:plmegarch}.  Under this point of view, the PLM-EGARCH model seems more general than the SFIEGARCH
one.  On the other hand, the left hand side of expression \eqref{eq:sigmat01} is more general then the left
hand side of \eqref{eq:plmegarch}. This is so because in the SFIEGARCH model, no restriction is made in the
order of the product $\beta(z)(1-z^s)^{d}$, allowing for the parameter $d$ to be fractional.

\begin{remark}
  It is immediate that, if $\{\ln(\sigma_t^2)\}_{t\in\Z}$ is a stationary process with finite mean, then
  $\omega = \E(\ln(\sigma_t^2))$, for all $t\in\Z$.  Also, if $d= 0$, we have the EGARCH$(p,q)$ model proposed
  by \cite{N1991} and, if $s=1$, we have the FIEGARCH$(p,d,q)$ process defined by \cite{BM1996}. A study on the
  theoretical properties of FIEGARCH$(p,d,q)$ process are presented in \cite{LP2013}.
\end{remark}

From Definition \ref{defn:sfiegarch}, one easily concludes that the existence of the stochastic process
$\{X_t\}_{t \in \mathds{Z}}$ depends on the existence of the stochastic process $\{\ln(\sigma_t^2)\}_{t\in\Z}$
which satisfies equation \eqref{eq:sigmat}.  The existence of a solution for equation \eqref{eq:sigmat} is
discussed in the sequel.

From now on, let $\lambda(\cdot)$ be the polynomial defined as
\begin{equation}
  \lambda(z) := \frac{\alpha(z)}{\beta(z)}(1-z^s)^{-d} =
  \sum_{k=0}^{\infty}\lambda_{d,k}z^k,
  \quad |z| < 1, \label{eq:polilambda}
\end{equation}
\noindent where $\alpha(\cdot)$ and $\beta(\cdot)$ are defined in \eqref{eq:alphabeta}.  Notice that, by
definition, $\beta(z)$ has no roots in the closed disk $\{z:|z|\leq 1\}$, and $\alpha(\cdot)$ and
$\beta(\cdot)$ have no common roots. Therefore, the function $\lambda(z)$ is analytic in the open disc
$\{z:|z|<1\}$ and, if $d \leq 0$, in the closed disk $\{z:|z|\leq 1\}$.  So, it has a unique power series
representation and the operator given in \eqref{eq:sigmat} can be rewritten as
$\frac{\alpha(\B)}{\beta(\B)}(1-\B^s)^{-d} = \sum_{k=0}^{\infty}\lambda_{d,k}\B^k = \lambda(\B)$.  This
representation is more convenient and will be used from now on.  In the following, we analyze the asymptotic
behavior of the coefficients $\lambda_{d,k}$, for all $k\in\Z$, defined by expression \eqref{eq:polilambda}.
This result is fundamental for proving the results regarding the existence, invertibility, stationarity and
ergodicity of SFIEGARCH processes.

It is immediate that, if $p\geq 0$ and $q=0$, one can rewrite \eqref{eq:polilambda}
as,
\[
\lambda(z)=  \left\{
  \begin{array}{cl}
    \displaystyle (1-z^s)^{-d} = \sum_{k=0}^\infty \pi_{d,k}z^k, & \mbox{if }  p = 0 =q;\vspace{0.2cm}\\
    \displaystyle \alpha(z)(1-z^s)^{-d} =   \sum_{k=0}^{\infty}\Big[\sum_{i=0}^{\min\{p,k\}}(-\alpha_i\pi_{d,k-i})\Big]z^k,
    & \mbox{if }  p >0 \mbox{ and } q = 0.
  \end{array}
\right.
\]
Thus, for all $r\in\{0,\cdots,s-1\}$ and all $k\in\N$,
\[
\lambda_{d,sk+r} = \pi_{d,sk+r}, \quad \mbox{if} \quad p=0=q,
\]
and, whenever $p>0$ and $q=0$, \vspace{-1\baselineskip}
\begin{align*}
  \lambda_{d,sk+r} =
  -\hspace{-15pt}\sum_{j=0}^{\min\{p,sk+r\}}\hspace{-15pt}\alpha_j\pi_{d,sk+r-j}
  = \left\{
    \begin{array}{cl}
      0,   & \mbox{if }   p < r; \vspace{0.2cm}  \\
      \displaystyle-\hspace{-15pt}\sum_{j = 0}^{\min\{\lfloor\frac{p-r}{s}\rfloor,k\}}\hspace{-15pt}\alpha_{sj+r}\pi_{d,sk-sj}, & \mbox{otherwise.}
    \end{array}
  \right.\nonumber
\end{align*}
Consequently, given $\r>0$,
\[
\sum_{k=0}^\infty|\lambda_{d,k}|^{\r} < \infty \quad \mbox{if and only if }
\quad \sum_{k=0}^\infty |\pi_{d,k}|^{\r} < \infty.
\]
Theorem \ref{thm:convergenceOrder} bellow shows that this result also holds in the general case $p\geq 0$ and
$q>0$.  The proofs of all results stated in this work are given in the Appendix.

\begin{remark}
  By Stirling's formula and from lemma 3.1 in \cite{KT1995}, one easily concludes that
  \begin{equation}\label{eq:stirling}
    \pi_{d,sk}  := \frac{\Gamma(k+d)}{\Gamma(d)\Gamma(k+1)} =
    \dfrac{1}{\Gamma(d)k^{1-d}} + O(k^{d-2}) \sim
    \frac{1}{\Gamma(d)k^{1-d}}, \quad \mbox{as } k\to \infty.
  \end{equation}
  Since  (integral convergence test)
  \[
  \sum_{k=1}^m k^{-\r(d-1)} \leq \int_{1}^m x^{-\r(d-1)} \mathrm{d}x = \frac{1}{1-r(d-1)}\Big|_{1}^m, \quad \mbox{for any } \r > 0,
  \]
  converges to a finite constant as $m \to \infty$ if and only if $1 - (1-d)\r<0$, it follows that
  $\sum_{k=0}^\infty |\pi_{d,k}|^{\r} < \infty$ if and only if $(1-d)\r > 1$.
\end{remark}

\begin{thm}\label{thm:convergenceOrder}
  Let $\lambda_{d,k}$, for $k\in\N$, be the coefficients of the polynomial $\lambda(\cdot)$, given by
  expression \eqref{eq:polilambda}.  Let $f(z):= \frac{\alpha(z)}{\beta(z)} = \sum_{k=0}^\infty f_kz^k$. Then,
  for each $r\in\{0,\cdots,s-1\}$ and any $\nu>0$, one has
  \begin{equation}
    \lambda_{d,sk+r} = \pi_{d,sk}\mathscr{K}(sk+r) + o(k^{-\nu}),
    \qquad \qquad\mbox{as} \ k \rightarrow \infty,  \label{eq:lambdaapprox01}
  \end{equation}
  where $\mathscr{K}(\cdot)$ satisfies $\displaystyle \lim_{k\to\infty}\sum_{r=0}^{s-1}\mathscr{K}(sk+r) =
  \frac{\alpha(1)}{\beta(1)} $. Thus, $\displaystyle \sum_{r=0}^{s-1}\lambda_{d,sk+r} \sim
  \pi_{sk}\frac{\alpha(1)}{\beta(1)}, \quad \mbox{as } k\rightarrow \infty.$
\end{thm}

Theorem \ref{thm:convergenceOrder2} presents an alternative asymptotic representation for the coefficients
$\lambda_{d,k}$, as $k$ goes to infinity.  While expression \eqref{eq:lambdaapprox01} is more convenient for
proving the asymptotic behavior of $\gamma_{\ln(\X_t^2)}(\cdot)$ (see Theorem \ref{thm:resultsXt}), expression
\eqref{eq:lambdaapprox} is useful for simulation purpose (see Remark \ref{remark:lambda}).

\begin{thm}\label{thm:convergenceOrder2}
  Let $\lambda_{d,k}$, for $k\in\N$, be the coefficients of the polynomial $\lambda(\cdot)$, given by
  expression \eqref{eq:polilambda}, with $d<0.5$.  Then, for each $r\in\{0,\cdots,s-1\}$, one can write
  \begin{align}
    \lambda_{d,sk+r} &= \frac{1}{\Gamma(d)k^{1-d}}\frac{\alpha(1)}{\beta(1)} -
    \frac{1}{\Gamma(d)k^{1-d}}\sum_{j =
      0}^{sk+r}f_{j}\Ind{\R\backslash\N}{\frac{|j-r|}{s}} + O(k^{d-2})\label{eq:lambdaapprox}\\
    & = O(k^{d-1}) + O(k^{d-1})\Ind{\N\backslash\{0,1\}}{s} + O(k^{d-2}) ,
    \qquad \qquad\mbox{as} \ k \rightarrow \infty, \nonumber
  \end{align}
  where $f(z)= \frac{\alpha(z)}{\beta(z)} = \sum_{k=0}^\infty f_kz^k$.
\end{thm}

\begin{remark}\label{remark:lambda}
  From Theorem \ref{thm:convergenceOrder} one observes that $\lambda_{d,sk+r}$ behaves asymptotically as the
  coefficient $\pi_{d,sk}$, as $k$ goes to infinity. This property is very useful to prove the results stated in
  Section \ref{thm:est_ergod}.  On the other hand, from Theorem \ref{thm:convergenceOrder2},
  \[
  \lambda_{d,k} \approx  \frac{s^{1-d}}{\Gamma(d)k^{1-d}}\frac{\alpha(1)}{\beta(1)},\quad \mbox{as} \quad k\to \infty.
  \]
  This approximation has a closed formula which also takes into account the magnitude of
  $\frac{\alpha(1)}{\beta(1)}$.  Although this is a rough approximation, it can be used to estimate a truncation
  point $m$ for $\lambda(\cdot)$ in Monte Carlo simulation studies. That is, given $\varepsilon > 0$, if one
  chooses $m \gg s${\scriptsize $
    \left[\frac{1}{|\Gamma(d)|\varepsilon}\frac{\alpha(1)}{\beta(1)}\right]^{\frac{1}{1-d}}$}, one gets
  $|\lambda_{d,m}| < \varepsilon$.
\end{remark}

In the following proposition we present a recurrence formula to calculate the coefficients $\lambda_{d,k}$,
for all $k\in\N$.  This recurrence formula is very useful in Monte Carlo simulation studies.

\begin{prop}\label{prop:coefs}
  Let $\lambda(\cdot)$ be the polynomial defined by \eqref{eq:polilambda}. Suppose $\alpha(\cdot)$ and
  $\beta(\cdot)$ have no common roots and $\beta(z) \neq 0$ in the closed disk $\{z: |z|\leq 1\}$. Then, the
  coefficients $\lambda_{d,k}$, for all $k\in\N$, are given by
  \[
  \lambda_{d,k} = \left\{
    \begin{array}{ccc}
      1, &\mbox{if } &  k=0,\vspace{.3cm}\\
      \displaystyle-\alpha_k + \sum_{i=0}^{k-1}\lambda_{d,i}\Big(\sum_{j=0}^{(k-i)\land q}
      \delta_{d,\frac{k-i-j}{s}}^*\beta_j\Big), &  \mbox{if } &   k \leq p;\vspace{.3cm}\\
      \displaystyle \sum_{i=0}^{k-1}\lambda_{d,i}\Big(\sum_{j=0}^{(k-i)\land q}
      \delta_{d,\frac{k-i-j}{s}}^*\beta_j\Big), &  \mbox{if } &   k > p,
    \end{array}
  \right.
  \]
  \noindent where $(k-i)\land q = \min\{k-i,q\}$ and, by definition,
  \begin{equation}\label{eq:newcoefs}
    \delta_{d,m}^* = \left\{
      \begin{array}{ccc}
        \delta_{d,m}, &  \mbox{if} &m \in \N,\vspace{.3cm}\\
        0,    &  \mbox{if} &  m \notin \N,
      \end{array}
    \right.
  \end{equation}
  with $\delta_{d,m}$, for all $m\in\N$, defined in \eqref{eq:pik}.
\end{prop}

The following proposition presents some properties of the stochastic process $\{g(Z_t)\}_{t \in
  \mathds{Z}}$. Although the proof is straightforward and follows immediately from the fact that $\{Z_t\}_{t \in
  \mathds{Z}}$ is a sequence of i.i.d.\! random variables, the proposition is fundamental to establish the
result in Lemma \ref{lemma:SARFIMA}, Corollary \ref{lemma:welldefined} and Theorem \ref{thm:invertibility}.

\begin{prop}\label{prop:gzt}
  Let $g(\cdot)$ be the function defined by \eqref{eq:gzt} and $\{Z_t\}_{t \in \mathds{Z}}$ be a sequence of
  i.i.d.\! random variables, with zero mean and variance equal to one.  Then, $\{g(Z_t)\}_{t \in \mathds{Z}}$ is
  a white noise process with i.i.d.\! random variables and its variance $\sigma_g^2$ is given by
  \begin{equation}\label{eq:sigma2gzt}
    \sigma_g^2 =  \theta^2 + \gamma^2 - [\gamma \E(|Z_0|)]^2 + 2\theta\gamma\E(Z_0|Z_0|).
  \end{equation}
  \noindent Moreover, the stochastic process $\{g(Z_t)\}_{t \in \mathds{Z}}$ is stationary (weak and strictly)
  and ergodic.
\end{prop}

\begin{remark}
  Henceforth, GED shall stands for the so-called Generalized Error Distribution \citep[see][]{N1991}. Whenever
  we consider $Z_0 \sim \mbox{GED}(\nu)$, and $\nu$ is the tail-thickness parameter, we assume that the random
  variable was normalized to have mean zero and variance equal to one.
\end{remark}

\begin{remark}
  If the random variable $Z_0$ is symmetric, then $\E(Z_0|Z_0|) = 0$ and \eqref{eq:sigma2gzt} is replaced by
  $\sigma_g^2 = \theta^2 + \gamma^2(1 - \E(|Z_0|)^2) $.  Besides, if $Z_0 \sim \mbox{GED}(\nu)$, then
  \[
  \sigma_g^2 = \theta^2 + \gamma^2\left( 1 - \frac{
      \Gamma(2/\nu)}{\sqrt{\Gamma(1/\nu)\Gamma(3/\nu)}}\right), \quad \mbox{for
    any} \quad \nu > 1.
  \]
  \noindent If  $\nu = 2$, one has the Gaussian case, that is,   $Z_0\sim \mathcal{N}(0,1)$ and $\sigma_g^2 = \theta^2 + \gamma^2( 1 - 2/\pi)$.
\end{remark}

\begin{figure}[!htb]
  \vspace{-1\baselineskip}\centering
  \subfloat[]{\includegraphics[width = 0.33\textwidth]{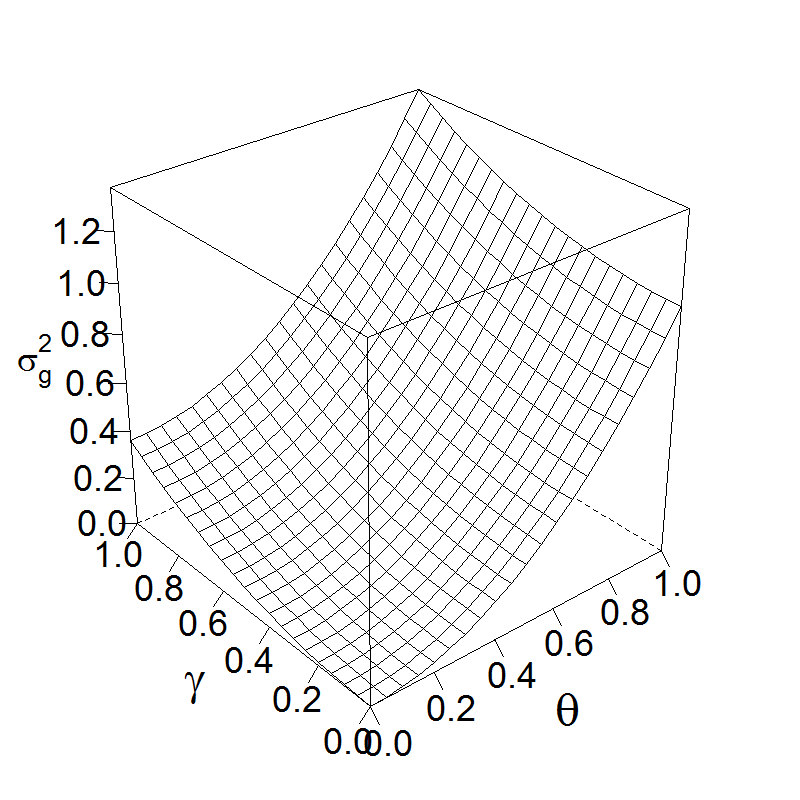}}
  \subfloat[]{\includegraphics[width = 0.33\textwidth]{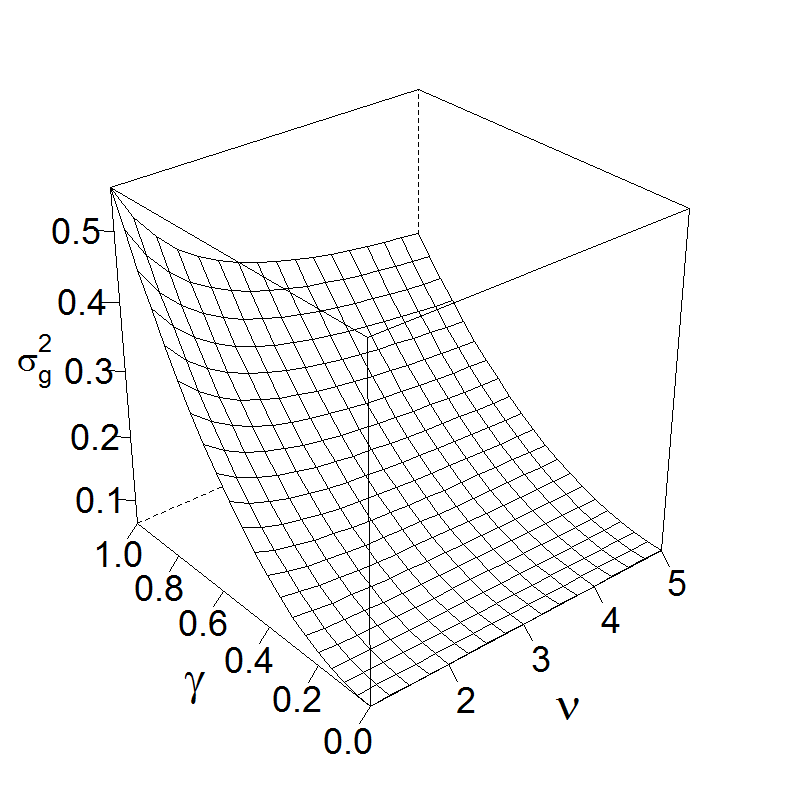}}
  \subfloat[]{\includegraphics[width = 0.33\textwidth]{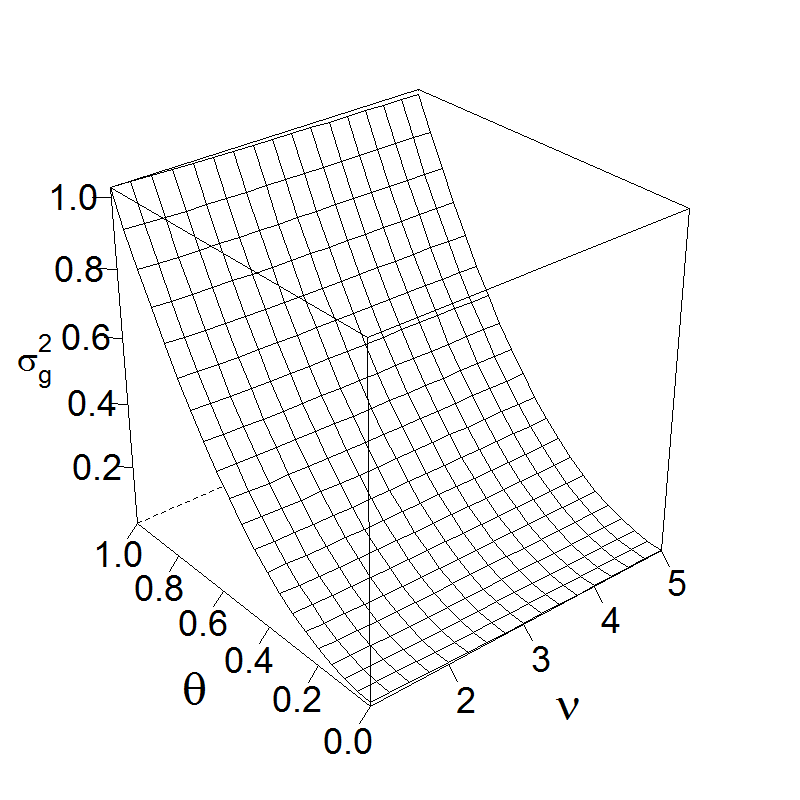}}
  \caption{ This figure presents the behavior of $\sigma_g^2$, the variance of the process $\{g(Z_t)\}_{t\in
      \Z}$, as a function of parameters $\theta$, $\gamma$ and $\nu$, when $Z_0 \sim \mbox{GED}(\nu)$.  Panel
    (a) considers $\sigma_g^2$ as a function of $\theta$ and $\gamma$ when $\nu = 2$, that is, when $Z_0\sim
    \mathcal{N}(0,1)$. Panel (b) shows $\sigma_g^2$ as a function of $\gamma$ and $\nu$ when $\theta =
    0.25$. Panel (c)  considers $\sigma_g^2$ as a function of $\theta$ and $\nu$ when $\gamma = 0.24$.}\label{fig:sigma2g}
\end{figure}

\begin{figure}[!ht]
  \centering
  \subfloat[]{ \includegraphics[width = 0.4\textwidth]{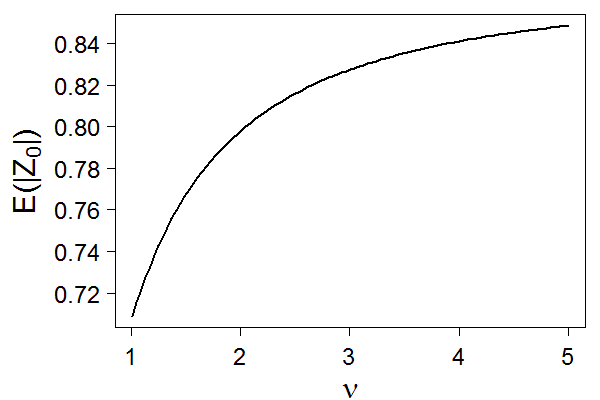}}\hspace{1cm}
  \subfloat[]{ \includegraphics[width = 0.4\textwidth]{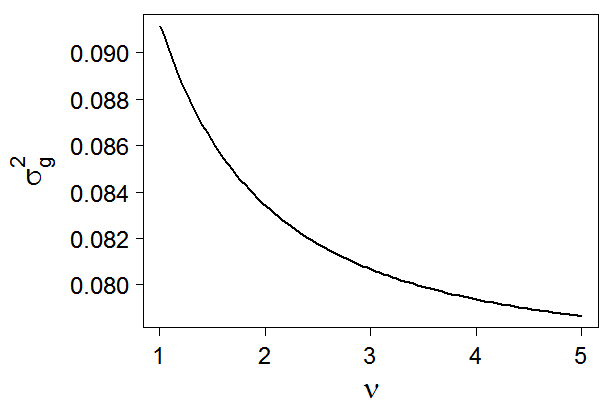}}
  \caption{ This figure considers $Z_0 \sim \mbox{GED}(\nu)$.  Panel (a) presents the graph of $\E(|Z_0|)$ as a
    function of $\nu$. Panel (b) shows the the graph of $\E(g(Z_0)^2) = \sigma_g^2$, the variance of the process
    $\{g(Z_t)\}_{t\in \Z}$, as a function of $\nu$ when $\theta = 0.25$ and $\gamma = 0.24$. }\label{fig:ezo}
\end{figure}

\begin{example}
  Figure \ref{fig:sigma2g} (a) shows the graphs of $\sigma_g^2$ as a function of $\theta$ and $\gamma$, when
  $Z_0\sim \mathcal{N}(0,1)$. Figures \ref{fig:sigma2g} (b) and (c) consider $Z_0\sim \mbox{GED}(\nu)$, with
  $\nu >1$ ($\nu = 2$ corresponds to the Gaussian case), and present the graphs of $\sigma_g^2$, respectively,
  as a function of $\gamma$ and $\nu$, for $\theta = 0.25$, and as a function of $\theta$ and $\nu$, for $\gamma
  = 0.24$.  Notice that for these graphs, $Z_0$ is a symmetric random variable, so we only consider positive
  values of $\theta$ and $\gamma$.

  From Figure \ref{fig:sigma2g} one observes that, although $\sigma_g^2$ is increasing in both $\theta$ and
  $\gamma$, for any $\nu >1$, it varies faster as $\theta$ increases than when $\gamma$ does (notice the scales
  for the ordinate axis). Moreover, for each $\gamma$ and $\theta$ fixed, $\sigma_g^2$ is decreasing in
  $\nu$. This is expected since $\E(|Z_0|)$ is increasing for $\nu\in(1,5.00]$.  This fact is illustrated in
  Figure \ref{fig:ezo} where the graphs of $\E(|Z_0|)$ and $\sigma_g^2$ as functions of $\nu$ are presented. In
  Figure \ref{fig:ezo} (b) we fixed $\theta = 0.25$ and $\gamma = 0.24$.  From this figure it is easy to see that
  $\sigma_g^2$ is indeed decreasing in $\nu$.
\end{example}

Lemma \ref{lemma:welldefined} provides the necessary and sufficient conditions for the existence of the
process $\{\ln(\sigma_t^2)\}_{t\in\Z}$.

\begin{lemma} \label{lemma:welldefined}
  Suppose that $\{Z_t\}_{t \in \mathds{Z}}$ is a sequence of i.i.d.\! random variables with zero mean and
  variance equal to one. Let $\{g(Z_t)\}_{t \in \mathds{Z}}$ be the process defined by \eqref{eq:gzt}, $\omega$
  be a real constant and $\lambda(\cdot)$ be the operator defined by \eqref{eq:polilambda}. Define
  \begin{equation}
    \ln(\sigma_t^2) - \omega = \sum_{k=0}^{\infty}\lambda_{d,k}g(Z_{t-1-k}), \quad \mbox{for all } t\in\Z.\label{eq:sigma_sum}
  \end{equation}
  \noindent Thus, the series \eqref{eq:sigma_sum} is well defined and converges a.s.\! if and only if $d <
  0.5$. Moreover, the series \eqref{eq:sigma_sum} converges absolutely a.s.  for $d\leq 0$.
\end{lemma}

Corollaries \ref{lemma:SARFIMA} and \ref{cor:xfinite} follow immediately from Lemma \ref{lemma:welldefined}
and show, respectively, that the process $\{\ln(\sigma_t^2)\}_{t\in\Z}$ is a causal SARFIMA process and that
$X_t$ is finite with probability one, for all $t\in\Z$.  We emphasize that, causality and invertibility are
defined in terms of convergence in the linear space $L^2$ \citep[see][]{P2007} and not in the linear space
$L^1$ \citep[as in][]{BD1991}.  The same approach is considered in \cite{B1985} and \cite{BP2007}.

\begin{cor}\label{lemma:SARFIMA}
  Let $\{\ln(\sigma_t^2)\}_{t\in\Z}$ be the stochastic process defined by expression \eqref{eq:sigmat}, with
  $d<0.5$. Then, $\{\ln(\sigma_t^2)\}_{t\in\Z}$ is a causal \emph{SARFIMA}$(p,0,q)\times(0,D,0)_s$ process, with
  $D = d$.
\end{cor}

\begin{cor}\label{cor:xfinite}
  Let $\{X_t\}_{t \in \mathds{Z}}$ be an \emph{SFIEGARCH$(p,d,q)_s$} process, with $d<0.5$. Then, the random
  variable $X_t$ is finite with probability one, for all $t\in\Z$.
\end{cor}

\cite{BL2009} show that a SARFIMA$(p,d,q)\times(P,D,Q)_s$ process is invertible whenever $|d+D| <
0.5$. Moreover, it is usually stated that an ARFIMA$(p,d,q)$ process is invertible for $|d| < 0.5$ \citep[see
for instance][]{H1981,BD1991}.  However, \cite{B1985} proves that, for an ARFIMA$(0,d,0)$, this range can be
extended to $d\in(-1, 0.5)$. \cite{BP2007} show that this result actually holds for any ARFIMA$(p,d,q)$.
Although the spectral density function of $\{\ln(\sigma_t^2) - \omega\}_{t\in\Z}$ does not satisfy all
conditions imposed in \cite{BP2007}, with some modifications in their proof, we show here that the results
still holds for a SARFIMA$(p,0,q)\times(0,D,0)_s$ process (see Theorem \ref{thm:invertibility}).

\begin{thm}\label{thm:invertibility}
  Let $\{X_t\}_{t \in \mathds{Z}}$ be an \emph{SFIEGARCH}$(p,d,q)_s$, defined by \eqref{eq:xt} and
  \eqref{eq:sigmat}, with $\gamma$ and $\theta$, not both equal to zero. Assume that $\alpha(z)\neq 0$, for all
  $|z|\leq 1$.  Let $Y_{t} := \ln(\sigma_{t}^2) - \omega$, for all $t\in\Z$. Then,
  \[
  \lim_{m\rightarrow \infty} \E\left(\left|\sum_{k=0}^m
      \tilde\lambda_{d,k}Y_{t-k} - g(Z_{t-1})\right|^{\p}\right) = 0, \quad
  \mbox{for all} \quad 0 <\p \leq 2,
  \]
  \noindent if and only if $d \in (-1,0.5)$, where $\{\tilde\lambda_{d,k}\}_{k\in\N}$ is the sequence of
  coefficients in the series expansion of $\tilde\lambda(z) := \lambda^{-1}(z)$, for $|z| < 1$, that is,
  \[
  \sum_{k=0}^{\infty}\tilde\lambda_{d,k}z^k = \tilde \lambda(z)
  :=  \lambda^{-1}(z)  = \dfrac{\beta(z)}{\alpha(z)}(1-z^s)^d, \quad |z|<1.
  \]
\end{thm}

\section{Stationarity and Ergodicity}\label{sec:estac_erg}

Here we show that for any SFIEGARCH$(p,d,q)_s$, with $\theta$ and $\gamma$ not both equal to zero and $d<0.5$,
the processes $\{X_t\}_{t \in \mathds{Z}}$ and $\{\sigma^2_t\}_{t\in\Z}$ are strictly stationary and ergodic
processes.  We also prove that if $\E([\ln(Z_0^2)]^2) < \infty$, the process $\{\ln(X_t^2)\}_{t\in\Z}$ is well
defined and it is stationary (weakly and strictly) and ergodic.  Weakly stationarity of the processes
$\{X_t\}_{t \in \mathds{Z}}$ and $\{\sigma_t^2\}_{t\in\Z}$ is also discussed.  For any stationary SFIEGARCH
process, we give the expressions for the autocovariance and autocorrelation functions of
$\{\ln(\sigma_t^2)\}_{t\in\Z}$ and $\{\ln(X_t^2)\}_{t\in\Z}$ and study their relation and asymptotic
behavior. We also provide expression for the asymmetry (also known as skewness) and kurtosis measures for any
stationary SFIEGARCH process $\{X_t\}_{t \in \mathds{Z}}$.

Lemma \ref{lemma:lsstationary} presents the conditions for the stationarity of the SARFIMA process
$\{\ln(\sigma_t^2)-\omega\}_{t\in\Z}$. This lemma is useful to prove Theorem \ref{thm:est_ergod} that presents
results on the stationarity of the processes $\{X_t\}_{t \in \mathds{Z}}$ and $\{\ln(X_t^2)\}_{t\in\Z}$.

\begin{lemma}\label{lemma:lsstationary}
  Let $\{\ln(\sigma_t^2)-\omega\}_{t\in\Z}$ be defined by \eqref{eq:sigmat}, with $\gamma$ and $\theta$ not
  both equal to zero. If $d<0.5$, the stochastic process $\{\ln(\sigma_t^2)-\omega\}_{t\in\Z}$ is stationary
  (strictly and weakly) and ergodic.
\end{lemma}

\begin{cor}\label{cor:sigmastationary}
  If $d<0.5$, the stochastic process $\{\sigma^2_t\}_{t\in\Z}$ is strictly stationary and ergodic.
\end{cor}

Theorem \ref{thm:est_ergod} shows that both  processes $\{X_t\}_{t \in \mathds{Z}}$ and
$\{\ln(X_t^2)\}_{t\in\Z}$ are strictly stationary and ergodic, whenever $d<0.5$
and $\E(|\ln(Z_0)|)<\infty$, regardless the distribution of the random variable
$Z_0$. This theorem also provides the necessary condition for
$\{\ln(X_t^2)\}_{t\in\Z}$ to be a weakly stationary process.

\begin{thm}\label{thm:est_ergod}
  Let $\{X_t\}_{t \in \mathds{Z}}$ be an \emph{SFIEGARCH}$(p,d,q)_s$, defined by \eqref{eq:xt} and
  \eqref{eq:sigmat}.  Suppose that $\gamma$ and $\theta$, given in \eqref{eq:gzt}, are not both equal to zero.
  If $d<0.5$
  \begin{enumerate}[{\rm \bf \hspace{-0.2cm}i)}]
  \item the stochastic process $\{X_t\}_{t \in \mathds{Z}}$ is strictly stationary and ergodic.

  \item if $\E(|\ln(Z_0^2)|) < \infty$, the process $\{\ln(X_t^2)\}_{t\in\Z}$ is well defined and it is
    strictly stationary and ergodic. Moreover, if $\E([\ln(Z_0^2)]^2) < \infty$ then it is also weakly stationary.
  \end{enumerate}
\end{thm}

Although both processes $\{X_t\}_{t \in \mathds{Z}}$ and $\{\sigma^2_t\}_{t\in\Z}$ are strictly stationary,
whenever $d<0.5$, they are not necessarily weakly stationary. This property depends on the distribution of
$Z_0$ and not only on the existence of the second moment for this random variable.  Theorem
\ref{thm:estacionar1} gives the condition for the existence of the $\r$-th moment, for any $\r>0$, for both
processes $\{X_t\}_{t \in \mathds{Z}}$ and $\{\sigma_t^2\}_{t\in\Z}$.

\begin{thm}\label{thm:estacionar1}
  Let $\{X_t\}_{t \in \mathds{Z}}$ be an \emph{SFIEGARCH}$(p,d,q)_s$ process, with $d<0.5$.  Assume that
  $\theta$ and $\gamma$ are not both equal to zero.  If there exists $\r > 2$ such that
  \begin{equation}\label{eq:convergencecriteria}
    \E(|Z_0|^{\r})<
    \infty \quad \mbox{ and} \quad
    \sum_{k = 0}^\infty \Big|\E\Big(\exp\Big\{\frac{\r}{2}\lambda_{d,k}g(Z_0)\Big\}\Big) -
    1\Big| < \infty,
  \end{equation}
  \noindent then, $\E(|X_t|^\m)<\infty$ and $\E(|\sigma_t|^\m)<\infty$, for all $t\in \Z$ and $ 0< \m\leq \r$.
\end{thm}

Assume that $d<0.5$ and $Z_0 \sim \mathcal{N}(0,1)$. Let $\Phi(\cdot)$ and $\erf(\cdot)$ be, respectively, the
standard Gaussian distribution and the error function, that is,
\[
\Phi(z) = \int_{-\infty}^z \frac{1}{\sqrt{2\pi}}e^{-\frac{t^2}{2}}\mathrm{d}t
\quad \mbox{and} \quad \erf(z) =
\frac{2}{\sqrt{\pi}}\int_{0}^z e^{-t^2} \mathrm{d}t = 2\Phi(z\sqrt{2}) - 1,
\quad \mbox{for all } z\in\R.
\]
It follows that,
\begin{align}
  \E\Big(\exp\Big\{\r\lambda_{d,k}g(Z_{0})\Big\}\Big)
  =&\int_{-\infty}^\infty  \exp\bigg\{\r\lambda_{d,k}\Big[\theta|z|+\gamma\Big(|z|-\sqrt{2/\pi}\Big)\Big]\bigg\}\frac{1}{\sqrt{2\pi}}e^{-\frac{z^2}{2}}\mathrm{d}z
  \nonumber\\
  = & \ \  e^{-\r\lambda_{d,k}\gamma\sqrt{2/\pi}}\Bigg[\exp\left\{\frac{(\r\lambda_{d,k})^2(\theta-\gamma)^2}{2}\right\}\Phi\big(\r\lambda_{d,k}(\gamma
  - \theta)\big)\nonumber\\
  &  \quad \hspace{49pt} \ +  \ \exp\left\{\frac{(\r\lambda_{d,k})^2(\theta+\gamma)^2}{2}\right\}\Phi\big(
  \r\lambda_{d,k}(\gamma+ \theta)\big)\Bigg]\nonumber\\
  =& \
  \frac{1}{2}e^{-\r\lambda_{d,k}\gamma\sqrt{2/\pi}}\Bigg(
  e^{\frac{(\r\lambda_{d,k})^2(\gamma-\theta)^2}{2}}\Bigg[1 + \erf\left(\frac{\r\lambda_{d,k}(\gamma-\theta)}{\sqrt{2}}\right)\Bigg]
  \nonumber\\
  &  \quad  \hspace{47pt} \ +  \, e^{\frac{(\r\lambda_{d,k})^2(\gamma+\theta)^2}{2}}\Bigg[1+\erf\left(\frac{\r\lambda_{d,k}(\gamma+\theta)}{\sqrt{2}}\right)\Bigg]
  \Bigg),\label{eq:gausscase}
\end{align}
\noindent for all $k\in\N$ and all $\r>0$. Since $e^x = 1+x+O(x^2)$ and $\erf(x) = \frac{2}{\sqrt{\pi}}x +
O(x^3)$, as $x\to 0$, one can rewrite \eqref{eq:gausscase} as,
\begin{align*}
  \E\Big(\!\exp\Big\{\r\lambda_{d,k}g(Z_{0})\Big\}\!\Big)
  = & \  \frac{1}{2}\left[1 -   \r\lambda_{d,k}\gamma\sqrt{\frac{2}{\pi}} +
    O(\lambda_{d,k}^2)\right]\!\!\left[2+2\r\lambda_{d,k}\gamma\sqrt{\frac{2}{\pi}}+
    O(\lambda_{d,k}^2) \right] =  1 + O(\lambda_{d,k}^2),
\end{align*}
as $k\to\infty$.  Thus, condition \eqref{eq:convergencecriteria} holds and hence, $\E(|X_t|^\r)<\infty$ and
$\E(|\sigma_t|^\r)<\infty$, for all $\r>0$. Corollary \ref{cor:estacionary} shows that this result also holds
if $Z_0\sim \mbox{GED}(\nu)$, for any $\nu>0$.

\begin{cor}\label{cor:estacionary}
  Let $\{X_t\}_{t \in \mathds{Z}}$ be an \emph{SFIEGARCH}$(p,d,q)_s$ process, with $d<0.5$.  Assume that
  $\theta$ and $\gamma$ are not both equal to zero.  Let $\{Z_t\}_{t \in \Z}$ be i.i.d.\! \emph{GED} with zero
  mean, variance equal to one, and tail-thickness parameter $v>1$.  Then, $\E(X_t^\r)<\infty$ and
  $\E([\sigma_t^2]^\r)<\infty$, for all $t\in \Z$ and $\r>0$.
\end{cor}

From expression (A.1.5) in \cite{N1991}, if $Z_0 \sim \mbox{GED}(\nu)$, with $\nu>1$, and $b_k :=
\frac{r}{2}\lambda_{d,k}$, for all $k\in\N$ and any $r>0$, then
\begin{equation}\label{eq:moments}
  \E(e^{b_kg(Z_0)}) =
  \exp\left\{\!-b_k\gamma \lambda  2^{1/\nu}\frac{\Gamma\big(\frac{2}{\nu}\big)}{\Gamma\big(\frac{1}{\nu}\big)}\!\right\}\sum_{j=0}^\infty\Big(\!b_k\lambda
  2^{1/\nu}\!\Big)^j\big[(\gamma+\theta)^j
  + (\gamma - \theta)^j\big]\frac{\Gamma(\frac{j+1}{\nu})}{2\Gamma(\frac{1}{\nu})\Gamma(j+1)},
\end{equation}
where $\lambda = \big[2^{1/\nu}\Gamma(1/\nu)/\Gamma(3/\nu)\big]^{1/2}$, for all $k\in\N$.  From expression
\eqref{eq:moments}, it is easy to see that $\E(X_t^r)$ is symmetric in $\theta$, for any $r>0$ and $\nu>1$.

\begin{figure}[!htbp]
  \vspace{-1\baselineskip}\centering
  \mbox{
    \subfloat[]{\includegraphics[width = 0.33\textwidth]{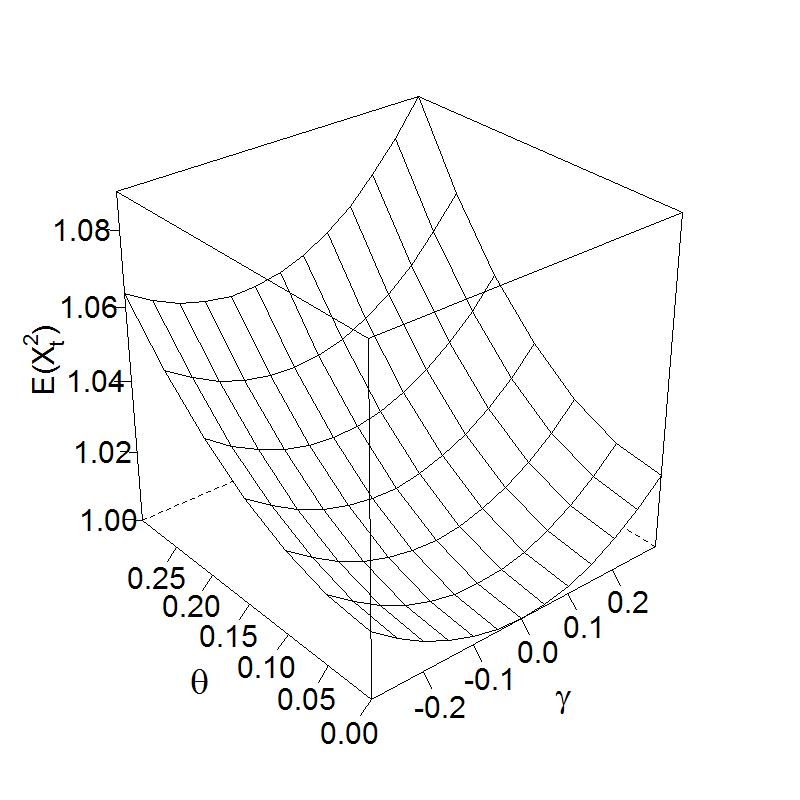}}
    \subfloat[]{\includegraphics[width = 0.33\textwidth]{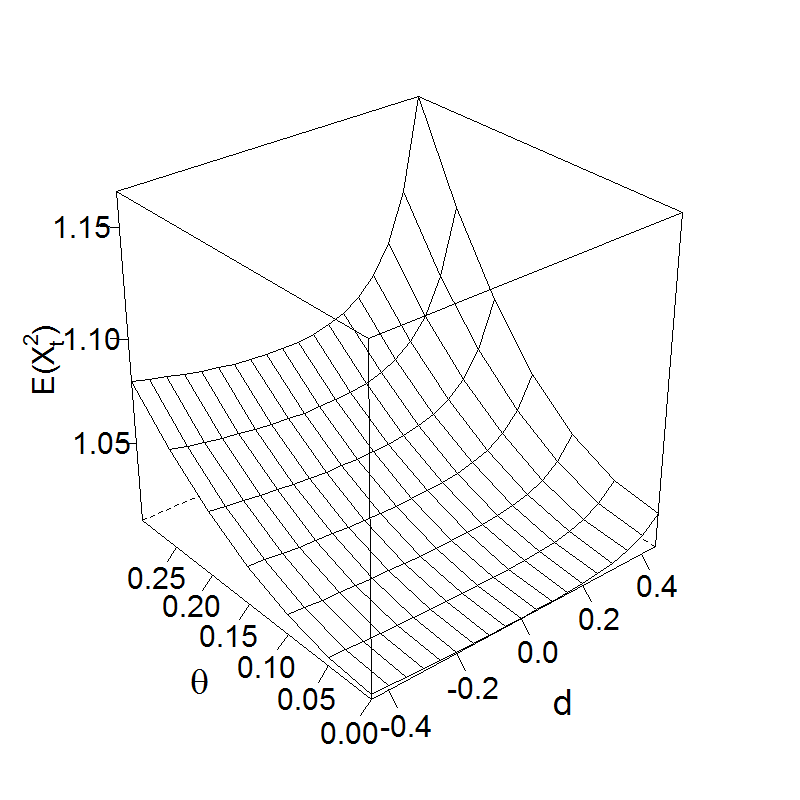}}
    \subfloat[]{\includegraphics[width = 0.33\textwidth]{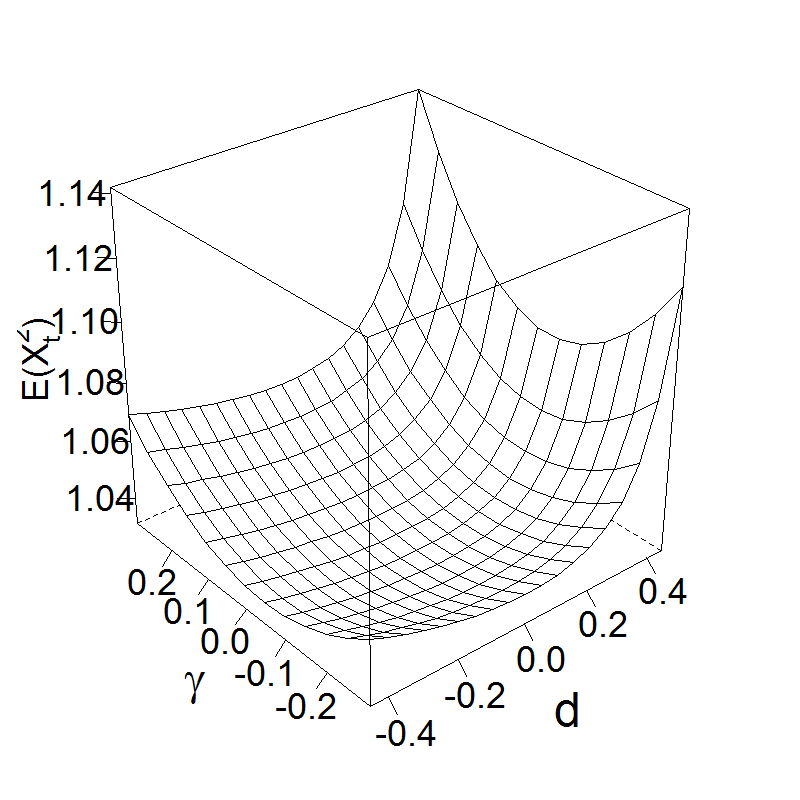}}
  }
  \caption{ This figure illustrates the behavior of $\E(X_t^2)$ with respect to the parameters $\theta,
    \gamma$ and $d$ when $\{X_t\}_{t\in\Z}$ is an SFIEGARCH$(0,d,0)_s$ processes, with $s=2$, $\omega = 0$ and
    $Z_0 \sim \mathcal{N}(0,1)$.  Panel (a) fixes $d = 0.25$ and shows $\E(X_t^2)$ as a function of $\theta$
    and $\gamma$. Panel (b) fixes $\gamma = 0.24$ and presents $\E(X_t^2)$ as a function of $\theta$ and
    $d$. Panel (c) fixes $\theta = 0.25$ and considers $\E(X_t^2)$ as a function of $\gamma$ and
    $d$. }\label{fig:moments}
\end{figure}

\begin{figure}[!htbp]
  \centering
  \subfloat[]{\includegraphics[width = 0.25\textwidth]{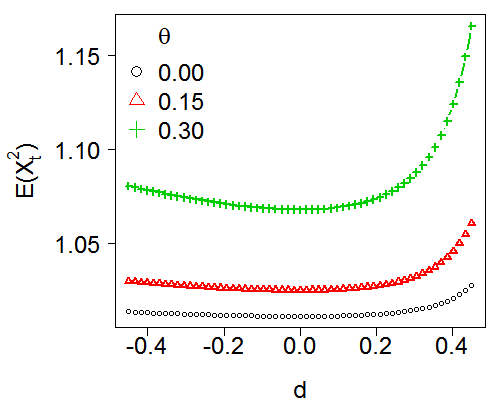}}
  \subfloat[]{\includegraphics[width = 0.25\textwidth]{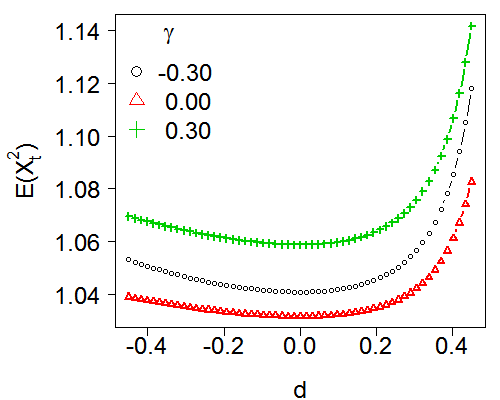}}
  \subfloat[]{\includegraphics[width = 0.25\textwidth]{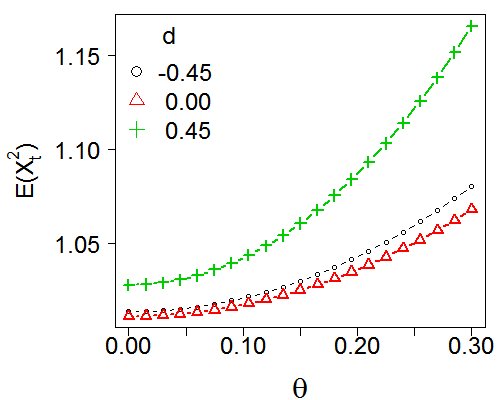}}
  \subfloat[]{\includegraphics[width = 0.25\textwidth]{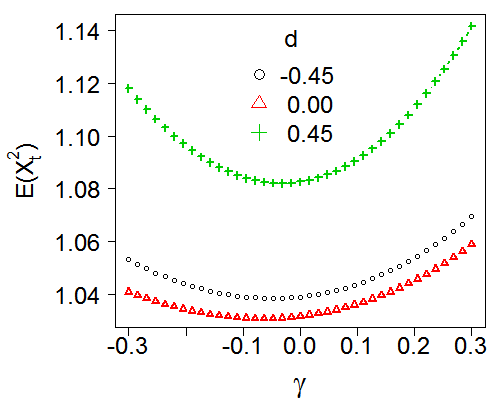}}
  \caption{ This figure shows the behavior of $\E(X_t^2)$ as a function of $d$, $\theta$ or $\gamma$, when
    $\{X_t\}_{t\in\Z}$ is an SFIEGARCH$(0,d,0)_s$ processes with $s=2$, $\omega = 0$ and
    $Z_0\sim\mathcal{N}(0,1)$.  Panel (a) shows $\E(X_t^2)$ as a function of $d$, for $\theta
    \in\{0,0.15,0.30\}$ and $\gamma = 0.24$. Panel  (b) considers $\E(X_t^2)$ as a function of $d$, for $\gamma
    \in\{-0.30,0.15,0.30\}$ and $\theta = 0.25$. Panel  (c) gives $\E(X_t^2)$ as a function of $\theta$, for $d
    \in\{-0.45,0,0.45\}$ and $\gamma = 0.24$. Panel  (d) presents $\E(X_t^2)$ as a function of $\gamma$, for $d
    \in\{-0.45,0,0.45\}$ and $\theta = 0.25$. }\label{fig:moments2}
\end{figure}

\begin{figure}[!htbp]
  \vspace{-2\baselineskip}\centering
  \mbox{
    \subfloat[]{\includegraphics[width = 0.3\textwidth]{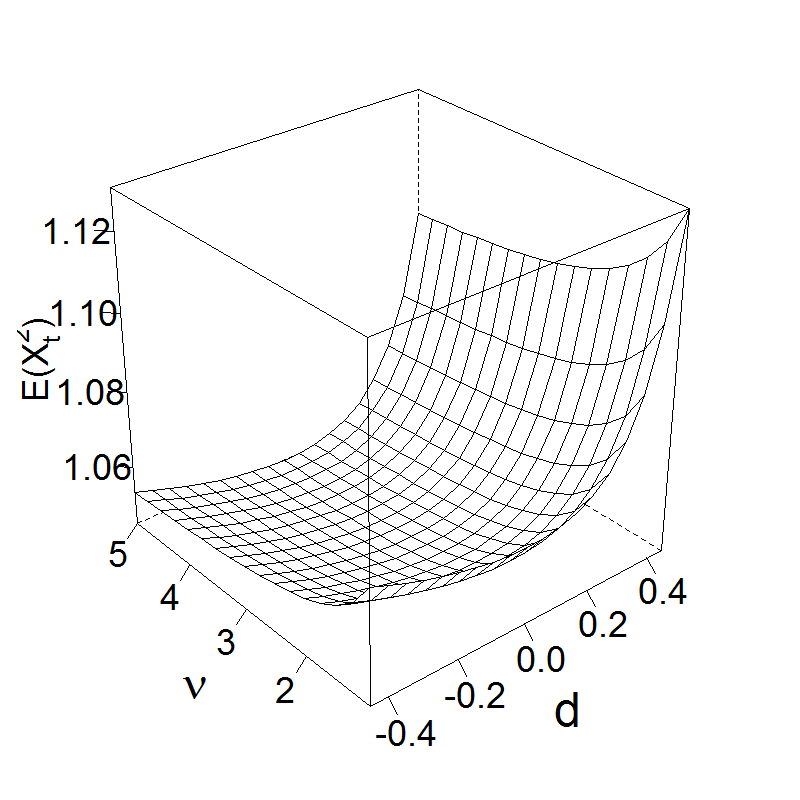}}
    \subfloat[]{\includegraphics[width = 0.33\textwidth]{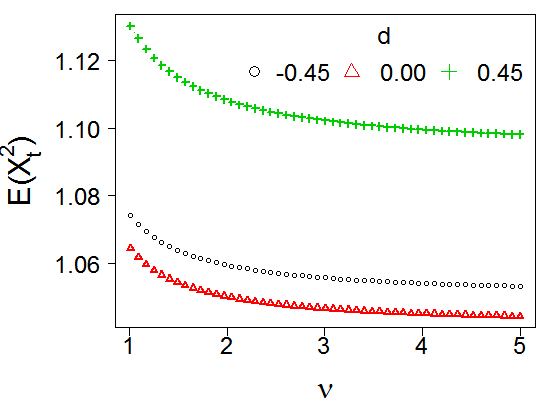}}
    \subfloat[]{\includegraphics[width = 0.33\textwidth]{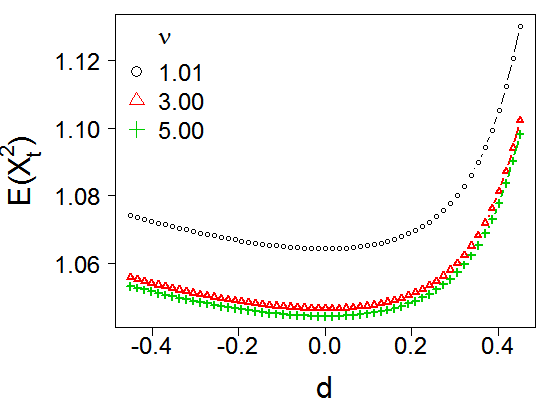}}
  }
  \caption{ This figure illustrates the behavior of $\E(X_t^2)$ as a function of $d$ and $\nu$, when
    $\{X_t\}_{t\in\Z}$ is an SFIEGARCH$(0,d,0)_s$ processes with $s=2$, $\theta=0.25$, $\gamma = 0.24$,
    $\omega = 0$ and $Z_0\sim\mbox{GED}(\nu)$. Panel (a) shows $\E(X_t^2)$ as a function of $\nu$ and $d$.
    Panel (b) gives $\E(X_t^2)$ as a function of $\nu$, for $d\in\{-0.45,0,0.45\}$. Panel  (c) considers $\E(X_t^2)$
    as a function of $d$, for $\nu\in\{1.01,3,5\}$. }\label{fig:moments3}
\end{figure}

\begin{example}
  Figures \ref{fig:moments} - \ref{fig:moments3} consider SFIEGARCH$(0,d,0)_s$ processes, with $s=2$. In these
  figures we analyze the behavior of $\E(X_t^2)$, with respect to the parameters $\theta, \gamma$ and $d$. We
  also study the behavior of $\E(X_t^2)$ with respect to the parameter $\nu$, when $Z_0 \sim \mbox{GED}(\nu)$.
  From expression \eqref{eq:moments} one observes that $\E(X_t^2)$ is symmetric in $\theta$, whenever $Z_0\sim
  \mbox{GED}(\nu)$, for any $\nu >1$. Therefore, for these figures we only consider positive values of
  $\theta$.  Figure \ref{fig:moments} (a) shows the behavior of $\E(X_t^2)$ as a function of $\theta$ and
  $\gamma$, for $d=0.25$. Figure \ref{fig:moments} (b) presents $\E(X_t^2)$ as a function of $\theta$ and $d$,
  for $\gamma = 0.24$.  Figure \ref{fig:moments} (c) shows $\E(X_t^2)$ as a function of $\gamma$ and $d$, for
  $\theta = 0.25$.  For all graphs in Figure \ref{fig:moments}, $s=2$, $\omega = 0$ and
  $Z_0\sim\mathcal{N}(0,1)$.

  From Figure \ref{fig:moments}, one observes that for $\theta$ or $\gamma$ fixed, $\E(X_t^2)$ slowly
  decreases for $d\in [-0.45,0]$ and increases for $d\in [0, 0.45]$.  This behavior can be better observed in
  Figures \ref{fig:moments2} (a) and (b) where the values of $\E(X_t^2)$, as a function of $d$, are plotted for
  $\theta \in\{0,0.15,0.30\}$ and $\gamma \in\{-0.30,0,0.30\}$, respectively.  Similarly, for each $d$ fixed,
  the function $\E(X_t^2)$ is decreasing for $\theta,\gamma \in [-0.3,0.0]$ (the function is symmetric in
  $\theta$) and it is increasing for $\theta, \gamma \in [0.0,0.3]$.  This behavior can be observed in Figures
  \ref{fig:moments2} (c) and (d) where $\E(X_t^2)$ is given, respectively, as a function of $\theta$ and
  $\gamma$, for $d \in\{-0.40,0,0.45\}$.
\end{example}

\begin{example}
  Figure \ref{fig:moments3} (a) shows the graph of $\E(X_t^2)$, as a function of $\nu$ and $d$, when $Z_0 \sim
  \mbox{GED}(\nu)$, with $\nu > 1$.  Figures \ref{fig:moments3} (b) and (c) present the graph of $\E(X_t^2)$,
  respectively, as a function of $\nu$, for $d\in\{-0.45,0,0.45\}$ and as a function of $d$, for
  $\nu\in\{1.01,3,5\}$.  For all graphs, $s = 2$, $\omega = 0$, $\theta = 0.25$ and $\gamma = 0.24$.  From
  Figure \ref{fig:moments3} one concludes that $\E(X_t^2)$ is a decreasing function of $\nu$ and, as a
  function of $d$, $\E(X_t^2)$ presents the same behavior as in the Gaussian case, that is, it is decreasing
  for $d\in[-0.45,0]$ and increasing for $d\in[0,0.45]$.
\end{example}

The following proposition presents the kurtosis and the asymmetry measures for any stationary SFIEGARCH
process.

\begin{prop}\label{prop:skewnessKurtosys}
  Let $\{X_t\}_{t \in \mathds{Z}}$ be a stationary \emph{SFIEGARCH}$(p,d,q)_s$ process with $\E(|X_t^4|)<
  \infty$.  The asymmetry and kurtosis measures of $\{X_t\}_{t \in \mathds{Z}}$ are given, respectively, by
  \begin{equation*}
    A_X=\E(Z_0^3)\, \frac{\displaystyle\prod_{k=0}^{\infty}\E\bigg(\exp\left\{
        \frac{3}{2}\lambda_{d,k}g(Z_{t})\right\}\bigg)}{\displaystyle\prod_{k=0}^{\infty}\Big[\E\Big(\exp\left\{
        \lambda_{d,k}g(Z_{t})\right\}\Big)\Big]^{3/2}}
    \quad \mbox{\normalsize and } \quad
    K_X=\E(Z_0^4)\, \frac{\displaystyle\prod_{k=0}^{\infty}\E\Big(\exp\left\{
        2\lambda_{d,k}g(Z_{t})\right\}\Big)}{\displaystyle\prod_{k=0}^{\infty}\Big[\E\Big(\exp\left\{
        \lambda_{d,k}g(Z_{t})\right\}\Big)\Big]^2}.
  \end{equation*}
\end{prop}

\begin{figure}[!htbp]
  \vspace{-1\baselineskip}\centering
  \subfloat[]{\includegraphics[width = 0.33\textwidth]{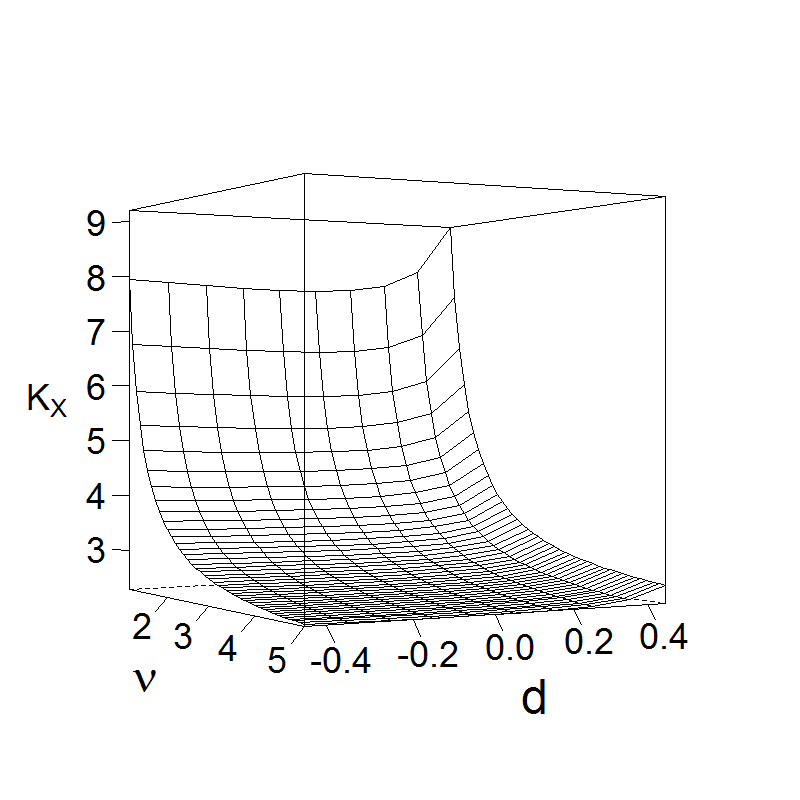}}
  \subfloat[]{\includegraphics[width = 0.38\textwidth]{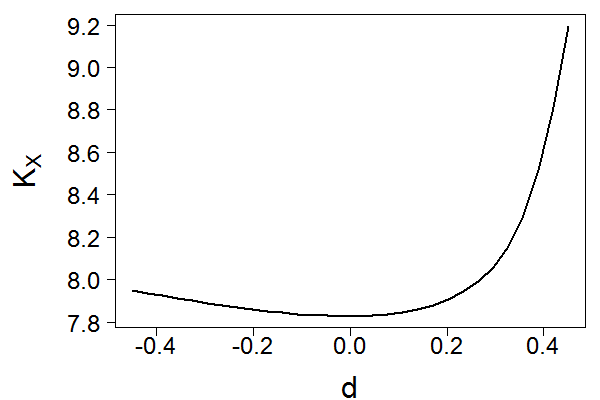}}
  \caption{ The kurtosis measure of an SFIEGARCH$(0,d,0)_s$ process with $s=2$, $\omega = 0$, $\theta = 0.25$,
    $\gamma = 0.24$ and $Z_0 \sim \mbox{GED}(\nu)$.  Panel (a) shows the graph of the kurtosis measure as a
    function of $d$ and $\nu$. Panel (b) presents the kurtosis measure as a function of $d$, for $\nu = 2$
    fixed. }\label{fig:kurtosisfig}
  \vspace{-1\baselineskip}
\end{figure}

\begin{example}
  From expression \eqref{eq:moments}, one easily concludes that the kurtosis measure $K_X$ is symmetric in
  $\theta$, whenever $Z_0 \sim \mbox{GED}(\nu)$, for any $\nu>1$.  Figure \ref{fig:kurtosisfig} (a) considers
  SFIEGARCH$(0,d,0)_s$ processes and shows the behavior of $ K_X$ as a function of $d$ and $\nu$, for $s=2$,
  $\omega = 0$, $\theta = 0.25$ and $\gamma = 0.24$. Figure \ref{fig:kurtosisfig} (b) presents the graph of $
  K_X$ as a function of $d$, for $\nu = 2$.  From the graphs in Figure \ref{fig:kurtosisfig} one concludes
  that the kurtosis measure is a decreasing function of $\nu$, for each $d$ fixed. Moreover, for each $\nu$
  fixed, it is decreasing for $d\in[-0.45,0]$ and increasing for $d\in[0, 0.45]$.
\end{example}

\begin{figure}[!htb]
  \vspace{-2\baselineskip} \centering
  \subfloat[]{\includegraphics[width = 0.33\textwidth]{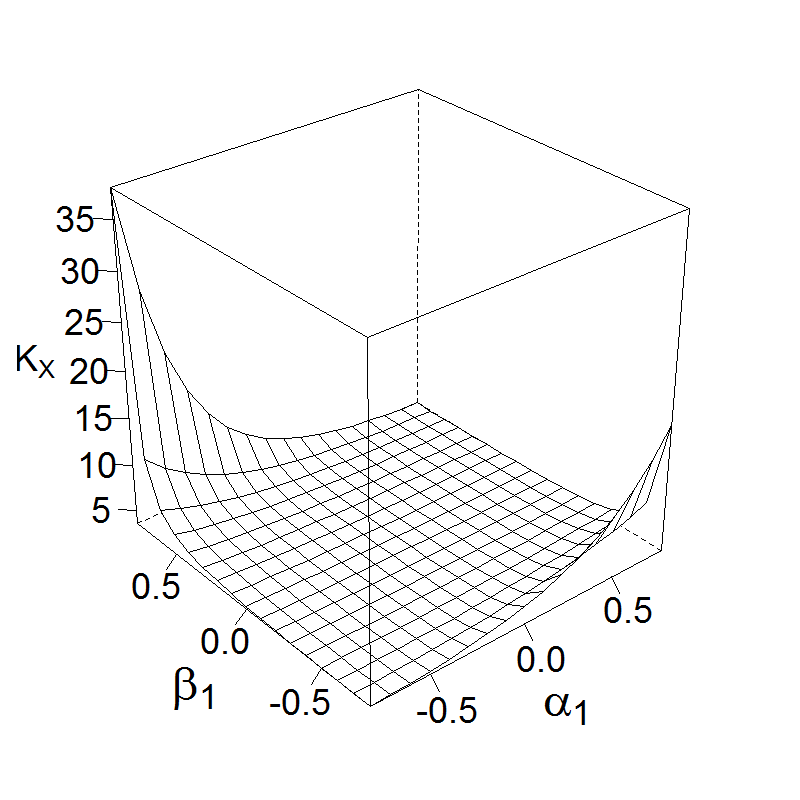}}
  \subfloat[]{\includegraphics[width = 0.33\textwidth]{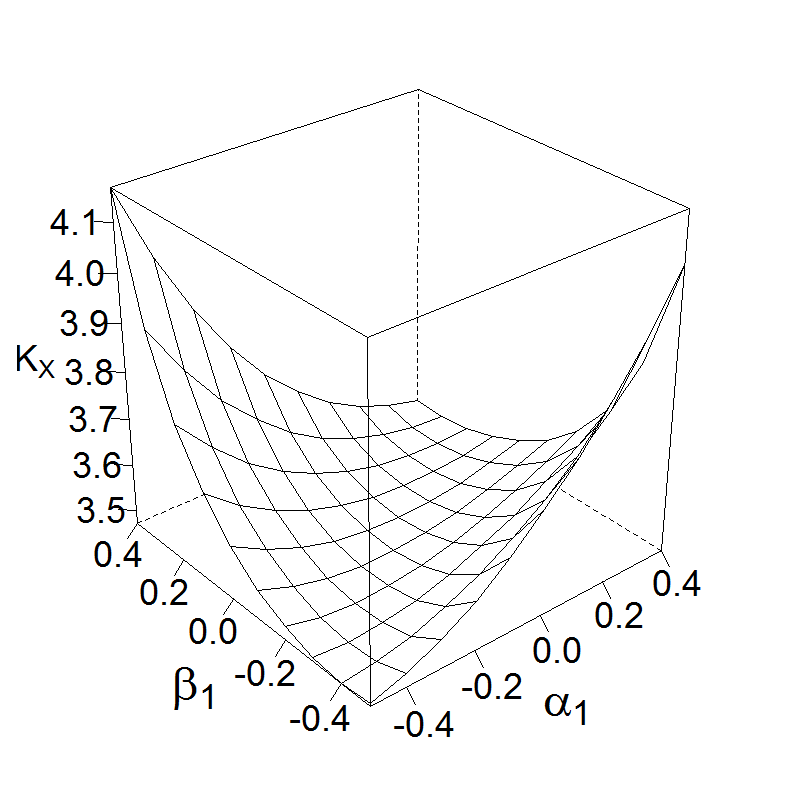}}
  \caption{ The kurtosis measure of an SFIEGARCH$(p,d,q)_s$ process with $p,q\in\{0,1\}$, $s=2$, $\omega = 0$,
    $\theta = 0.25$, $\gamma = 0.24$ and $Z_0 \sim \mbox{GED}(\nu)$.  Panel (a) shows the kurtosis measure as
    a function of $\alpha_1$ and $\beta_1$ for $\alpha_1,\beta_1 \in [-0.8,0.8]$. Panel (b) considers $\alpha_1$ and
    $\beta_1$ only in the interval $[-0.4,0.4]$ (for a better visualization).  }\label{fig:kurtosisfig2}
\end{figure}

\begin{example}
  Figure \ref{fig:kurtosisfig2} considers SFIEGARCH$(p,d,q)_s$ processes, with $p,q\in\{0,1\}$, $d = 0.25$ and
  the same values of $s,\omega, \theta$ and $\gamma$ as in Figure \ref{fig:kurtosisfig}.  This figure presents
  the behavior of $ K_X$ as a function of $\alpha_1$ and $\beta_1$.  The cases $\alpha_1 = \beta_1$ (the
  polynomials have a common root) are actually equivalent to the case $\alpha_1 = 0 =\beta_1$ and, in this
  case, one has an SFIEGARCH$(0,d,0)_s$ process.  While Figure \ref{fig:kurtosisfig2} (a) shows the graphs of
  $K_X$ for $\alpha_1,\beta_1 \in [-0.8,0.8]$, Figure \ref{fig:kurtosisfig2} (b) considers only the interval
  $[-0.4,0.4]$.
\end{example}

From Figure \ref{fig:kurtosisfig2} one observes that the behavior of $K_X$ depends on the sign of both
$\alpha_1$ and $\beta_1$.  Also, it goes from increasing (when $\beta_1 = -0.8$) to decreasing (when $\beta_1
= 0.8$) in $\alpha_1$. Also, by comparing the graphs in Figures \ref{fig:kurtosisfig2} (a) and (b), it is easy
to see that the kurtosis measure presents small variation on its value for $(\alpha_1,\beta_1) \in
[-0.4,0.4]\times [-0.4,0.4]$ (in this region $3 < K_X < 4.5 $). Moreover, for $\alpha_1 \in[-0.8,-0.4]$ and
$\beta_1\in[-0.8,0.8]$ the kurtosis measure varies faster than in the case $\beta_1 \in[-0.8,-0.4]$ and
$\alpha_1\in[-0.8,0.8]$.

Lemma \ref{lemma:expressioncovln} and Corollary \ref{cor:rhoconv} present the autocovariance function of the
process $\{\ln(\sigma_t^2)\}_{t\in\Z}$ and its asymptotic behavior. Although, in practice, this stochastic
process cannot be observed and hence, the sample autocovariance structure cannot be analysed, the results
presented in these theorems are necessary to prove Theorem \ref{thm:resultsXt}.

\begin{lemma}\label{lemma:expressioncovln}
  Let $\{X_t\}_{t \in \mathds{Z}}$ be an \emph{SFIEGARCH}$(p,d,q)_s$ process, defined by \eqref{eq:xt} and
  \eqref{eq:sigmat}, with $ d < 0.5$.  Suppose that $\gamma$ and $\theta$ are not both equal to zero. Then,
  the autocovariance function $\gamma_{\ln(\sigma_t^2)}(\cdot)$ of the process $\{\ln(\sigma_t^2)\}_{t\in\Z}$
  is given by
  \begin{equation}
    \gamma_{\ln(\sigma_t^2)}(sh + r) =
    \sum_{k \in\Z} \gamma_{\mbox{\tiny$A$}}(sk + r)\gamma_{\mbox{\tiny$V$}}(sh - sk),
    \quad \mbox {\normalsize  for all } h\in\Z \quad \mbox{\normalsize  and }r\in\{0,\cdots, s-1\}, \label{eq:covls}
  \end{equation}
  \noindent where   $\gamma_{\mbox{\tiny $A$}}(\cdot)$   and
  $\gamma_{\mbox{\tiny $V$}}(\cdot)$ are given, respectively,  by
  \begin{equation}\label{eq:covarma}
    \gamma_{\mbox{\tiny$A$}}(h) = \sum_{i = 0}^{\infty}f_if_{i+|h|}, \quad
    \mbox{for all } h\in\Z, \quad \mbox{with } f(z) := \frac{\alpha(z)}{\beta(z)} =\sum_{k =0}^{\infty}
    f_kz^k,
  \end{equation}
  and
  \vspace{-0.5\baselineskip}
  \begin{equation}
    \gamma_{\mbox{\tiny $V$}}(hs) =  \sigma_g^2
    \frac{(-1)^{h}\Gamma(1-2d)}{\Gamma(1-d+h)\Gamma(1-d-h)}
    \quad \mbox{ and } \quad
    \gamma_{\mbox{\tiny $V$}}(hs + r) = 0, \label{eq:covln}
  \end{equation}
  for all  $h\in\Z$ and $r\in\{0,\cdots, s-1\}$, with $\sigma_g^2$ given by \eqref{eq:sigma2gzt}.
\end{lemma}

From Lemma \ref{lemma:expressioncovln} one easily concludes that, if $p=0=q$, then
$\gamma_{\ln(\sigma_t^2)}(sh + r)\neq 0$ if an only if $r = 0$.  Moreover, if $p>0$ and $q=0$,
$\gamma_{\ln(\sigma_t^2)}(sh + r) = 0$, for all $h\in\Z$, whenever $\left\lceil-\frac{p+r}{s}\right\rceil >
\left\lfloor\frac{p-r}{s}\right\rfloor$. This is so because $\gamma_{\mbox{\tiny$A$}}(sk+r) \neq 0$ if and
only if $|sk+r| \leq p$, that is, $-\frac{p+r}{s} \leq k \leq \frac{p-r}{s}$ and $k\in\Z$. Thus, for $p\geq 0$
and $q = 0$, one can rewrite \eqref{eq:covls} as
\[
\gamma_{\ln(\sigma_t^2)}(sh + r) = \left\{
  \begin{array}{cc}
    \displaystyle\sum_{\left\lceil-\frac{p+r}{s}\right\rceil \leq k \leq \left\lfloor\frac{p-r}{s}\right\rfloor}  \gamma_{\mbox{\tiny$A$}}(sk + r)\gamma_{\mbox{\tiny$V$}}(sh
    -sk), &  \mbox{if} \quad \left\lceil-\frac{p+r}{s}\right\rceil \leq \left\lfloor\frac{p-r}{s}\right\rfloor;\\
    0,  & \mbox{otherwise},
  \end{array}
\right.
\]
\noindent for all $r\in\{0,\cdots,s-1\}$, and $h\in\Z$.  In this case, it is obvious that
$\sum_{h=0}^\infty|\rho_{\ln(\sigma_t^2)}(h)| < \infty$ if and only if,
$\sum_{h=0}^\infty|\rho_{\mbox{\tiny$V$}}(h)| < \infty$.  Corollary \ref{cor:rhoconv} presents the asymptotic
behavior of $\gamma_{\ln(\sigma_t^2)}(h)$, as $h\to\infty$, which leads to the conclusion that this result
actually holds for any $p$ and $q$.

\begin{cor}\label{cor:rhoconv}
  Let $\{X_t\}_{t \in \mathds{Z}}$ be an \emph{SFIEGARCH}$(p,d,q)_s$ process, defined by \eqref{eq:xt} and
  \eqref{eq:sigmat}, with $d < 0.5$.  Suppose that $\gamma$ and $\theta$ are not both equal to zero.  Let
  $\gamma_{\ln(\sigma_t^2)}(\cdot)$ be the autocovariance function of the process
  $\{\ln(\sigma_t^2)\}_{t\in\Z}$. Then, for all $r\in\{0,\cdots,s-1\}$,
  \begin{equation}\label{eq:covlsinfty}
    \gamma_{\ln(\sigma_t^2)}(sh + r) =
    \gamma_{\mbox{\tiny$V$}}(sh)\mathscr{G}(sh+r)  + o(h^{-\nu}), \quad
    \mbox{as} \quad h\to \infty,
  \end{equation}
  \noindent where $\nu$ is any positive number and $\mathscr{G}(\cdot)$ is a real function satisfying
  \[
  \lim_{h\to \infty} \sum_{r=0}^{s-1}\mathscr{G}(sh+r) = \sum_{k\in\Z
  }\gamma_{\mbox{\tiny$A$}}(k),\quad \mbox{with $\gamma_{\mbox{\tiny$A$}}(\cdot)$ given by \eqref{eq:covarma}.}
  \]
\end{cor}

Theorem \ref{thm:resultsXt} presents the autocovariance function $ \gamma_{\ln(\X_t^2)}(\cdot)$ of the process
$\{\ln(X_t^2)\}_{t\in\Z}$, where $\{X_t\}_{t \in \mathds{Z}}$ is an SFIEGARCH process. This theorem also gives
the asymptotic behavior of $ \gamma_{\ln(\X_t^2)}(sh+r)$, for all $r\in\{0,\cdots,s-1\}$, as $h$ goes to
infinity.

\begin{thm}\label{thm:resultsXt}
  Let $\{X_t\}_{t \in \mathds{Z}}$ be an \emph{SFIEGARCH}$(p,d,q)_s$ process, defined by \eqref{eq:xt} and
  \eqref{eq:sigmat}, with $ d < 0.5$. Suppose that $\gamma$ and $\theta$ are not both equal to zero and $
  \var(\ln(Z_0^2)) := \sigma^2_{\ln(Z_t^2)} < \infty$. Then, the autocovariance function of the process
  $\{\ln(X_t^2)\}_{t\in\Z}$ is given by
  \begin{equation}\label{eq:covlx}
    \gamma_{\ln(\X_t^2)}(h)  = \gamma_{\ln(\sigma_t^2)}(h) +C_1\lambda_{d,|h|-1}\Ind{\Z^*}{h} +\sigma^2_{\ln(Z_t^2)}\Ind{\{0\}}{h},\quad \mbox{for all } h\in\Z,
  \end{equation}
  \noindent where $C_1=\cov(g(Z_0),\ln(Z_0^2))$ and $\gamma_{\ln(\sigma_t^2)}(\cdot)$ is given in Lemma
  \ref{lemma:expressioncovln}. Thus, for all $\nu>0$,
  \begin{equation}\label{eq:gammainfty}
    \gamma_{\ln(\X_t^2)}(sh+r) = \gamma_{\mbox{\tiny$V$}}(sh)\mathscr{G}(sh+r) +
    \pi_{d,s\lfloor\frac{sh+r-1}{s}\rfloor}\mathscr{K}(sh+r-1) + o(h^{-\nu}), \quad
    \mbox{as } h\to\infty,
  \end{equation}
  for all $r\in \{0,\cdots,s-1\}$, where $\mathscr{G}(\cdot)$ and $\mathscr{K}(\cdot)$ are given, respectively,
  in Corollary \eqref{cor:rhoconv} and Theorem \ref{thm:convergenceOrder}.
\end{thm}

\begin{figure}[!htbp]
  \vspace{-1\baselineskip}  \centering
  \mbox{
    \subfloat[]{\includegraphics[width = 0.33\textwidth]{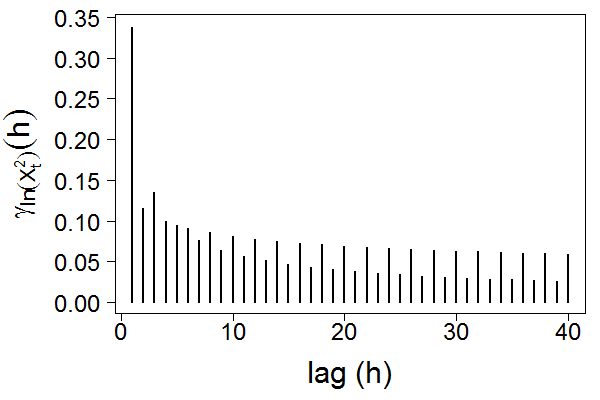}}
    \subfloat[]{\includegraphics[width = 0.33\textwidth]{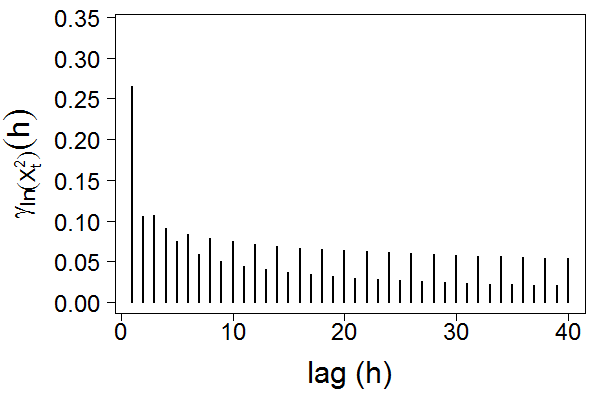}}
    \subfloat[]{\includegraphics[width = 0.33\textwidth]{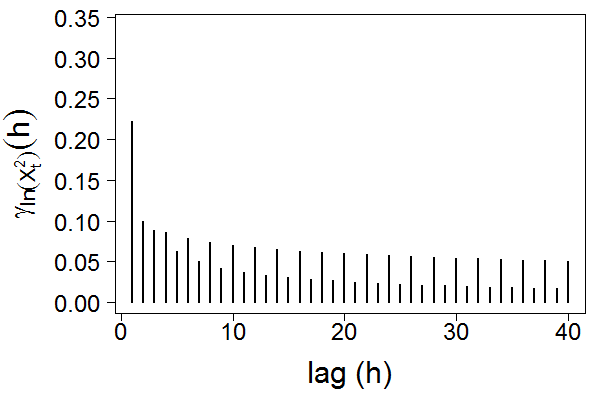}}
  }
  \caption{ Theoretical autocovariance function $\gamma_{\ln(\X_t^2)}(\cdot)$, for $h\in\{1,\cdots,100\}$,
    corresponding to the process $\{\ln(X_t^2)\}_{t\in\Z}$, when $\{X_t\}_{t \in \mathds{Z}}$ is an
    SFIEGARCH$(0,d,0)_s$ process with $d = 0.4$, $s = 2$, $\omega = 0$, $\theta = 0.25$, $\gamma = 0.24$ and
    $Z_0 \sim \mbox{GED}(\nu)$.  Panel (a) considers $\nu = 1.01$. Panel (b) assumes $\nu = 2$ (Gaussian
    case). Panel  (c) fixes $\nu = 5$. }\label{fig:autocov}
\end{figure}

\begin{example}
  From expressions \eqref{eq:covls} and \eqref{eq:covlx}, one concludes that, if $Z_0$ is a symmetric random
  variable, then $\gamma_{\ln(\X_t^2)}(\cdot)$ is symmetric in $\theta$. Figures \ref{fig:autocov} (a) - (d)
  show the graphs of $ \gamma_{\ln(\X_t^2)}(h)$, for $h\in\{1,\cdots,100\}$, where $\{X_t\}_{t \in
    \mathds{Z}}$ is an SFIEGARCH$(0,d,0)_s$ process, with $d = 0.4$, $s = 2$, $\omega = 0$, $\theta = 0.25$,
  $\gamma = 0.24$ and $Z_0 \sim \mbox{GED}(\nu)$, for $\nu\in\{1.01,2,3,5\}$, respectively. All graphs are
  presented in the same scale for a better visualization. The corresponding values of
  $\gamma_{\ln(\X_t^2)}(0)$ are, respectively, 6.7228, 5.0978, 4.6445 and 4.3556.  From Figure
  \ref{fig:autocov}, one observes that, for each fixed $h$, $\gamma_{\ln(\X_t^2)}(h)$ decreases as $\nu$
  increases.
\end{example}

\begin{figure}[!h]
  \vspace{-1.5\baselineskip}\centering
  \subfloat[]{\includegraphics[width = 0.33\textwidth]{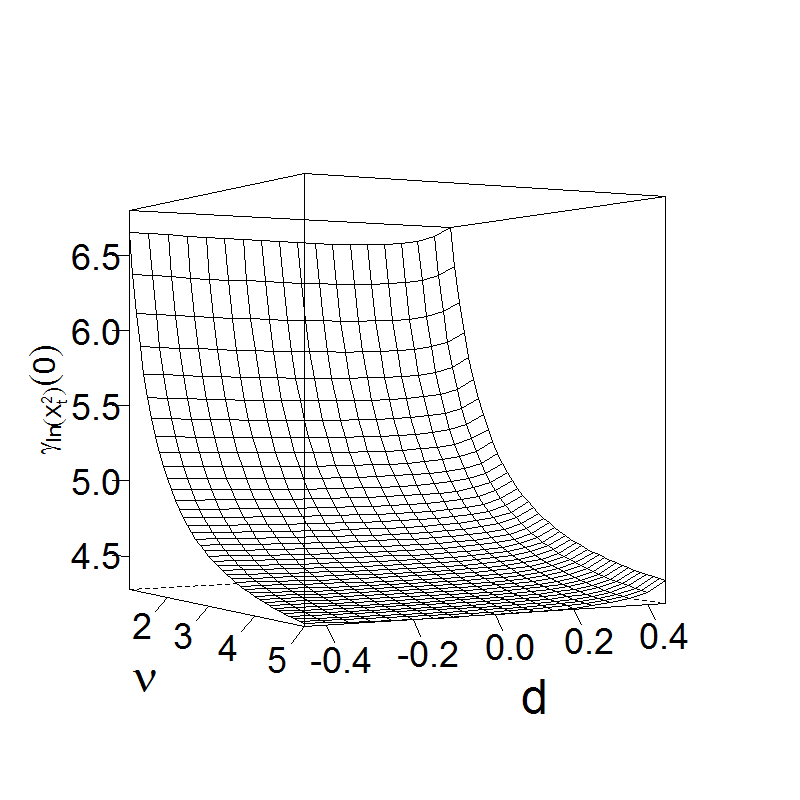}}
  \hspace{0.5cm}
  \subfloat[]{\includegraphics[width = 0.38\textwidth]{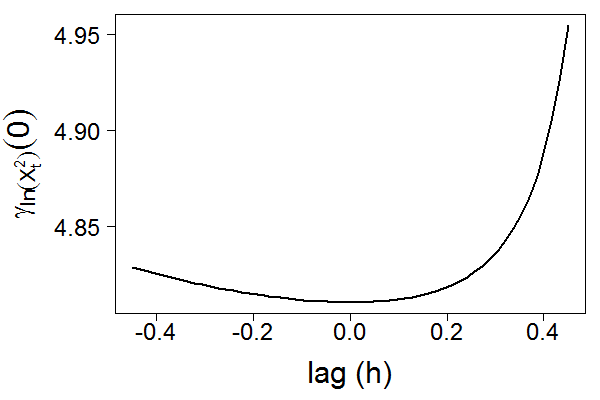}}
  \caption{ This figure shows the graphs of $\sigma_{\ln(\X_t^2)}^2 :=\gamma_{\ln(\X_t^2)}(0)$, the variance of
    the process $\{\ln(X_t^2)\}_{t\in\Z}$, when $\{X_t\}_{t \in \mathds{Z}}$ is an SFIEGARCH$(0,d,0)_s$ process
    with $s = 2$, $\omega = 0$, $\theta = 0.25$, $\gamma = 0.24$ and $Z_0 \sim \mbox{GED}(\nu)$.  Panel (a)
    considers the variance as a function of $d$ and $\nu$. Panel (b) shows the variance as a function of $d$, for
    $\nu = 2$ fixed.  }\label{fig:autocov2}
\end{figure}

\begin{example}
  Figure \ref{fig:autocov2} (a) presents $\gamma_{\ln(\X_t^2)}(0) = \var(\ln(X_t^2))$ as a function of $\nu$
  and $d$, where $\{X_t\}_{t \in \mathds{Z}}$ is an SFIEGARCH$(0,d,0)_s$ process, with $s = 2$, $\omega = 0$,
  $\theta = 0.25$ and $\gamma = 0.24$. Figure \ref{fig:autocov2} (b) presents $\gamma_{\ln(\X_t^2)}(0)$ as a
  function of $d$, for $\nu = 2$. From the graphs in Figure \ref{fig:autocov2} one observes that the variance
  of $\{\ln(X_t^2)\}_{t\in\Z}$ decreases with $\nu$. For each $\nu$ fixed, $\gamma_{\ln(\X_t^2)}(0)$ is
  decreasing for $d\in[-0.45,0]$ and increasing for $d\in[0,0.45]$.
\end{example}

The following corollary compares the asymptotic behavior of $\sum_{r=0}^{s-1}\gamma_{\ln(\X_t^2)}(sh+r)$ when
$d<0$ and $d>0$.

\begin{cor}\label{cor:asymptcov}
  Let $\gamma_{\ln(\X_t^2)}(\cdot)$ be the autocovariance function of the process $\{\ln(X_t^2)\}_{t\in\Z}$,
  given in Theorem \ref{thm:resultsXt}. Then,
  \begin{equation}\label{eq:asymptcov}
    \sum_{r=0}^{s-1}\gamma_{\ln(\X_t^2)}(sh+r)  \sim \left\{
      \begin{array}{cc}
        \mathscr{C}_1(h) h^{d-1}, &  \mbox{ if } d<0;\vspace{0.3cm}\\
        \mathscr{C}_2(h)h^{2d-1},  & \mbox{ if } d > 0,
      \end{array}
    \right.
  \end{equation}
  with
  \vspace{-0.7\baselineskip}
  \begin{equation}\label{eq:c1c2}
    \mathscr{C}_1(h) := \sigma^2_g
    \frac{\Gamma(1-2d)}{\Gamma(1-d)\Gamma(d)}\sum_{r=0}^{s-1}\mathscr{G}(sh+r)
    \quad \mbox{and} \quad \mathscr{C}_2(h) := \frac{1}{\Gamma(d)}\sum_{r=0}^{s-1}\mathscr{K}(sh+r-1),
  \end{equation}
  where $\mathscr{G}(\cdot)$ and $\mathscr{K}(\cdot)$ are given, respectively, in Corollary \ref{cor:rhoconv}
  and Theorem \ref{thm:convergenceOrder} and $\sigma_g^2$ is given by \eqref{eq:sigma2gzt}.
\end{cor}

\section{Spectral Representation}\label{sec:spectralrep}

Recently economists have noticed that volatility of high frequency financial time series shows long range
dependence merged with periodic behavior due to some operating features of financial markets.  Periodic
components are represented as marked peaks at some frequencies in the periodogram function.  It is a well
known result that the periodogram function is an estimator of the spectral density function. Therefore, in
order to choose the best model for a time series, one should know how does the spectral density function
behaves in order to gather information from the periodogram function.

Here we present the spectral density function of both processes $\{\ln(\sigma_t^2)\}_{t\in\Z}$ and
$\{\ln(X_t^2)\}_{t\in\Z}$.  It is easy to see that expression \eqref{eq:spectraldensity}, in Theorem
\ref{thm:resultsXt2}, is similar to the expression (2.5) from \cite{HEA2005}. In this paper, the authors
present the asymptotic properties of some semiparametric estimators for the long-memory parameter for a class
of stochastic process which includes LMSV (\!\emph{Long Memory Stochastic Volatility}) and FIEGARCH models.

\begin{thm}\label{thm:resultsXt2}
  Let $\{X_t\}_{t \in \mathds{Z}}$ be an \emph{SFIEGARCH}$(p,d,q)_s$ process, defined by \eqref{eq:xt} and
  \eqref{eq:sigmat}, with $ d < 0.5$.  Suppose that $\gamma$ and $\theta$, given in \eqref{eq:gzt}, are not
  both equal to zero and that $\alpha(z)\neq 0$ in the closed disk $\{z:|z|\leq 1\}$. If $\var(\ln(Z_t^2)) :=
  \sigma^2_{\ln(Z_t^2)} < \infty$, for all $t\in\Z$, the spectral density function of
  $\{\ln(X_t^2)\}_{t\in\Z}$ is given by
  \begin{equation}
    f_{\ln(\X_t^2)}(\lambda) = f_{\ln(\sigma_t^2)}(\lambda) +  \frac{C_1}{\pi}\Re\big(e^{-\im \lambda}\Lambda(\lambda)\big) + f_{\ln(Z_t^2)}(\lambda), \quad \mbox{for all }  \lambda \in [0,\pi], \label{eq:spectraldensity}
  \end{equation}
  where $f_{\ln(\sigma_t^2)}(\cdot)$ is given in \eqref{eq:fln}, $f_{\ln(Z_t^2)}(\lambda) =
  \frac{\sigma^2_{\ln(Z_t^2)}}{2\pi}$ is the spectral density function of $\{\ln(Z_t^2)\}_{t\in\Z}$, $C_1 =
  \cov\big(g(Z_0),\ln(Z_0^2)\big)$, $\Lambda(z) := \lambda(e^{-\im z})$ and $\lambda(\cdot)$ is defined in
  \eqref{eq:polilambda}.
\end{thm}

Notice that the spectral density function is symmetric around $\pi$. Hence, in what follows, although the
graphs consider the interval $[0,2\pi]$, one only needs to pay attention to the interval $[0,\pi]$.  Moreover,
all graphs are presented in the same scale and they are truncated in the $y$-axis for a better visualization.

\begin{figure}[!htbp]
  \centering
  \mbox{
    \hspace{-2pt}\subfloat[]{\includegraphics[width = 0.25\textwidth]{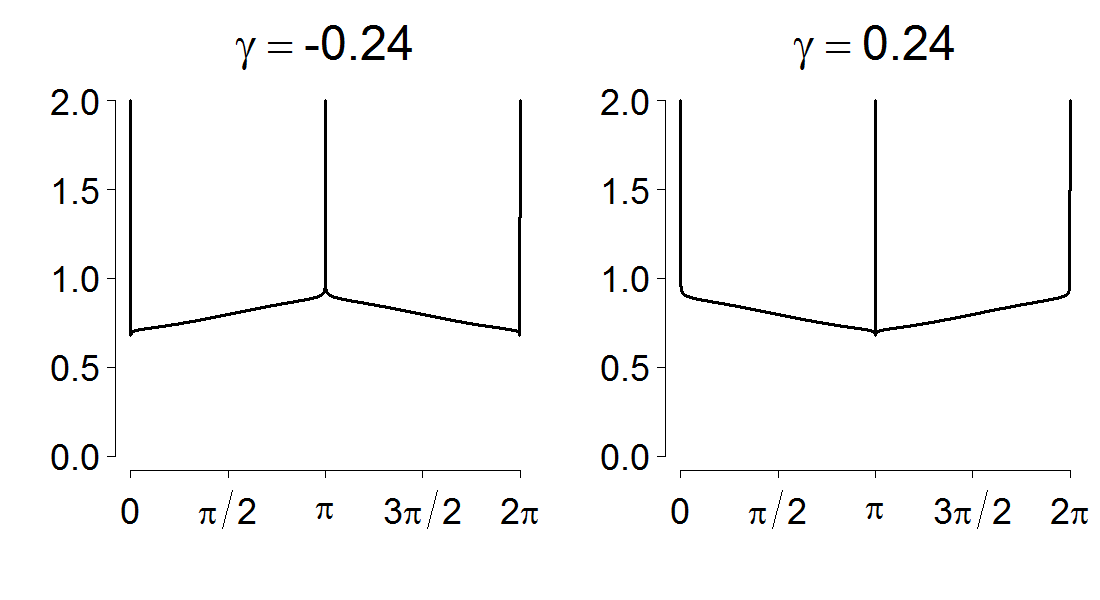}}
    \subfloat[]{\includegraphics[width = 0.25\textwidth]{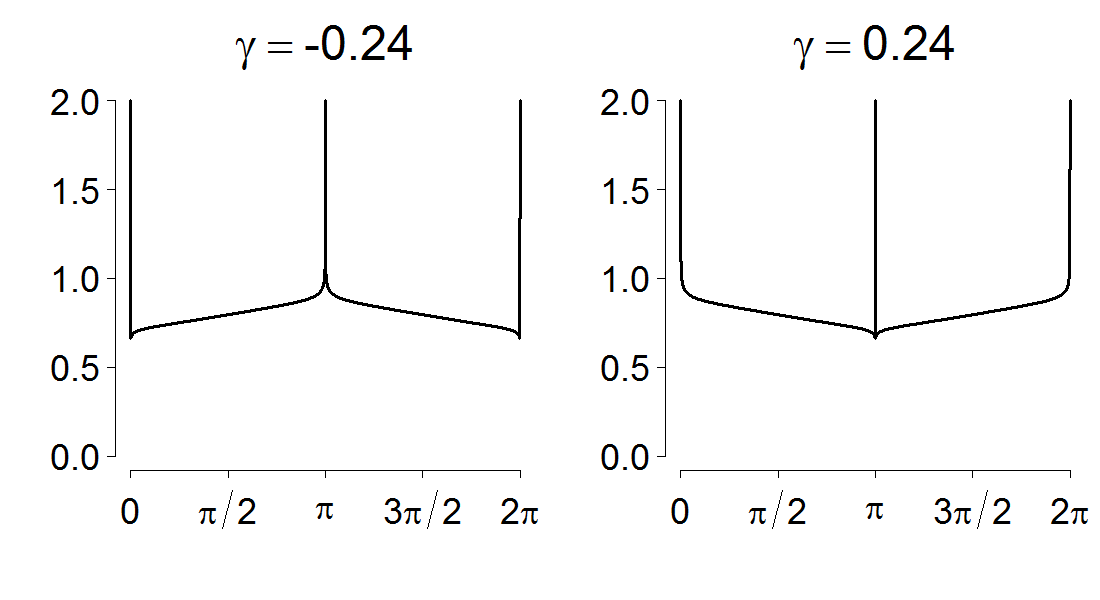}}
    \subfloat[]{\includegraphics[width = 0.25\textwidth]{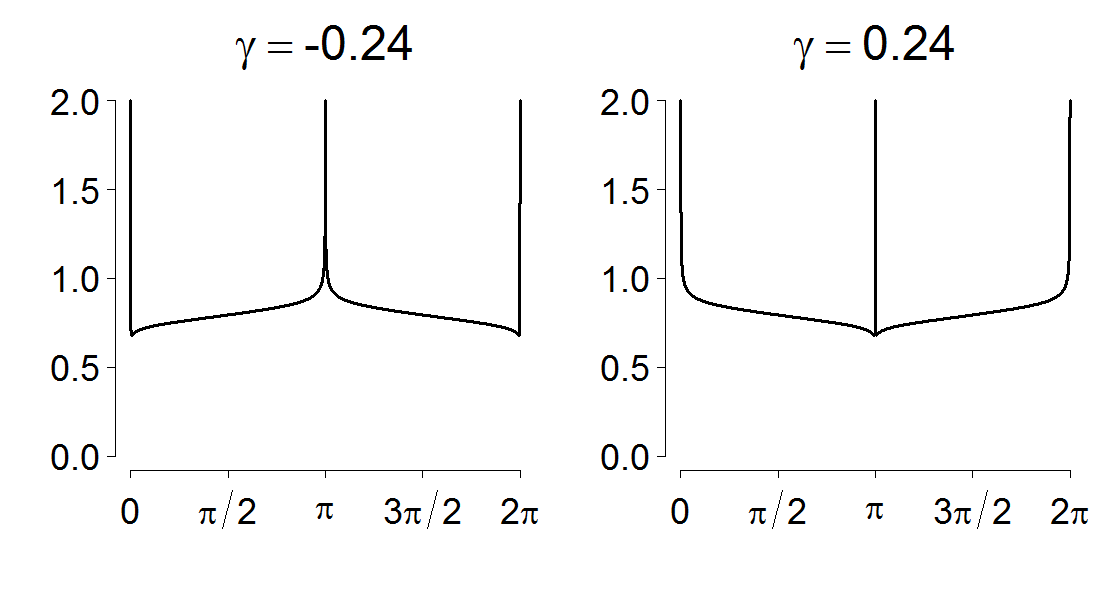}}
    \subfloat[]{\includegraphics[width = 0.25\textwidth]{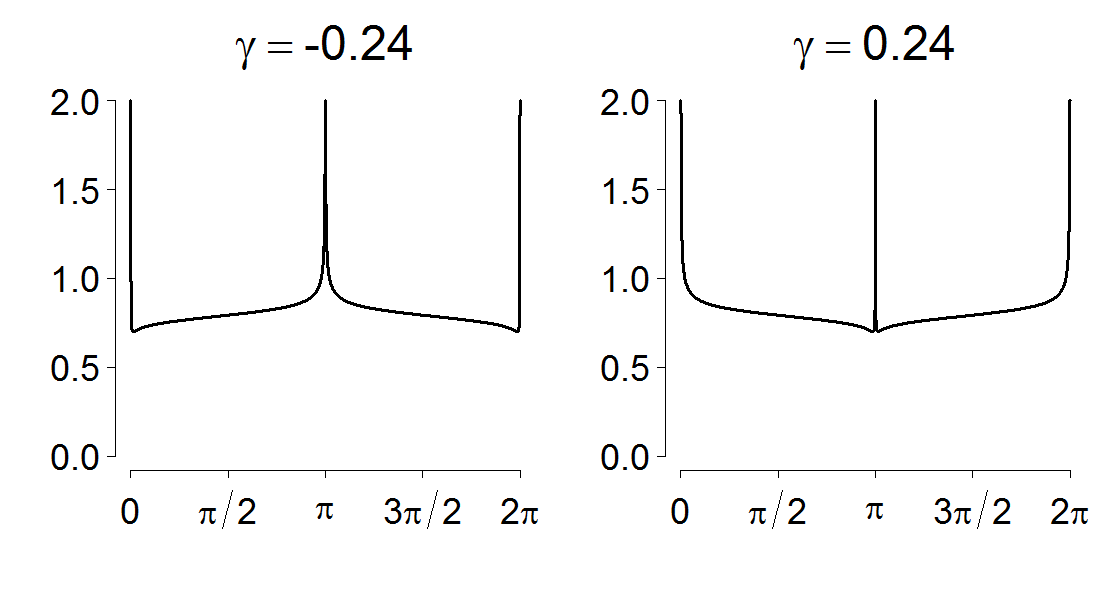}}
  }
  \caption{ Theoretical spectral density function of the process $\{\ln(X_t^2)\}_{t\in\Z}$, when $\{X_t\}_{t \in
      \mathds{Z}}$ is an SFIEGARCH$(0,d,0)_s$ process with $s=2$, $\theta = 0.25$ and $\gamma \in\{-0.24,0.24\}$
    (in each panel, from left to right).  The parameter $d$ is set as follows: in (a) $d = 0.1$, in (b) $d=0.2$,
    in (c) $d=0.3$ and in (d) $d=0.4$. } \label{fig:spectral}
\end{figure}

\begin{figure}[!htbp]
  \centering
  \mbox{
    \hspace{-4pt}\subfloat[]{\includegraphics[width = 0.25\textwidth]{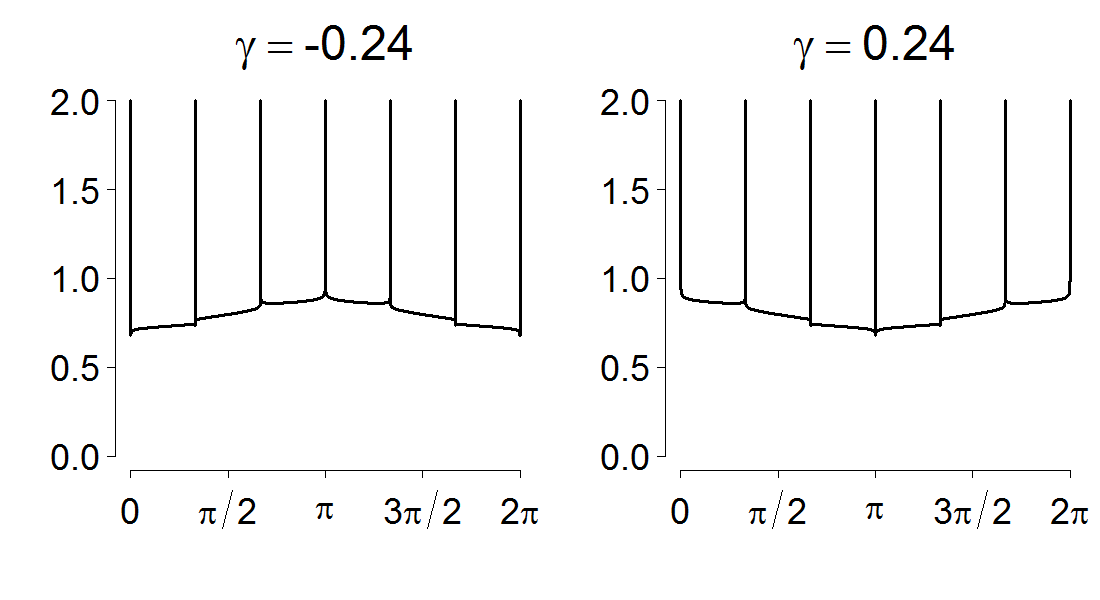}}
    \subfloat[]{\includegraphics[width = 0.25\textwidth]{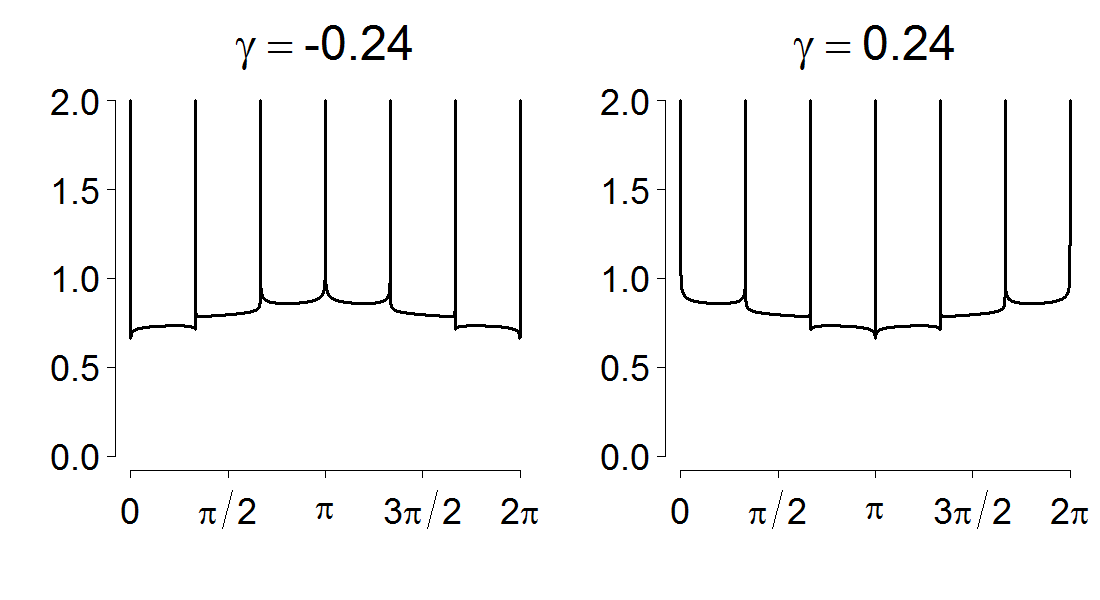}}
    \subfloat[]{\includegraphics[width = 0.25\textwidth]{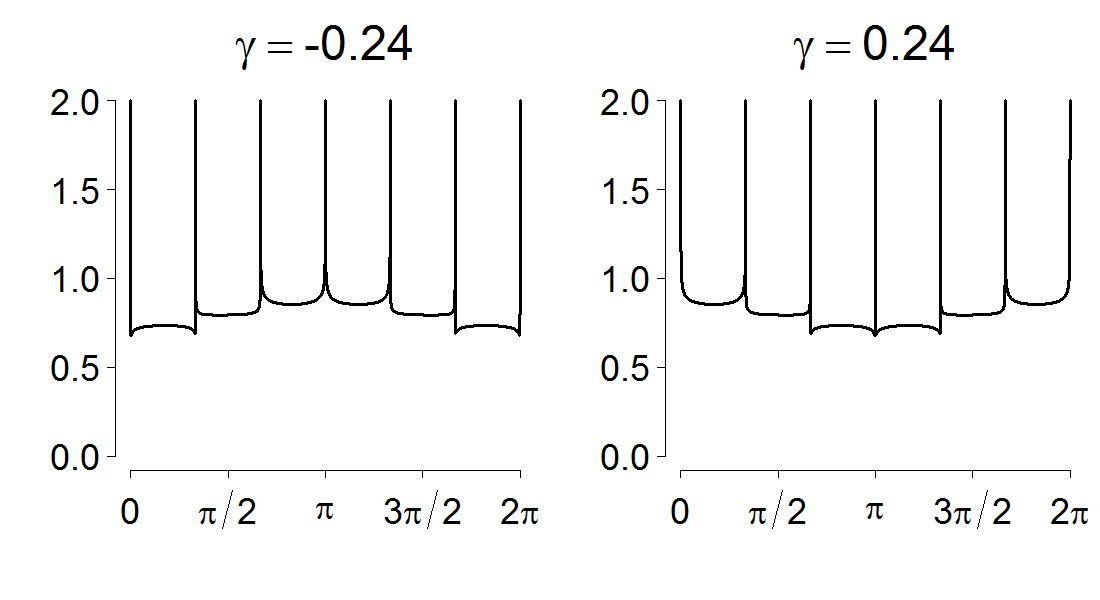}}
    \subfloat[]{\includegraphics[width = 0.25\textwidth]{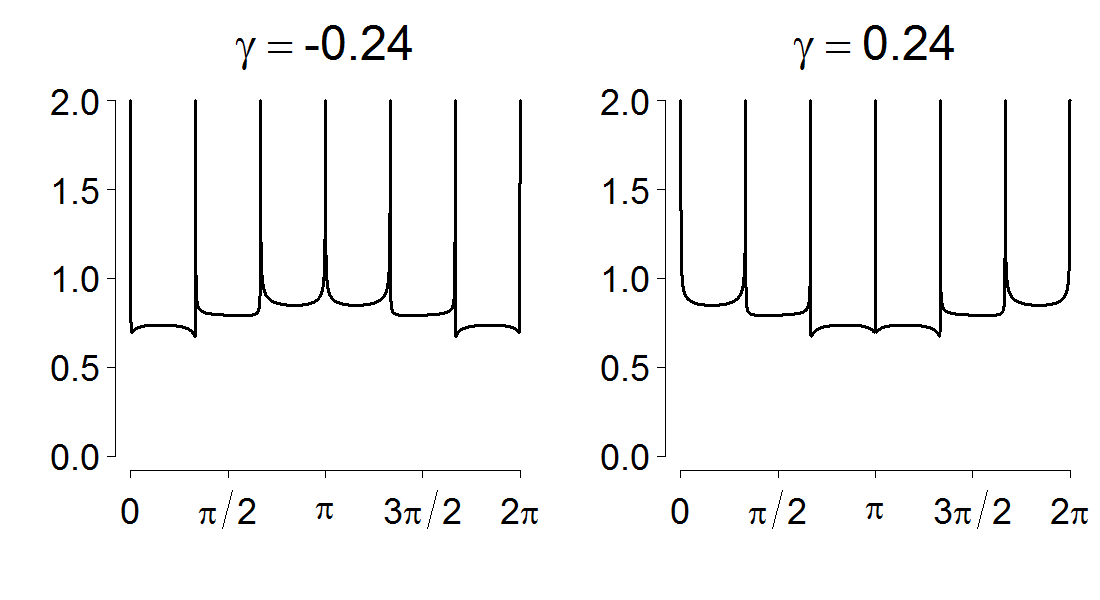}}
  }
  \caption{ Theoretical spectral density function of the process $\{\ln(X_t^2)\}_{t\in\Z}$, when $\{X_t\}_{t \in
      \mathds{Z}}$ is an SFIEGARCH$(0,d,0)_s$ process with $s=6$, $\theta = 0.25$ and $\gamma \in\{-0.24,0.24\}$
    (in each panel, from left to right).  The parameter $d$ is set as follows: in (a) $d = 0.1$, in (b) $d=0.2$,
    in (c) $d=0.3$ and in (d) $d=0.4$.}  \label{fig:spectral2}
\end{figure}

\begin{example}
  Figures \ref{fig:spectral} and \ref{fig:spectral2} show the spectral density function of the process
  $\{\ln(X_t^2)\}_{t\in\Z}$, where $\{X_t\}_{t \in \mathds{Z}}$ is an SFIEGARCH$(0,d,0)_s$ with different
  parameter values and $Z_0\sim \mathcal{N}(0,1)$. Since $Z_0$ is a symmetric random variable, $\E(Z_0|Z_0|) =
  0$ and the function $f_{\ln(\X_t^2)}(\cdot)$ is symmetric in $\theta$. Thus, in both figures, we fixed
  $\theta = 0.25$.  In Figure \ref{fig:spectral} we fix $s=2$ and, in Figure \ref{fig:spectral2}, we consider
  $s = 6$.  For each figure, $\gamma \in\{ -0.24,0.24\}$ and $d\in\{0.1,0.2,0.3,0.4\}$.  From Figures
  \ref{fig:spectral} and \ref{fig:spectral2}, one observes that, for each fixed $d$ and $s$, the behavior of
  the function completely changes as $\gamma$ changes from $-0.24$ to $0.24$ (left to right).  While for
  $\gamma = -0.24$ the function attains its minimum in the region close to zero, for $\gamma = 0.24$ the
  minimum is attained close to $\pi$.
\end{example}

\begin{figure}[!htb]
  \centering
  \subfloat[]{\includegraphics[height=0.1\textheight]{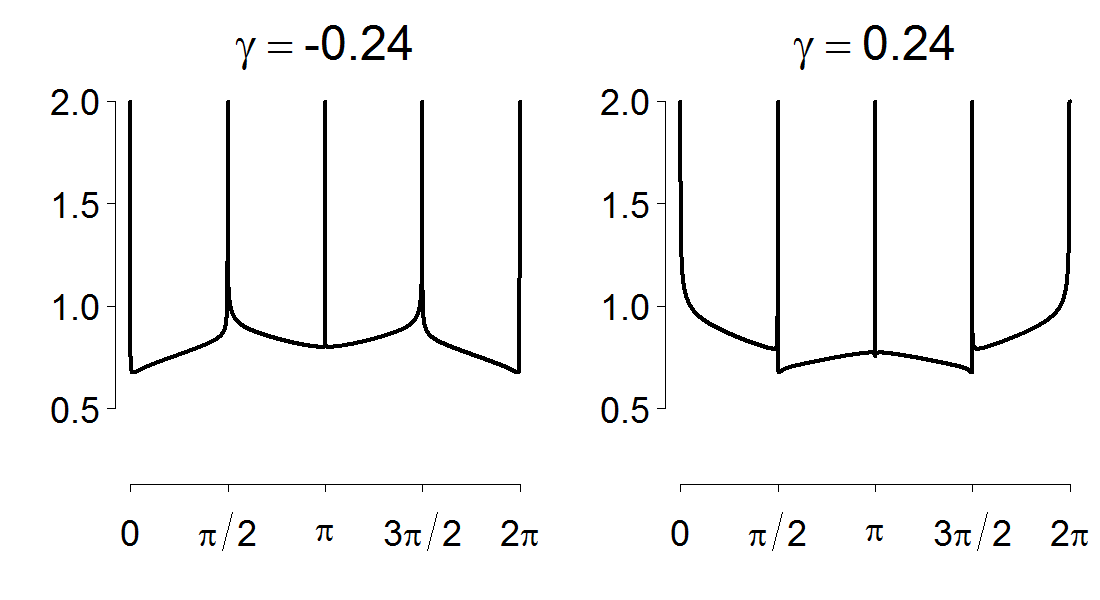}}\hspace{0.5cm}
  \subfloat[]{\includegraphics[height=0.1\textheight]{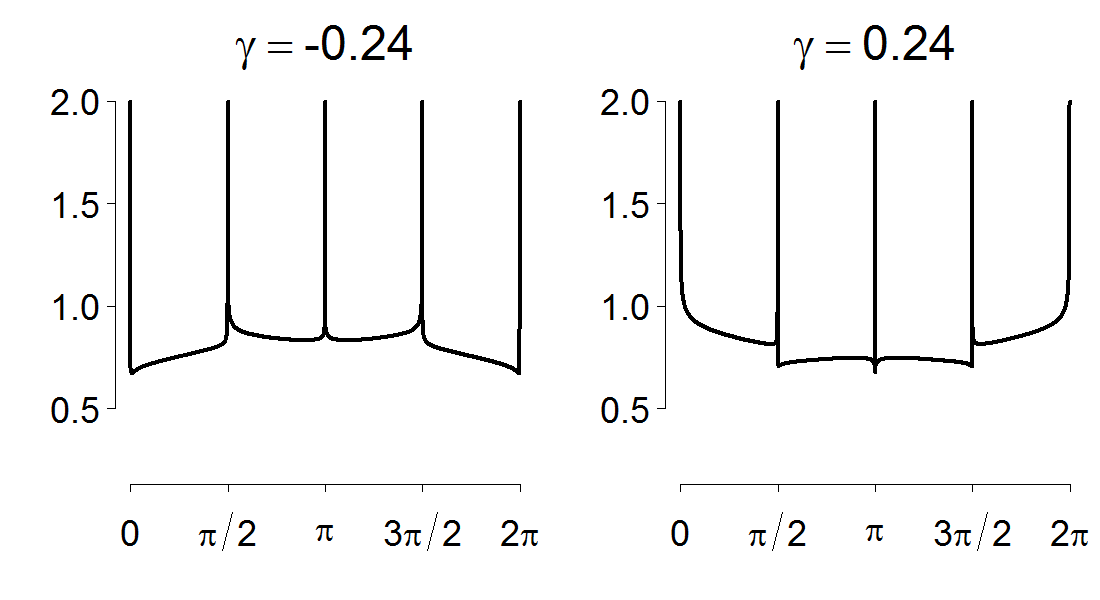}}\hspace{0.5cm}
  \subfloat[]{\includegraphics[height=0.1\textheight]{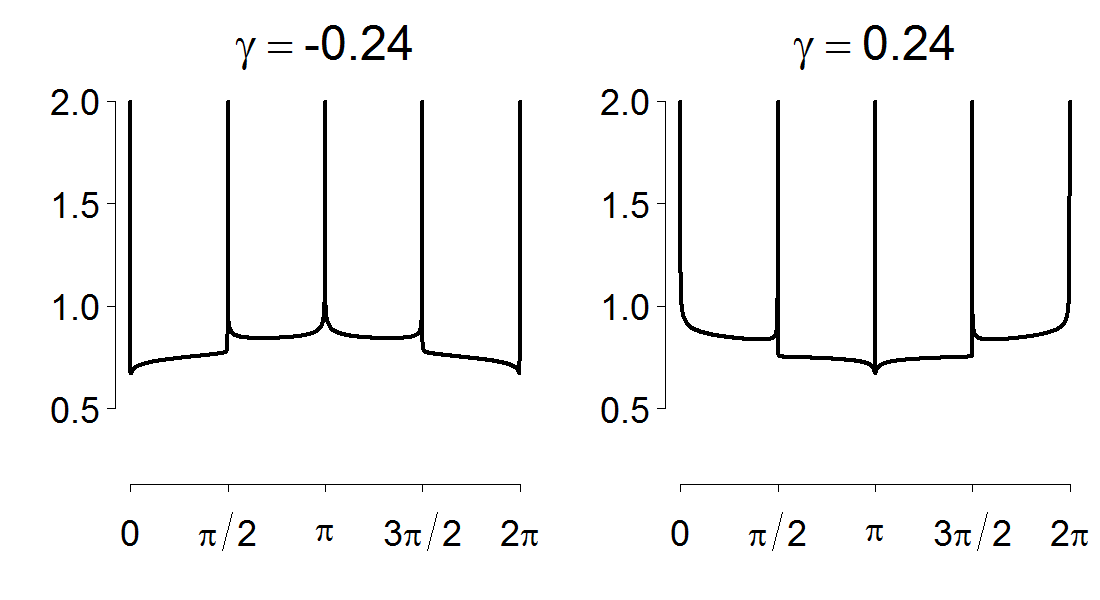}}\\
  \subfloat[]{\includegraphics[height=0.1\textheight]{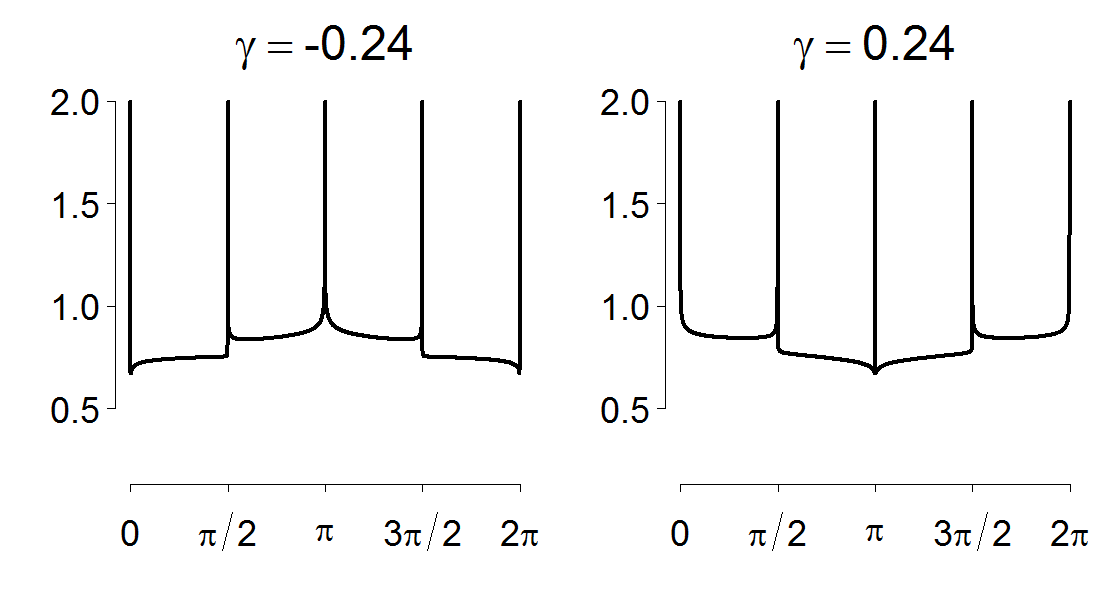}}\hspace{0.5cm}
  \subfloat[]{\includegraphics[height=0.1\textheight]{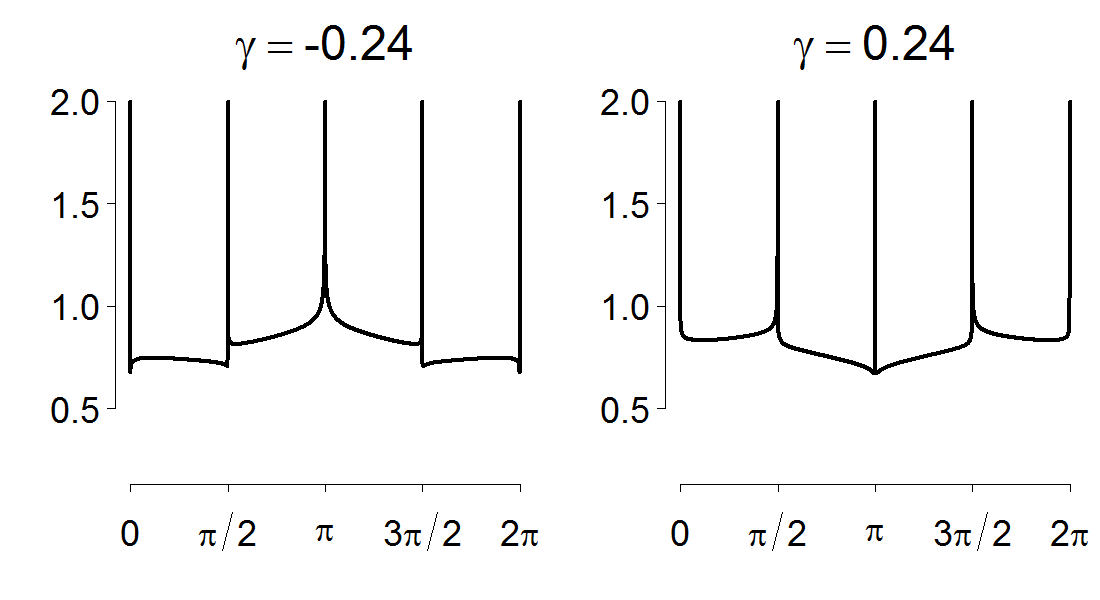}}\hspace{0.5cm}
  \subfloat[]{\includegraphics[height=0.1\textheight]{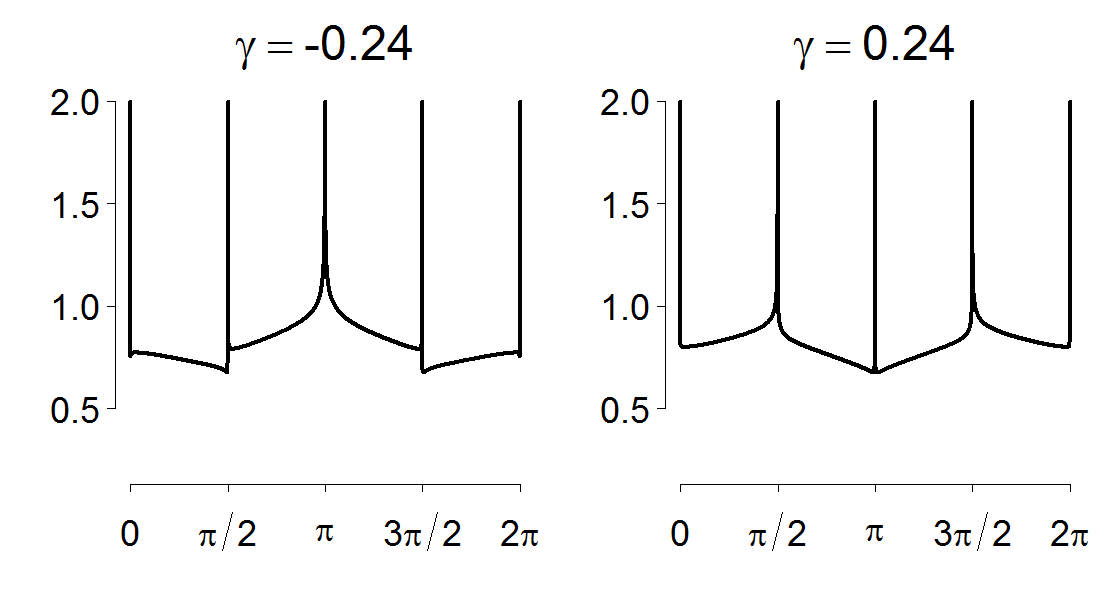}}

  \caption{ Theoretical spectral density function of $\{\ln(X_t^2)\}_{t\in\Z}$, when $\{X_t\}_{t \in \mathds{Z}}$
    is an SFIEGARCH$(1,d,0)_s$, with $s=4$, $\omega = 0$, $d=0.25$, $\theta = 0.25$, $\gamma \in\{-0.24,0.24\}$
    (in each panel, from left to right). The parameter $\alpha_1$ is set as follows: in (a) $\alpha_1 = -0.9$,
    in (b) $\alpha_1 = -0.5$, in (c) $\alpha_1 = -0.1$, in (d) $\alpha_1 = 0.1$, in (e) $\alpha_1 = 0.5$ and in
    (f) $\alpha_1 = 0.9$.  } \label{fig:spectral3}
\end{figure}

\begin{figure}[!htbp]
  \centering
  \begin{tabular}{ccccc}
    {\tiny (-0.5, -0.9)} &
    {\tiny (-0.5, -0.1)} &
    {\tiny (-0.5, 0.1)} &
    {\tiny (-0.5, 0.5)} &
    {\tiny (-0.5, 0.9)} \vspace{-0.1cm}\\

    \includegraphics[width=0.19\textwidth]{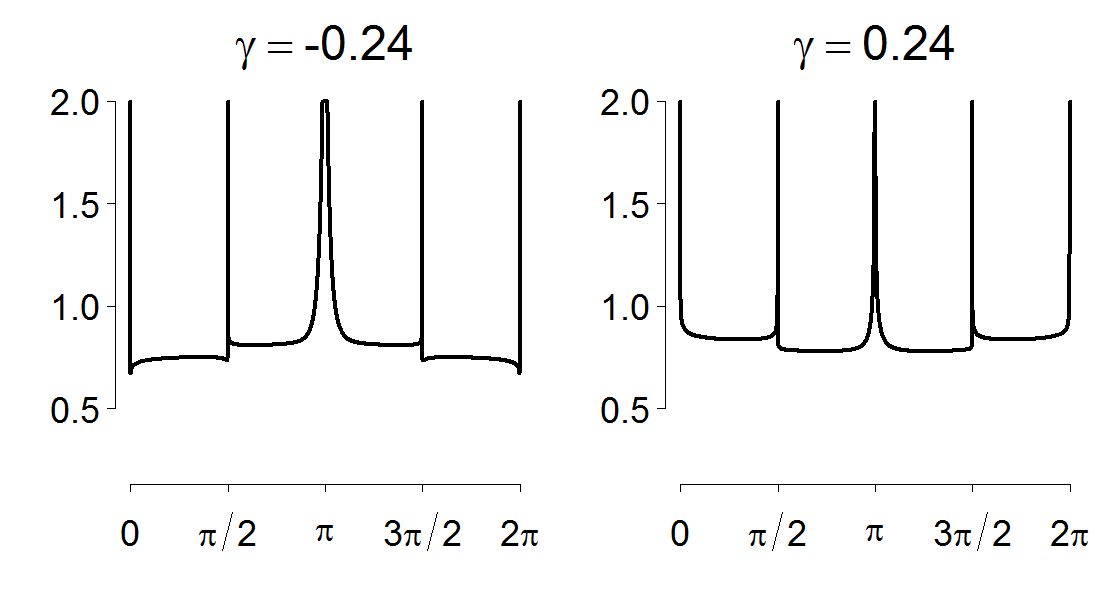}
    &
    \hspace{-12pt}\includegraphics[width=0.19\textwidth]{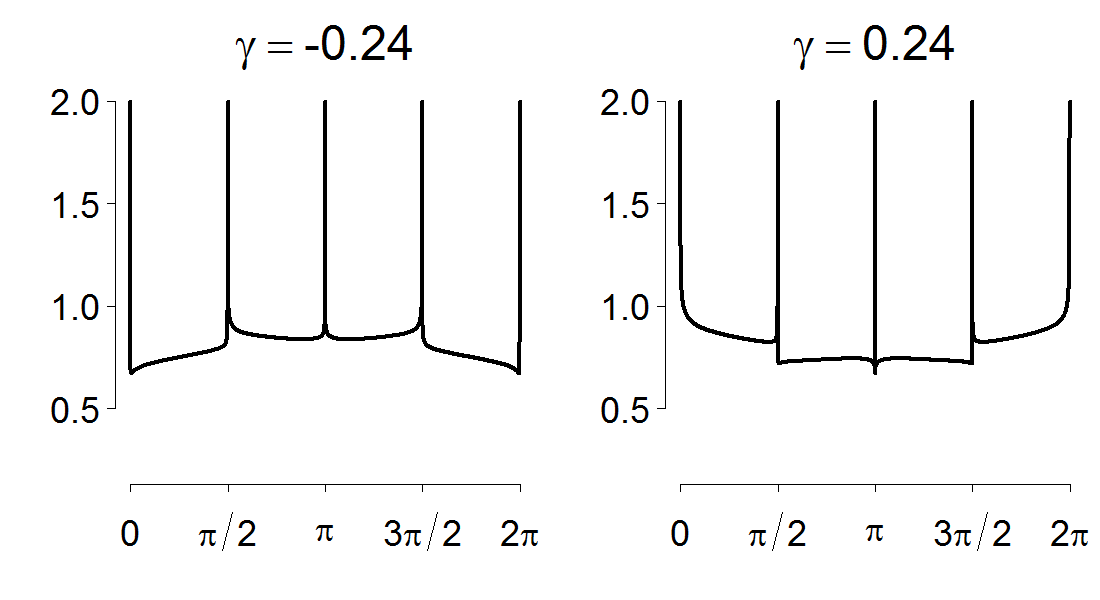}
    &
    \hspace{-12pt}\includegraphics[width=0.19\textwidth]{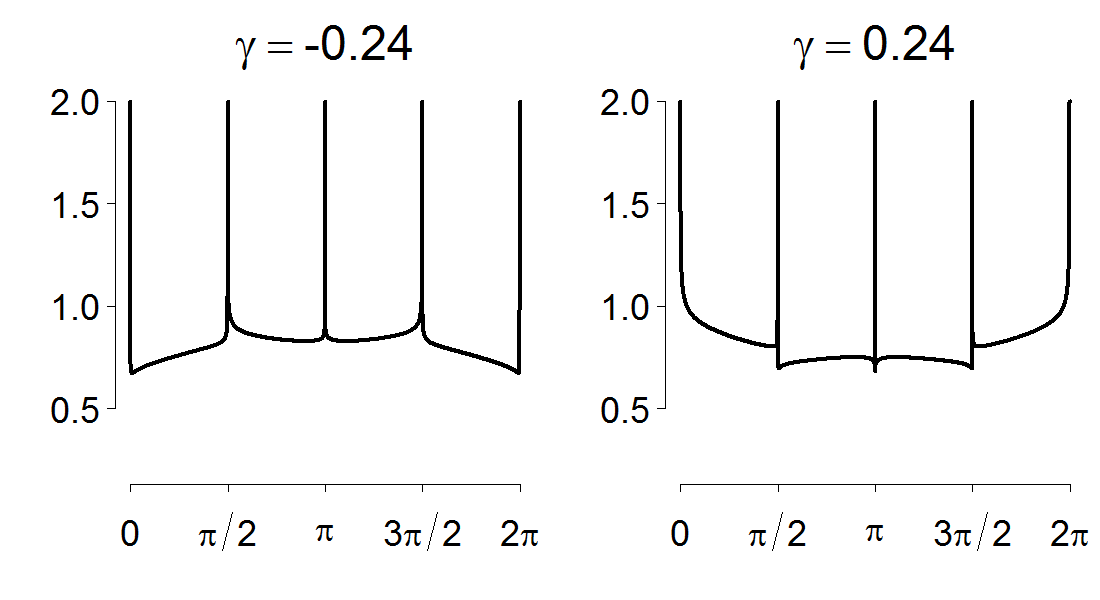}
    &
    \hspace{-12pt}\includegraphics[width=0.19\textwidth]{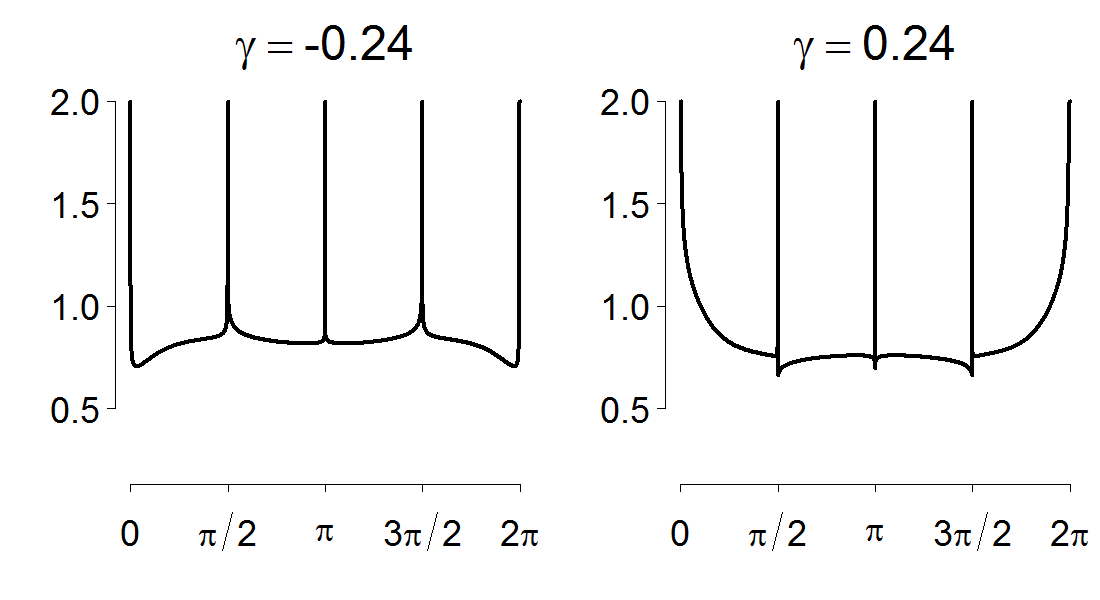}
    &
    \hspace{-12pt}\includegraphics[width=0.19\textwidth]{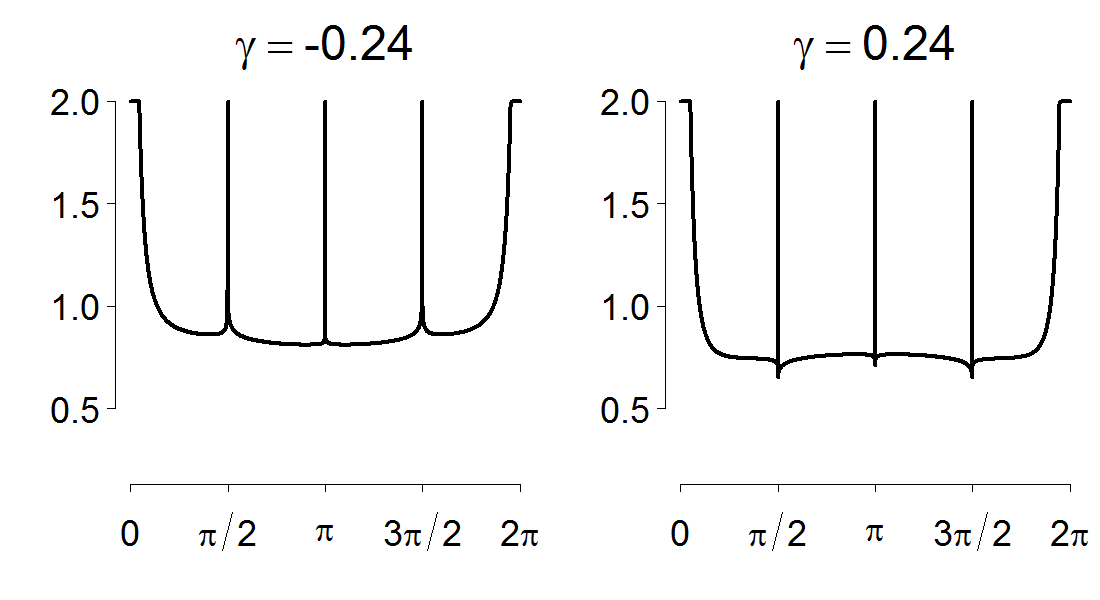}\\
    {\footnotesize (a)} &
    {\footnotesize (b)} &
    {\footnotesize (c)} &
    {\footnotesize (d)} &
    {\footnotesize (e)} \vspace{0.4cm} \\

    {\tiny (0.5, -0.9)} &
    {\tiny (0.5, -0.5)} &
    {\tiny (0.5, -0.1)} &
    {\tiny (0.5, 0.1)} &
    {\tiny (0.5, 0.9)} \vspace{-0.1cm}\\

    \includegraphics[width=0.19\textwidth]{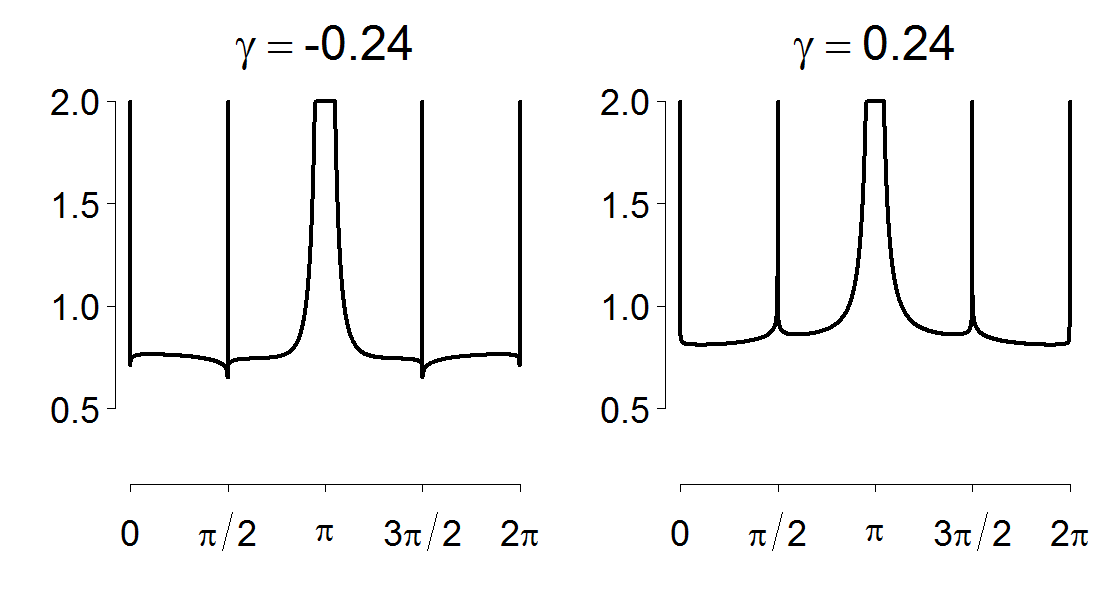}
    &
    \hspace{-12pt}\includegraphics[width=0.19\textwidth]{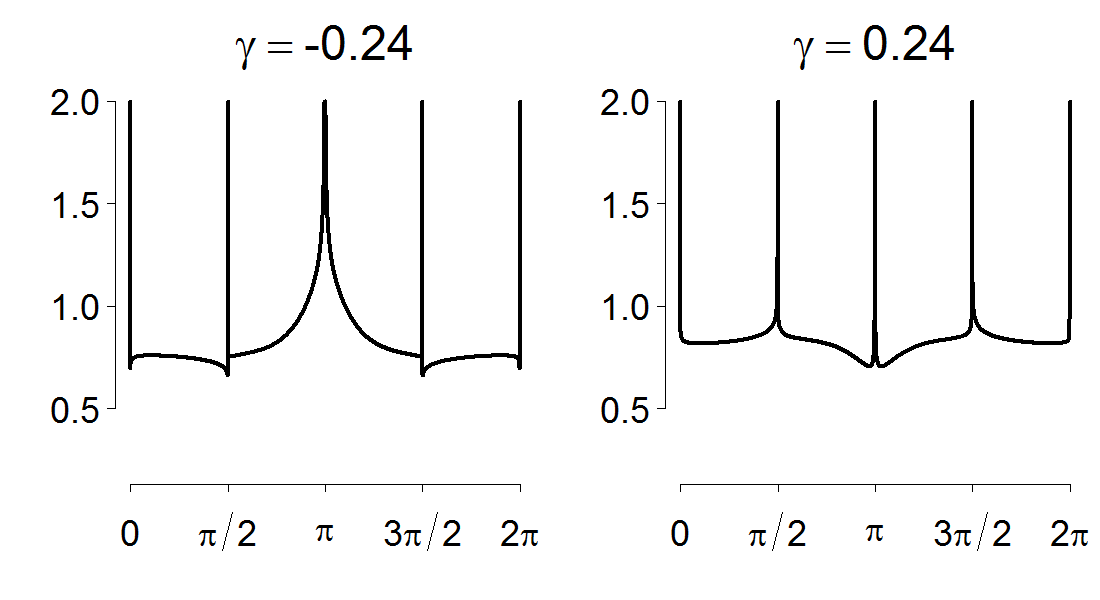}
    &
    \hspace{-12pt}\includegraphics[width=0.19\textwidth]{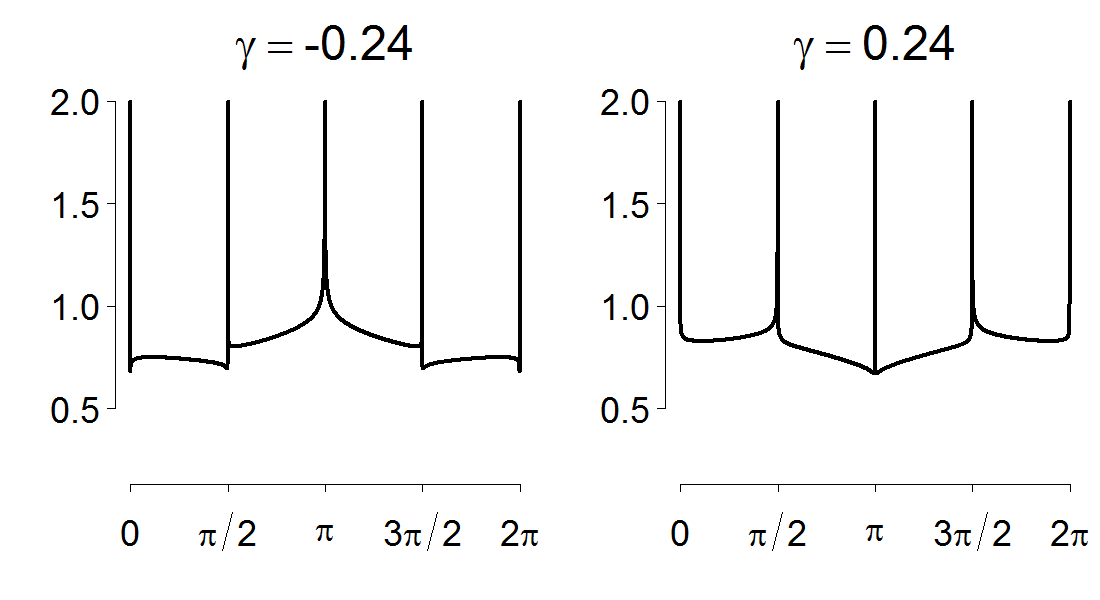}
    &
    \hspace{-12pt}\includegraphics[width=0.19\textwidth]{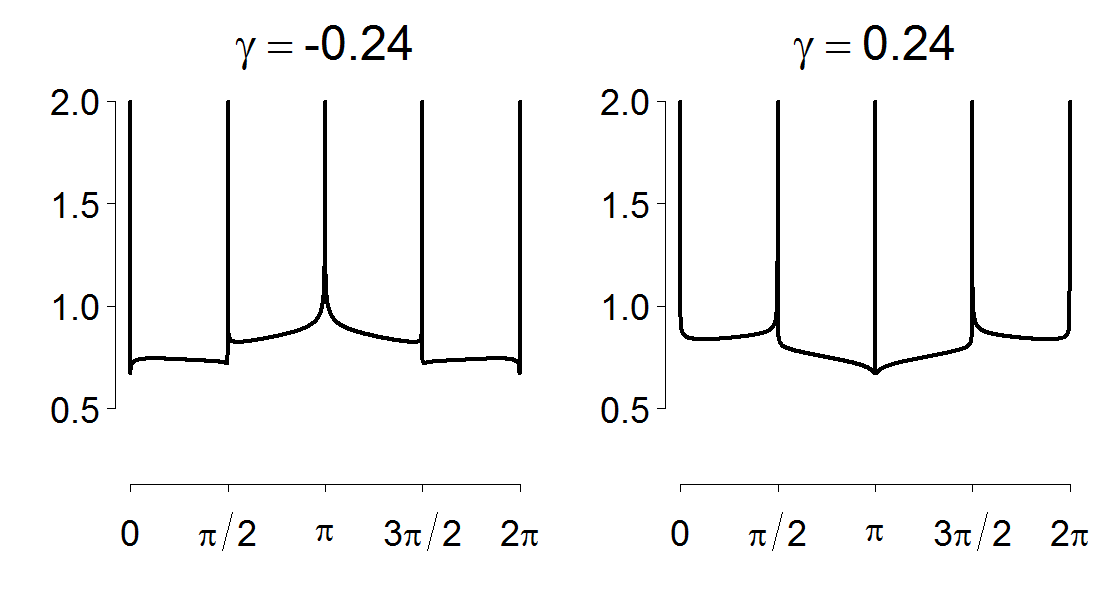}
    &
    \hspace{-12pt} \includegraphics[width=0.19\textwidth]{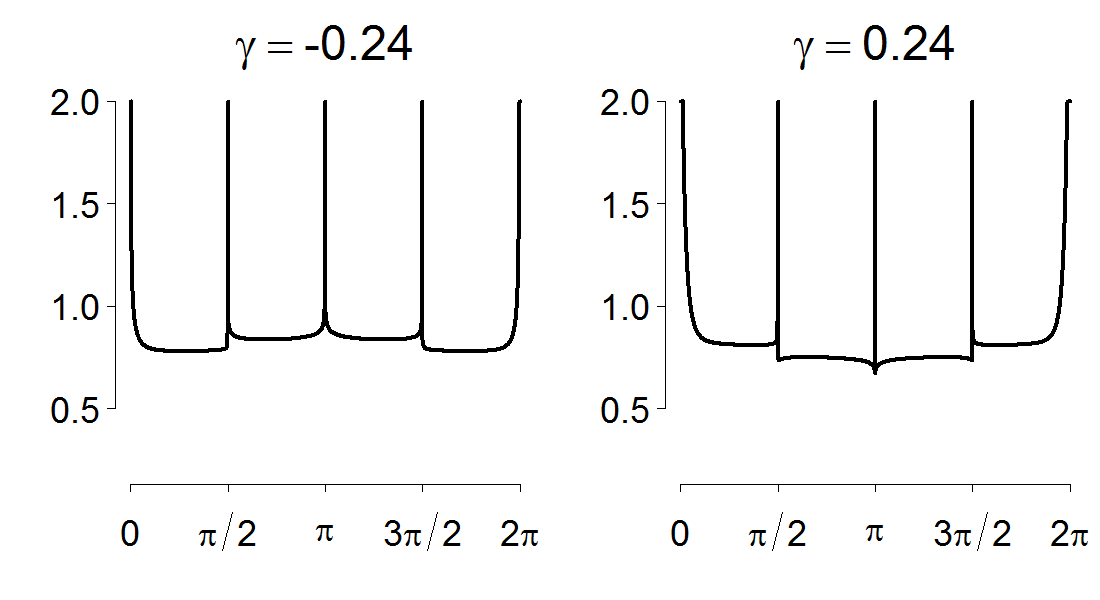}\\

    {\footnotesize (f)} &
    {\footnotesize (g)} &
    {\footnotesize (h)} &
    {\footnotesize (i)} &
    {\footnotesize (j)} \\

  \end{tabular}
  \caption{ Theoretical spectral density function of $\{\ln(X_t^2)\}_{t\in\Z}$, when $\{X_t\}_{t \in \mathds{Z}}$ is an
    SFIEGARCH$(1,d,1)_s$, with $s=4$, $\omega = 0$, $d=0.25$, $\theta = 0.25$.  For each panel $\gamma
    \in\{-0.24,0.24\}$ (from left to right).  For all panels in the first row $\alpha_1 = -0.5$ and for all
    panels in the second row $\alpha = 0.5$.  For each row $\beta_1 \in\ \{-0.9,-0.5,-0.1,0.1,0.5,0.9\}$ (from
    left to right panel).} \label{fig:spectral4}
\end{figure}

\begin{example}
  Figures \ref{fig:spectral3} - \ref{fig:spectral4} present the spectral density function of
  $\{\ln(X_t^2)\}_{t\in\Z}$, where $\{X_t\}_{t \in \mathds{Z}}$ is an SFIEGARCH$(p,d,q)_s$, with
  $Z_0\sim\mathcal{N}(0,1)$, $p,q\in\{0,1\}$ (not both equal to zero), $s=4$, $d=0.25$, $\theta = 0.25$,
  $\gamma \in\{-0.24,0.24\}$ and $\alpha_1,\beta_1 \in\ \{-0.9,-0.5,-0.1,0.1,0.5,0.9\}$. For these figures the
  parameters values increase from left to right and from top to bottom.
\end{example}

\begin{figure}[!htb]
  \centering
  \subfloat[]{\includegraphics[height=0.1\textheight]{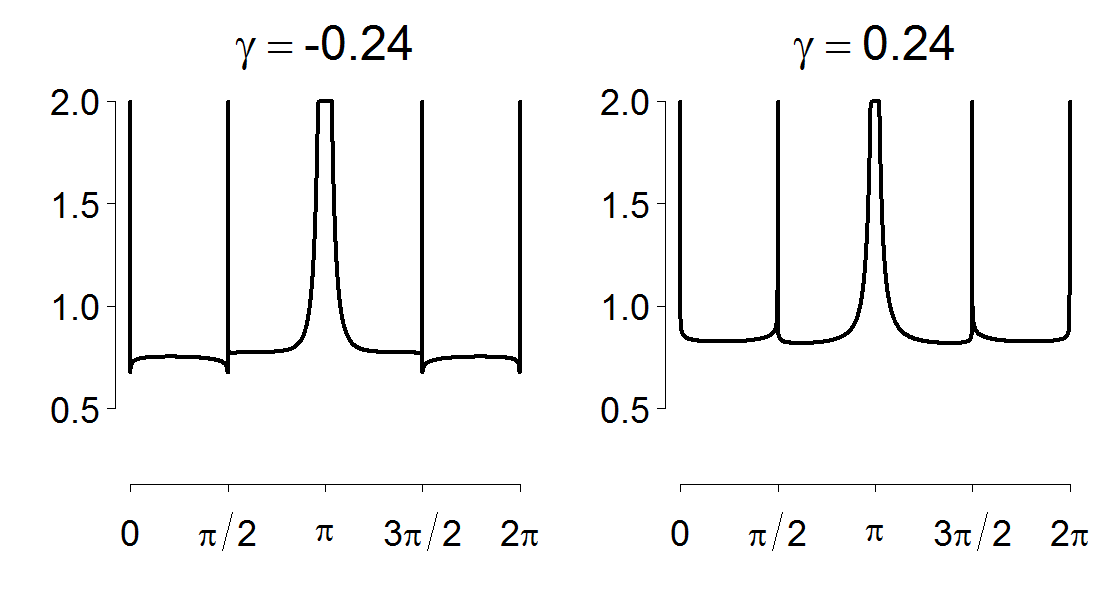}}\hspace{0.5cm}
  \subfloat[]{\includegraphics[height=0.1\textheight]{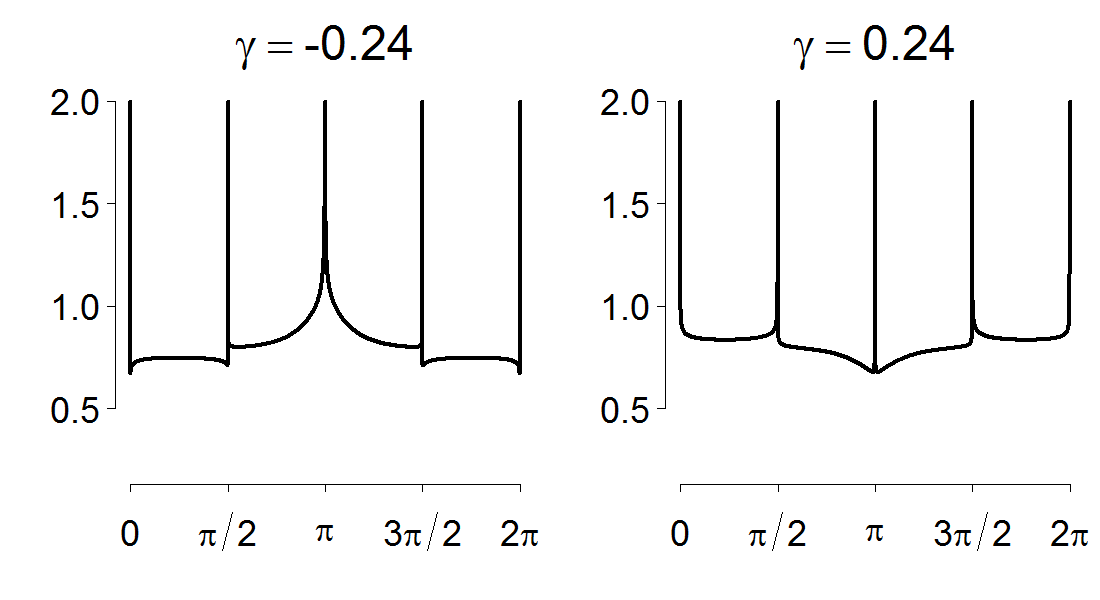}}\hspace{0.5cm}
  \subfloat[]{\includegraphics[height=0.1\textheight]{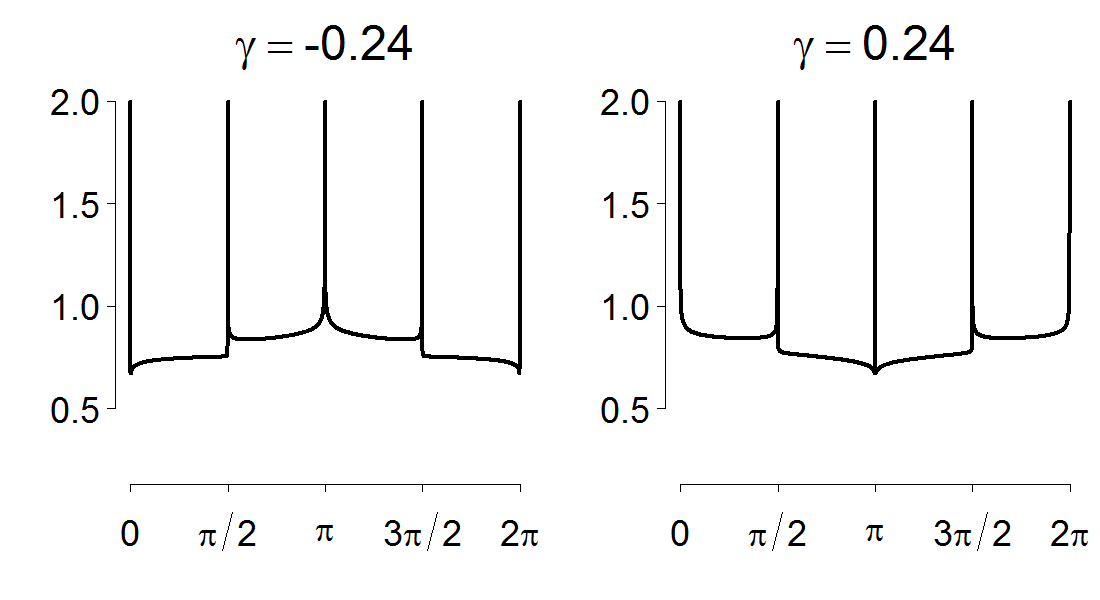}}\\
  \subfloat[]{\includegraphics[height=0.1\textheight]{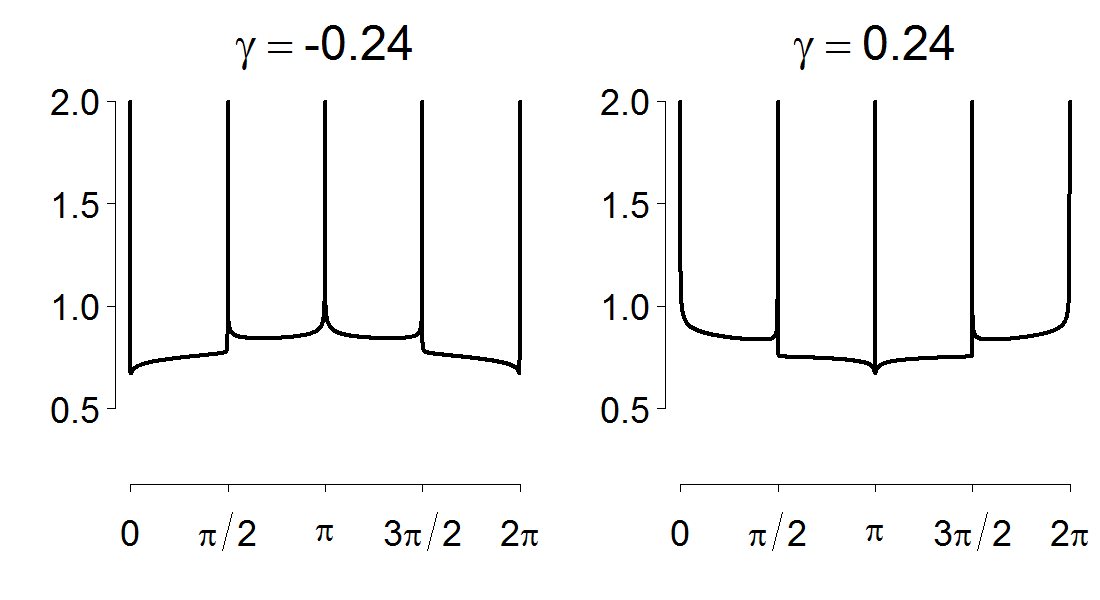}}\hspace{0.5cm}
  \subfloat[]{\includegraphics[height=0.1\textheight]{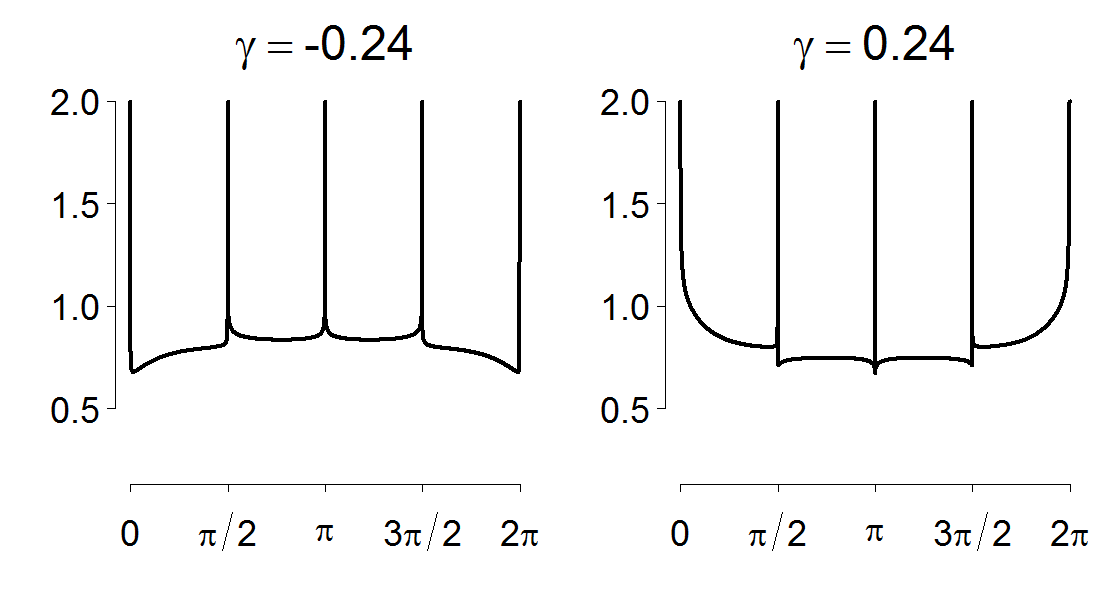}}\hspace{0.5cm}
  \subfloat[]{\includegraphics[height=0.1\textheight]{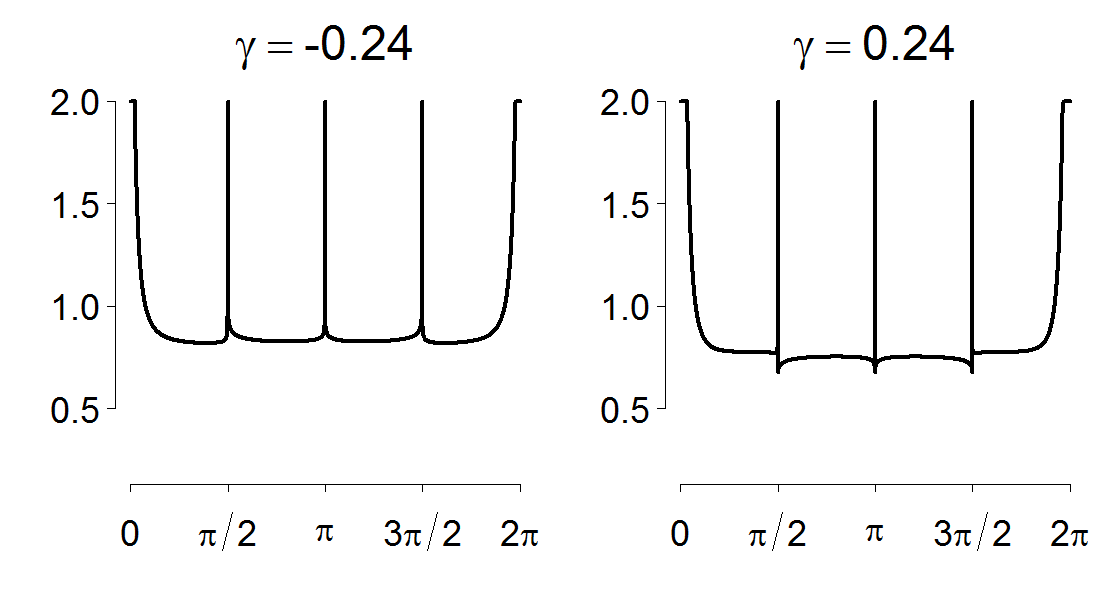}}

  \caption{ Theoretical spectral density function of $\{\ln(X_t^2)\}_{t\in\Z}$, when $\{X_t\}_{t \in \mathds{Z}}$ is an
    SFIEGARCH$(1,d,0)_s$, with $s=4$, $\omega = 0$, $d=0.25$, $\theta = 0.25$, $\gamma \in\{-0.24,0.24\}$ (in
    each panel, from left to right). The parameter $\beta_1$ is set as follows: in (a) $\beta_1 = -0.9$, in (b)
    $\beta_1 = -0.5$, in (c) $\beta_1 = -0.1$, in (d) $\beta_1 = 0.1$, in (e) $\beta_1 = 0.5$ and in (f)
    $\beta_1 = 0.9$.}
\end{figure}

Although Figure \ref{fig:spectral4} presents only the graphs for $\alpha_1 \in\{-0.5,0.5\}$, in the sequel we
discuss the behavior of $f_{\ln(\X_t^2)}(\cdot)$, for all $\alpha_1 \in\{-0.9, -0.5, -0.1, 0.1, 0.5,
0.9\}$. The remaining graphs are available upon request.  From Figures \ref{fig:spectral3} -
\ref{fig:spectral4}, one concludes the following:

\begin{itemize}
\item if $p=1$ and $q = 0$
  \begin{itemize}
  \item the region where $f_{\ln(\X_t^2)}(\cdot)$ attains its minimum depends not only on the sign of
    $\gamma$, but also on the sign of $\alpha_1$. If $\gamma < 0$, the minimum is attained either close to
    zero or close to $\pi/2$. If $\gamma >0$, the minimum is attained either close to $\pi/2$ or close to
    $\pi$;

  \item for $\gamma$ fixed, the graph of $f_{\ln(\X_t^2)}(\cdot)$ slowly changes its behavior in the regions
    around the seasonal frequencies, as $\alpha_1$ increases;

  \item almost no difference is observed in the graphs of $f_{\ln(\X_t^2)}(\cdot)$ when $\alpha_1$ changes
    from $-0.1$ to $0.1$. Actually, for these values of $\alpha_1$, the function behaves as in the case
    $p=0=q$ (the graphs can be obtained upon request);

  \item the graph of $f_{\ln(\X_t^2)}(\cdot)$ for $\alpha_1 = -0.9$ and $\gamma = -0.24$ is very similar to
    the graph of the same function for $\alpha_1 = 0.9$, $\gamma = 0.24$.
  \end{itemize}

\item if $p=0$ and $q = 1$
  \begin{itemize}
  \item the region where $f_{\ln(\X_t^2)}(\cdot)$ attains its minimum depends on the sign of $\gamma$ and on
    the sign of $\beta_1$. The behavior is similar to the case $p=1$ and $q = 0$;

  \item for $\gamma$ fixed, the changes in the graph of $f_{\ln(\X_t^2)}(\cdot)$, in the regions around the
    seasonal frequencies, are much more visible than in the case $p=1$ and $q=0$;

  \item for $\beta_1\in\{-0.1,0.1\}$ the graphs are very similar to the case $p=0$ and $q=1$, with
    $\alpha_1\in\{-0.1,0.1\}$.

  \item The graphs of $f_{\ln(\X_t^2)}(\cdot)$ for $\beta_1 = -0.5$ are almost identical to the graphs of this
    function in the case $p=1$ and $q=0$ with $\alpha_1 = 0.5$. The same similarity is observed between the
    graphs of this function for $\beta_1 = 0.5$ and $\alpha_1 = -0.5$.
  \end{itemize}

\item if $p=1$ and $q = 1$
  \begin{itemize}
  \item for $\gamma$ and $\theta$ fixed, the changes in the graph of $f_{\ln(\X_t^2)}(\cdot)$ are more visible
    as $\beta_1$ changes than when $\alpha_1$ does;

  \item for $\alpha_1,\beta_1\in\{-0.1,0.1\}$ the graphs are very similar to the case $p=0$ and $q=1$ (same
    occurs with $p=1$ and $q=0$);

  \item generally, the graphs do not show any peculiar characteristic that is not present in the cases $p=0$
    or $q=0$.
  \end{itemize}
\end{itemize}

\section{Forecasting}\label{sec:forecasting}

Let $\{X_t\}_{t \in \mathds{Z}}$ be an SFIEGARCH$(p,d,q)_s$ process and $\{\mathcal{F}_{t}\}_{t\in\mathds{Z}}$
be the filtration defined by $\mathcal{F}_{t} := \sigma(\{Z_s\}_{s\leq t})$.  Notice that, by considering the
same argument as in the proof of lemma 1 in \cite{LP2013}, one can show that $\{X_t\}_{t \in \mathds{Z}}$ is a
martingale difference with respect to $\{\mathcal{F}_{t}\}_{t\in\mathds{Z}}$.  In this case, the best
predictor (in terms of the mean square error measure) for $X_{t+h}$, given $\mathcal{F}_{t}$, is
$\mathds{E}(X_{t+h}|\mathcal{F}_t) = 0$, for all $h > 0$ and $t\in\Z$.

Since the $h$-step ahead predictor for $\{X_t\}_{t \in \mathds{Z}}$ is always zero (the mean of the process),
the aim of this section is to derive  expressions for the $h$-step ahead forecast for the processes
$\{\sigma_{t}^2\}_{t\in\mathds{Z}}$, $\{X_{t}^2\}_{t\in\mathds{Z}}$, $\{\ln(\sigma_{t}^2)\}_{t\in\mathds{Z}}$
and $\{\ln(X_{t}^2)\}_{t\in\mathds{Z}}$, for any $h> 0 $.  The approach considered in this work is slightly
different from \cite{LP2013}.  In \cite{LP2013} two different $h$-step ahead predictors for $\sigma_{t+h}^2$
(consequently, for  $X_{t+h}^2$) were proposed, both based on the $h$-step ahead predictor for
$\ln(\sigma_{t+h}^2)$.  Here we provide the exact formula for $\mathds{E}(\sigma_{t+h}^2|\mathcal{F}_t) $, for
any $h>0$ and $t\in \Z$, and the relation between this expression and the ones in \cite{LP2013}.

\begin{remark}
  In the sequel we consider the following notation, which is the same as in \cite{LP2013}.  Let $Y_t$, for
  $t\in\mathds{Z}$, denote any random variable.  Then
  \begin{itemize}
  \item the symbol ``\^{}'' denotes the $h$-step ahead forecast defined in terms of the conditional
    expectation, that is, $\hat Y_{t+h} = \mathds{E}(Y_{t+h}|\mathcal{F}_{t})$.  Notice that this is the best
    linear (or non-linear) predictor in terms of mean square error value;

  \item the symbols ``\~{}'' (e.g. $\tilde Y_{t+h}$) and ``\v{}'' (e.g. $\check Y_{t+h}$) denote alternative
    estimators;

  \item $\hat{\ln}(Y_{t+h})$ denotes the $h$-step ahead forecast of $\ln(Y_{t+h})$ (analogously,  for ``\~{}''
    and ``\v{}'');

  \item we follow the approach usually considered in the literature and denote the $h$-step ahead forecast of
    $Y_{t+h}^2$ as $\hat Y_{t+h}^2$ instead of $\widehat{Y_{t+h}^2}$ and,  to avoid confusion, we
    will denote the square of $\hat Y_{t+h}$ as $(\hat Y_{t+h})^2$ (analogously,  for ``\~{}'' and ``\v{}'').
  \end{itemize}
\end{remark}

To obtain the predictors for the processes $\{\sigma_{t}^2\}_{t\in\mathds{Z}}$ and
$\{X_{t}^2\}_{t\in\mathds{Z}}$ observe that the i.i.d. property of $\{Z_t \}_{t\in \Z}$ implies that
$\mathds{E}(Z_{n+h}^2|\mathcal{F}_n) = \mathds{E}(Z_{n+h}^2) = 1$,   $\sigma_{n+1}^2$ is
$\mathcal{F}_n$-measurable and  $\sigma_n^2$ and $Z_n^2$ are independent, for all $n\in \Z$ and $h > 0$.
Therefore, the $h$-step ahead forecast for $X_{n+h}^2$ given $\mathcal{F}_{n}$ is given by
\[
\hat{X}_{n+h}^2 := \mathds{E}(X_{n+h}^2|\mathcal{F}_{n}) =
\mathds{E}(\sigma_{n+h}^2|\mathcal{F}_{n}) := \hat{\sigma}_{n+h}^2, \quad \mbox{ for all} \quad h>0.
\]
In particular, $\hat{\sigma}_{n+1}^2 =\sigma_{n+1}^2$ and hence the $1$-step ahead forecast of $X_{n+1}^2$,
given $\mathcal{F}_{n}$,  is simply $\hat{X}_{n+1}^2 =\sigma_{n+1}^2$.  Theorem \ref{thm:predictor} provides the
general formula for $\hat{\sigma}_{n+h}^2$, when $h> 1$.

\begin{thm}\label{thm:predictor}
  Let $\{X_t\}_{t \in \mathds{Z}}$ be a stationary \emph{SFIEGARCH}$(p,d,q)_s$ process with
  $\mathds{E}(\sigma_t^2) < \infty$.  Then, for any fixed $n\in\mathds{Z}$, the $h$-step ahead forecast of
  $\sigma_{n+h}^2$ (consequently, for $X_{n+h}^2$), given $\mathcal{F}_{n}$,  can be expressed as
  \begin{equation}\label{eq:predictor}
    \hat\sigma_{n+h}^2  = e^{\omega}\!\!\prod_{k=h-1}^\infty\!\!\exp\big\{\lambda_{d,k}g(Z_{n+h-1-k})\big\}
    \prod_{\ell=0}^{h-2}\E\big(\exp\big\{\lambda_{d,\ell}g(Z_{0})\big\}\big), \quad \mbox{for all } h > 1.
  \end{equation}
  Moreover, if\, $\E(\sigma_{n+h}^4) < \infty$, the mean square errors  of forecast for $\sigma_{n+h}^2$ and
  $X_{n+h}^2$, respectively denoted as $mse(\sigma_{n+h}^2)$ and $ mse(X_{n+h}^2)$, are given by
  \begin{equation}\label{eq:mse_sigma}
    mse(\sigma_{n+h}^2)  =
    \Bigg[\prod_{k = 0}^{h-2}\E\big(\exp\{2\lambda_{d,k}g(Z_0)\}\big) -  \prod_{\ell = 0}^{h-2}\!\Big[\E\big(\exp\big\{\lambda_{d,\ell}g(Z_0)\big\}\big)\Big]^2\Bigg]\prod_{j = h-1}^\infty\!\E\big(\exp\big\{2\lambda_{d,j}g(Z_0)\big\}\big)
  \end{equation}
  and
  \[
  mse(X_{n+h}^2) = \E(\sigma_0^4)\big[\E(Z_0 ^4) - 1\big] + mse(\sigma_{n+h}^2), \quad \mbox{for all } h > 1.
  \]
\end{thm}

\cite{LP2013} propose two $h$-step ahead predictors for the process $\{\sigma_t^2\}_{t\in\Z}$ in the context
of FIEGARCH processes.  The first one, denoted by $\check{\sigma}_{n+h}^2 $, was obtained through the relation
$\check{\sigma}_{n+h}^2 := \exp\{\hat{\ln}(\sigma_{n+h}^2)\}$, where $ \hat{\ln}(\sigma_{n+h}^2)$ is the
$h$-step ahead predictor for $\ln(\sigma_{n+h}^2)$.  The second predictor, denoted by $\tilde{\sigma}_{n+h}^2$,
was derived upon considering an order 2 Taylor's expansion of the exponential function.  The authors also
showed that $\check \sigma_{n+h}^2$ and $\tilde \sigma_{n+h}^2$ satisfy
\begin{equation}\label{eq:relation2}
  \tilde \sigma_{n+h}^2  = \left\{
    \begin{array}{ccc}
      \displaystyle  \exp\{\hat\ln(\sigma_{n+h}^2)\} = \check{\sigma}_{n+h}^2, & \mbox{if} & h =1;\vspace{0.2cm}\\
      \displaystyle \exp\{\hat\ln(\sigma_{n+h}^2)\}\bigg[ 1
      +\frac{1}{2}\sigma^2_g\sum_{k=0}^{h-2}\lambda_{d,k}^2\bigg] =  \check{\sigma}_{n+h}^2\bigg[ 1
      +\frac{1}{2}\sigma^2_g\sum_{k=0}^{h-2}\lambda_{d,k}^2\bigg], & \mbox{if } &  h > 1.
    \end{array}
  \right.
\end{equation}

With obvious identifications, the same predictors $\check \sigma_{n+h}^2$ and $\tilde
\sigma_{n+h}^2$ can be defined for SFIEGARCH processes.  However, by following the same steps as in
\cite{LP2013}, it can be shown that both $\check \sigma_{n+h}^2$ and $\tilde \sigma_{n+h}^2$ are biased
estimators for $\sigma_{n+h}^2$. On the other hand $\hat\sigma_{n+h}^2$, given in \eqref{eq:predictor}, not
only is an unbiased estimator but also it is the best predictor for $\sigma_{n+h}^2$ in terms of the mean square
error measure.

Now, to obtain the $h$-step ahead predictor for $\ln(X_{n+h}^2)$ observe that, from Definition
\ref{defn:sfiegarch} and from the i.i.d. property of $\{Z_t\}_{t\in\Z}$,
\[
\hat\ln(X_{n+h} ^2) :=  \E\big(\!\ln(X_{n+h} ^2)|\mathcal{F}_n\big)
= \E\big(\!\ln(\sigma_{n+h} ^2) |\mathcal{F}_n\big)  + \E\big(\!\ln(Z_{n+h}^2)|\mathcal{F}_n\big)
:= \hat \ln(\sigma_{n+h}^2) + \E(\ln(Z_0^2)),
\]
for all $n \in \Z$ and $h > 0$.  The expressions for $\hat{\ln}(\sigma_{n+h}^2) :=
\mathds{E}(\ln(\sigma_{n+h}^2)|\mathcal{F}_{t})$ and for the mean square errors of forecast for
$\ln(\sigma_{n+h}^2)$ and $\ln(X_{n+h}^2)$ are given in Theorem \ref{thm:forecast_log}.

\begin{thm}\label{thm:forecast_log}
  Let $\{X_t\}_{t \in \mathds{Z}}$ be an \emph{SFIEGARCH}$(p,d,q)_s$ process.  Then, for any fixed
  $n\in\mathds{Z}$, the $h$-step ahead forecast $\hat{\ln}(\sigma_{n+h}^2)$ of $\ln(\sigma_{n+h}^2)$, given
  $\mathcal{F}_{n}$, can be expressed as
  \begin{equation}\label{eq:forecast_log}
    \hat{\ln}(\sigma_{n+h}^2) = \omega + \sum_{k=0}^{\infty}\lambda_{d,k+h-1}\,g(Z_{n-k}),\
    \mbox{ for all } h>0.
  \end{equation}
  Moreover, if $h = 1$, the mean square errors of forecast for $\ln(\sigma_{n+h}^2)$ and $\ln(X_{n+h}^2)$ are both
  equal to zero and, for any $h> 1$,  are given, respectively, by
  \begin{equation}\label{eq:mse_log}
    mse(\ln(\sigma_{n+h}^2))  = \sigma_g^2\sum_{k=0}^{h-2}\lambda_{d,k}^2  \quad \mbox{and} \quad mse(\ln(X_{n+h}^2))  =     mse(\ln(\sigma_{n+h}^2)),
  \end{equation}
  where $\sigma^2_g = \mathds{E}\big([g(Z_0)]^2\big)$.
\end{thm}

\begin{remark}
  \cite{LP2013} consider the $h$-step ahead predictor for $\ln(X_{n+h} ^2)$ defined as $\check\ln(X_{n+h} ^2)
  = \hat \ln(\sigma_{n+h}^2)$, for any $n\in\Z$ and $h>0$.  This is an unbiased estimator if and
  only if $\E(\ln(Z_0^2)) = 0$.
\end{remark}

From \eqref{eq:predictor} and \eqref{eq:forecast_log}, one concludes that, $\hat\sigma_{n+1}^2 =
\check\sigma_{n+1}^2 = \tilde\sigma_{n+1}^2$ and
\begin{equation}\label{eq:relation3}
  \hat \sigma_{n+h}^2  =   \check{\sigma}_{n+h}^2 \prod_{\ell=0}^{h-2}\E\big(\exp\big\{\lambda_{d,\ell}g(Z_{0})\big\}\big) =
  \tilde{\sigma}_{n+h}^2 \bigg[ 1
  +\frac{1}{2}\sigma^2_g\sum_{k=0}^{h-2}\lambda_{d,k}^2\bigg]^{-1}\prod_{\ell=0}^{h-2}\E\big(\exp\big\{\lambda_{d,\ell}g(Z_{0})\big\}\big),
\end{equation}
for any fixed $n\in \Z$ and  $h > 1$.

Now, let $\{r_t\}_{t\in\Z}$ be the stochastic process defined by
\begin{equation}\label{eq:rt}
  r_t = \mu + \sum_{k = 0}^\infty \psi_kX_{t-k} := \mu + \psi(\B)X_t, \quad \mbox{for all } t\in\Z,
\end{equation}
where $\mu \in\R$, $\{\psi_k\}_{t\in\Z}$ is a sequence of real numbers satisfying $\sum_{k =0}^\infty \psi_k^2
< \infty$ and $\{X_t\}_{t\in\Z}$ is an SFIEGARCH$(p,d,q)_s$ with $\displaystyle \sup_{t\in\Z} \{\var(X_t^2)\}
< \infty$.  Notice that, $\sum_{k =0}^\infty \psi_k^2 < \infty$ and $\displaystyle \sup_{t\in\Z}
\{\var(X_t^2)\} < \infty$ imply $\cov(X_k, X_j) = 0$, for all $k\neq j$, and \vspace{-0.5\baselineskip}
\[
\E\bigg(\bigg|\sum_{k = m}^n \psi_kX_{t-k}\bigg|^2\bigg) \leq  \sup_{t\in\Z} \{\var(X_t^2)\} \sum_{k = m}^n\psi_{k}^2, \quad \mbox{for all } m \leq  n.
\]
Therefore,  \eqref{eq:rt} converges in $L^2$ (Cauchy convergence criterion) and hence $\{r_t\}_{t\in\Z}$ is well
defined.

\begin{example}
  Let $\phi(\cdot)$ and $\varphi(\cdot)$ be the polynomials of order $p$ and $q$, with no common roots,
  respectively defined by $\phi(z) := \sum_{k=0}^{p_1}(-\phi_k)z^k$ and $\varphi(z)
  :=\sum_{j=0}^{q_1}(-\varphi_j)z^j$, with $\phi_0 = \varphi_0 = -1$.  Given a weakly stationary
  SFIEGARCH$(p,d,q)_s$ process $\{X_t\}_{t\in\Z}$, define $\{r_t\}_{t\in\Z}$ by
  \[
  \phi(\B)( r_t - \mu) =  \varphi(\B)X_t, \quad \mbox{for all } t\in\Z.
  \]
  Observe that $\{r_t\}_{t\in\Z}$ is an ARMA$(p_1,q_1)$ process and hence it can be rewritten as in equation
  \eqref{eq:rt}, whenever $\phi(z) \neq 0$ in the closed disk $\{z: |z| \leq 1\}$, with $\{\psi_k\}_{k\in \Z}$
  uniquely defined through
  \[
  \sum_{k = 0}^\infty \psi_kz^k = \psi(z) := \frac{\varphi(z)}{\phi(z)}, \quad  |z| \leq 1.
  \]
\end{example}

Theorem \ref{thm:arma_sfiegarch} provides the $h$-step ahead forecast for the process $\{r_t^2\}_{t\in\Z}$,
with $r_t$ defined in \eqref{eq:rt}.  Similar equations can be derived if the assumption that
$\{X_t\}_{t\in\Z}$ is an SFIEGARCH$(p,d,q)_s$ is replaced by any other ARCH-type model.  This result is
applied in Section \ref{sec:analysisObserved} to compare the forecasting performance of the different models
considered in the time series analysis.

\begin{thm}\label{thm:arma_sfiegarch}
  Let $\{r_t\}_{t\in\Z}$ be defined by \eqref{eq:rt}. Then, for any fixed $n\in\ Z$, the $h$-step ahead
  forecast $\hat r_{n+h}^2$ of $r_{n+h}^2$, given $\mathcal{F}_{n}$, can be expressed as
  \begin{equation}\label{eq:rt2forecast}
    \hat r_{n+h}^2 = \mu^2  +  \sum_{k = 0}^{h-1} \psi_k^2 \hat \sigma_{n+h-k}^2 +
    \sum_{j = h}^\infty\sum_{\ell = h}^\infty \psi_j\psi_\ell X_{n+h -j}X_{n+h-\ell} +
    2\mu\sum_{i = h}^\infty\psi_iX_{n+h-i}, \quad \mbox{for any } h > 0,
  \end{equation}
  where $\hat \sigma_{n+1}^2 = \sigma_{n+1}^2$ is given by \eqref{eq:sigmat} and $\hat \sigma_{n+h}^2$ is the
  $h$-step ahead forecast of $X_{n+h}^2$ given in \eqref{eq:predictor}, for all $h>1$.
\end{thm}

\section{An Application}\label{sec:analysisObserved}

In this section we analyze the behavior of the intraday volatility of the S\&P500 US stock index log-return
time series in the period from December 13, 2004 to  October 10, 2009.  The trading hours are from
8:30 am to 3:15 pm (Chicago time) and the intraday frequency of the original index time series is 15 minutes,
which gives a total of  33993 observations (1259 days).   The fifteen-minute log-returns (see Remark
\ref{remark:return}) are aggregated to obtain a one-hour\footnote{ Up to this day, we only had access to
  the stock index time series with sampling frequency equal to 15 minutes. Therefore, a thorough study on the
  existence of microstructural noise could not be performed.  With the data set available we create a
  signature volatility plot \citep[see, for instance,][]{AEA2000} by considering only sample frequencies which
  are multiples of 15 minutes (e.g., 15, 30 and 45 minutes).  The minimum and maximum sample frequency values
  considered were, respectively, 15 minutes (which is the original sampling frequency) and 405 minutes (which
  corresponds to one-day log-returns).  Under this scenario, the average realized volatilities were all close
  to each other regardless the sampling frequencies considered. } log-return time series
$\{R_t\}_{t=1}^{8498}$.

\begin{remark}\label{remark:return}
  The fifteen-minutes  log-return time series  $\{R_t^{(15)}\}_{t=1}^n$ is defined as
  \[
  R_t^{(15)} = 100\times\ln\bigg(\frac{P_t}{P_{t-1}}\bigg), \quad \mbox{for all } t \in \{1,\cdots, n\},
  \]
  where $\{P_t\}_{t=0}^{n}$ is the index time series, with $n = 33992$ ($t = 0$ corresponds to the first  available
  observation).  The one-hour log-return time series
  $\{R_t\}_{t=1}^{8498}$ is obtained by letting $R_t = R_{\tau-3}^{(15)} + \cdots + R_\tau^{(15)}$, with $\tau =
  4t$, for all $t\in\{1,\cdots, 8498\}$.
\end{remark}

\begin{figure}[!ht]
  \centering
  \subfloat[]{ \includegraphics[width = 0.47\textwidth]{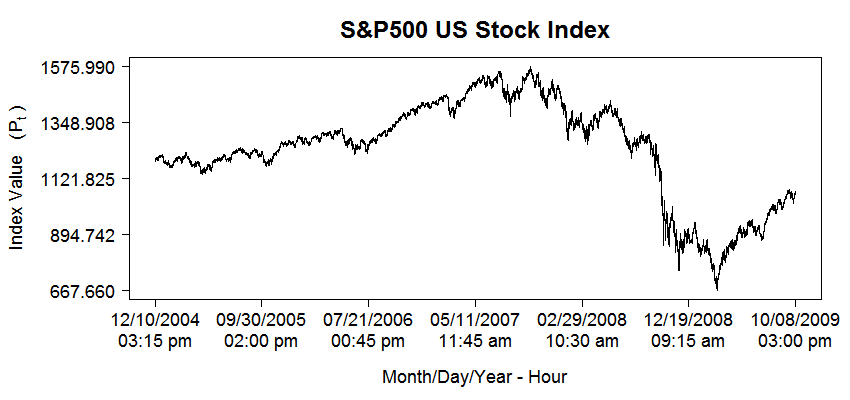} }
  \subfloat[]{  \includegraphics[width = 0.47\textwidth]{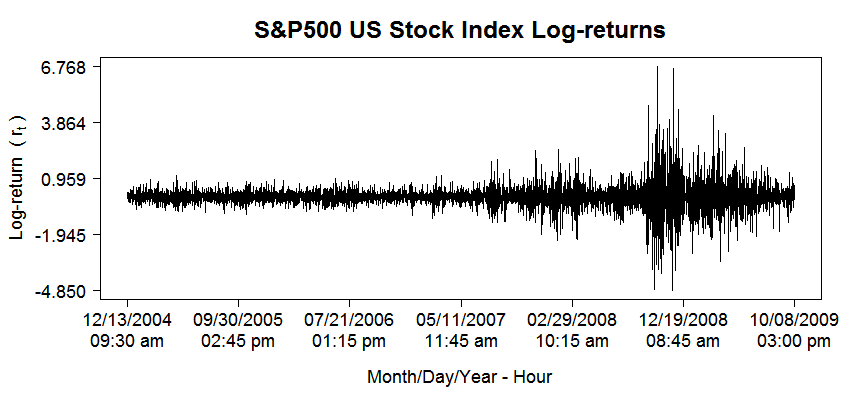}}
  \caption{ This figure considers the S\&P500 US stock index time series in the period from December 13, 2004 to
    October 10, 2009.  The intraday frequency of the index time series is 15 minutes, which gives a total of
    33993 observations (1259 days).  Panel (a) shows the original index time series. Panel  (b) presents the one-hour
    log-return time series, with $n = 8498$ observations, obtained by aggregating the fifteen-minute
    returns.}\label{fig:timeseries}
\end{figure}

Figure \ref{fig:timeseries} (a) and (b) show, respectively, the S\&P500 US stock index log-return time series
and the one-hour log-return time series $\{R_t\}_{t=1}^{8498}$, in the studied period.     Notice that, for any
$t \in \{1,\cdots, 1259\}$ the times 3:15 pm (closing time) from day $t$ and 8:30 am (opening time) from
day $t+1$ are equivalent (there is no trading between these two periods).  Therefore, there are 27 index values
and, consequently, 27 available 15-minutes log-returns  for each trading day.    It is easy to see that, by
applying the aggregation equation described in Remark \ref{remark:return}, the trading time associated to
one-hour log-returns for two consecutive days are not necessarily the same.   In particular,  the  following
holds
\begin{itemize}
\item the first available index value corresponds to December 13, 2004 8:30 am (or equivalently, December
  12, 2004 3:15 pm).   Consequently, the first one-hour log-return for day 1  corresponds to 9:30 am;

\item the last one-hour log-return for day 1 corresponds to 2:30 pm  and, since the trading day ends at 3:15
  pm, the next one-hour log-return will be the aggregation of  15-minutes log-returns for 2:30 pm, 2:45
  pm, 3:00 pm, 3:15 pm (or equivalently, 8:30 am from day 2)  and 8:45 am  from day 2.  Consequently,
  the first one-hour log-return for day 2 corresponds to 8:45 am;

\item whenever the first one-hour log-return corresponds to 9:30 am,  there are only 6 one-hour log-returns
  for the corresponding day;

\item for every four days there is  one day with  6 one-hour log-returns followed by 3 consecutive days with 7 one-hour
  log-returns. This fact may or may not induce a cyclical behavior.
  \end{itemize}

From Figure \ref{fig:timeseries} (b) it is clear that the volatility in the period from 2004 to 2007 is much
lower\footnote{This fact is already known in the literature and it is beyond the scope of this work to discuss
  possible causes for this behavior.}  than in the period from 2007 to 2009.  Given that high volatility is
more concerning and more difficult to model than low volatility  we shall discard the first 3996 observations
(592 days).  We also reserve the last 40 days of data (270 observations) to analyze the out-of-sample
forecasting performance of the fitted models.  The remaining time series corresponds to the period from March
21, 2007 to August 13, 2008 and has 4232 observations.   This time series shall be denoted by
$\{r_t\}_{t=1}^{4232}$, where $r_t := R_{t+3996}$, for all $t = 1,\cdots,4232$.

\begin{remark}
  The two highest peaks in Figure \ref{fig:timeseries} (a) correspond, respectively, to October 10, 2008 and
  November 21, 2008. The two lowest values in Figure \ref{fig:timeseries} (a) correspond, respectively, to
  October 06, 2008 and November 20, 2008.  October 10, 2008 is the day with the highest trading volume ever
  for the S\&P 500 index.  In this day, the trading volume for the SPY SPDR surpassed 871 million shares
  (see, for instance, AMEX:SPY daily prices for October 2008 from Yahoo! Finance).  On November 20, 2008
  the S\&P 500 index closed at 752.44, its lowest since early 1997.
\end{remark}

The descriptive statistics for the log-return time series $\{r_t\}_{t=1}^{4232}$ are given in Table
\ref{table:statistics}.  For comparison, this table also shows the descriptive statistics for time series
$\{R_t\}_{t=1}^{8498}$.  From Table \ref{table:statistics} one observes that both time series
$\{R_t\}_{t=1}^{8498}$ and $\{r_t\}_{t=1}^{4232}$ have mean approximately equal to zero but high skewness
values, which usually indicates a non-symmetric distribution.   However, we shall notice that, for these
time series, the high skewness values could be due to the presence of some outliers instead of
non-symmetry.  In fact, upon replacing all values higher than eight standard deviations by the sample mean
of the corresponding time series, the skewness values for $\{R_t\}_{t=1}^{8498}$ and $\{r_t\}_{t=1}^{4232}$
are, respectively, 0.0398 and -0.0018, which reinforces our claim.  Nevertheless, the possible outliers
are not removed in the analysis to be performed in the sequel.  The aim of this approach is to observe
whether the SFIEGARCH model captures or not this feature.

\begin{table}[!ht]
  \centering
  \renewcommand{\arraystretch}{1.1}
  \caption{ Descriptive statistics for  the S\&P500 one-hour log-return time series $\{R_t\}_{t=1}^{8498}$
    and for the time series $\{r_t\}_{t=1}^{4232}$, where $r_t := R_{t+3996}$, for all $t\in \Z$.  }\label{table:statistics}
  \begin{tabular}{ccccccc}
    \hline
   Period & $n$ & & mean      & st. dev. &  kurtosis & skewness \\
 \cline{1-2}\cline{4-7}
  2007 - 2009  &  4232   && -0.0078 & 0.6931 & 14.5986 &  0.4156\\
  2004 - 2009 &  8498   && -0.0013 & 0.5168  & 23.6251 &  0.4609\\
    \hline
    \multicolumn{6}{l}{\footnotesize {\bf Note:} st. dev. := standard deviation.}
  \end{tabular}
\end{table}

\begin{figure}[!ht]
  \centering
  \subfloat[]{ \includegraphics[width = 0.33\textwidth]{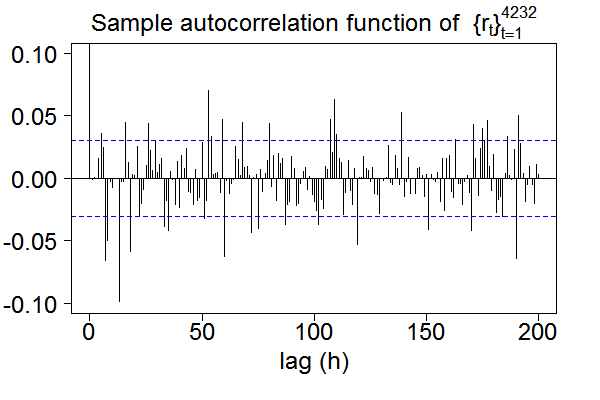}}
  \subfloat[]{ \includegraphics[width = 0.33\textwidth]{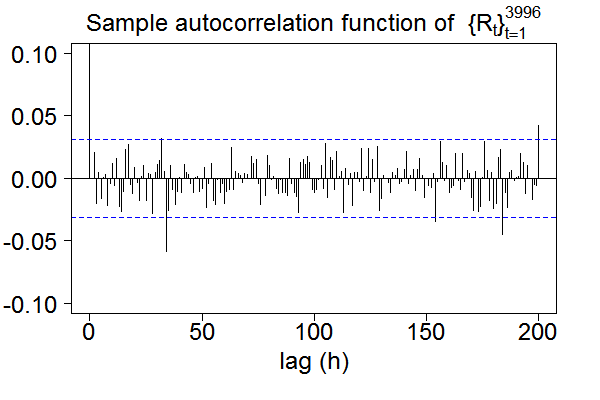} }
  \subfloat[]{\includegraphics[width = 0.33\textwidth]{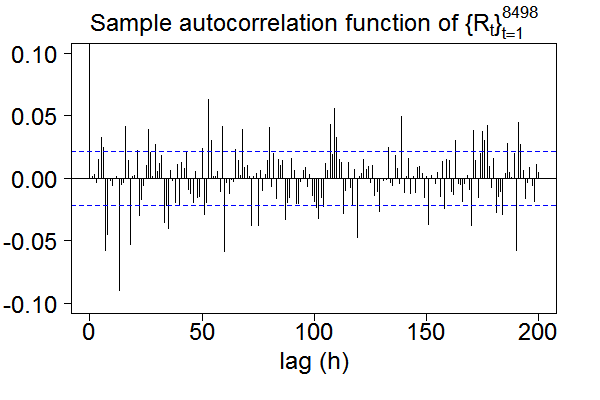}}
  \caption{ This figure shows the sample autocorrelation function for the S\&P500 US stock index one-hour
    log-return time series $\{R_t\}_{t=1}^{8498}$ and its sub-samples.  For $\{R_t\}_{t=1}^{8498}$ the time
    index $t = 1$ corresponds to December 13, 2004 at 09:30 am, $t = 3997$ corresponds to March 21, 2007 at
    09:30 am and $t = 8498$ corresponds to October 10, 2009 at 03:00 pm.  Panel (a) considers the time series
    $\{r_t\}_{t=1}^{4232}$ with $r_t := R_{t+3996}$. Panel (b) corresponds to the sub-sample $\{R_t\}_{t=1}^{3996}$.
    Panel  (c) shows the sample autocorrelation function of $\{R_t\}_{t=1}^{8498}$.  For better visualization the
    y-axis only shows the interval $[-0.10, 0.10]$.  In all graphs the dashed lines correspond to $\pm
    1.96/\sqrt{n}$, where $n$ is the sample size of the corresponding time series. }\label{fig:sp500acf}
\end{figure}

Figure \ref{fig:sp500acf} (a), (b) and (c) shows, respectively, the sample autocorrelation functions for the
log-return time series $\{r_t\}_{t=1}^{4232}$, $\{R_t\}_{t=1}^{3996}$ and $\{R_t\}_{t=1}^{8498}$.  From this
figure we observe that the log-return time series presents small (notice the scales) autocorrelations (but
significatively different from zero) for some lags $h>0$.  Upon comparing Figure \ref{fig:sp500acf} (a), (b)
and (c), one observes that while the sample autocorrelation functions for $\{r_t\}_{t=1}^{4232}$ and
$\{R_t\}_{t=1}^{8498}$ seem identical, the difference is remarkable when considering the sample
autocorrelation functions for $\{R_t\}_{t=1}^{3996}$ and $\{R_t\}_{t=1}^{8498}$.  This indicates that the
correlation in the series is mainly due to the last observations of the time series $\{R_t\}_{t=1}^{8498}$.

 From Figures \ref{fig:sp500acf} (a) and (c) one observes that the autocorrelation value with higher
magnitude is associated to the lag 13 (roughly 2 days).  Next, in order of magnitude (including the values not
reported in Figure \ref{fig:sp500acf}), are the autocorrelations associated to lags 53, 7, 263, 190, 109, 60, 18
(roughly 2 and a half days), 203, 218, 119, 139, 191 and 8.  The remaining autocorrelation values are all
smaller (in magnitude) than the one associated to lag 8 and, therefore, very close to the confidence limits
(this includes the autocorrelation values with lag higher than 200, which are not reported in Figure
\ref{fig:sp500acf}). Recall that the aggregation rule considered implies that for every 4 days there is one
day with only 6 one-hour log-returns followed by 3 days with 7 one-hour log-returns.  Moreover, notice that
53, 109, 139, 190, 191, 218 are very close to multiples of 27 (total number of observations in 4 days).
Furthermore,  while 60 and 263 differ from multiples of 27 by approximately 7, 119 and 203 differ from
approximately 13.

The facts just mentioned, indicate a short-memory cyclical behavior of length 27.  On the other hand, there is
also evidence that a single seasonal polynomial may not be enough to remove the correlation. For this reason we
shall consider a constrained\footnote{By constrained we mean that some  $\phi_i $ and $\theta_j$ will be fixed
  as zero. }  ARMA$(p_1,q_1)$ model for the log-return time series. Under this
assumption we have
\begin{equation}\label{eq:armapq}
\phi(\B) (r_t - \mu) = \varphi(\B)X_t, \quad \mbox{for all } t\in\Z,
\end{equation}
where $\mu \in \R$, $\phi(z) = \sum_{k=0}^{p_1}(-\phi_k)z^k$, $\varphi(z) = \sum_{j=0}^{q_1}(-\varphi_j)z^j$,
 with $\phi_0 = \varphi_0 = -1$, and $\{X_t\}_{t\in\Z}$ is a white noise process.    Notice that, by letting
 $p_1$ and $q_1$ be large enough, equation  \eqref{eq:armapq} also covers  the seasonal ARMA class of model,
  denoted  by SARMA$(p_1,q_1)\times(P,Q)_s$ (see Remark \ref{remark:sarima}).
 \begin{remark}\label{remark:sarima}
   For any $d, D \in \N$, let $(1-\B)^d$ and $(1-\B^s)^D$ be, respectively, the non-seasonal and seasonal
   difference operators  iterated, respectively, $d$ and $D$ times.  Let $A(z) = \sum_{k=0}^P(-A_k)z^k$,
   $a(z) = \sum_{k=0}^p(-a_k)z^k$, $M(z) = \sum_{k=0}^Q(-M_k)z^k$ and $m(z) = \sum_{k=0}^q(-m_k)z^k $ be
   polynomials, respectively, of order $P$, $p$, $Q$ and $q$, with $A_0 = a_0 = M_0 = m_0 = -1$.  A seasonal
   autoregressive integrated moving average model, denoted by SARIMA$(p,d,q)\times(P,D,Q)_s$, is defined by
   \citep[for more details and for  the definition of a SARFIMA process,  see][]{BL2009}
\begin{equation}\label{eq:sarima}
A(\B^s)a(\B)(1-\B)^d (1-\B^s)^D(Y_t - \mu) = M(\B^s)m(\B)\varepsilon_t, \quad \mbox{for all } t\in\Z,
\end{equation}
where $\mu \in \R$ and $\{\varepsilon_t\}_{t\in\Z}$ is a white noise process with zero mean and variance
$\sigma_{\varepsilon}^2$.  In particular, when $d = D = 0$, \eqref{eq:sarima} is called a
SARMA$(p,q)\times(P,Q)_s$ model.  It is immediate that by letting $d = D = 0$, $\phi(z) = A(\B^s)a(\B)$ and
$\varphi(z) = M(\B^s)m(\B)$, \eqref{eq:sarima} can be rewritten as an ARMA$(p_1, q_1)$ model, given in
\eqref{eq:armapq},  with $p_1 = P+p$ and $q_1 = Q+q$.
\end{remark}

\begin{figure}[!ht]
  \centering
  \subfloat[]{\includegraphics[width = 0.45\textwidth]{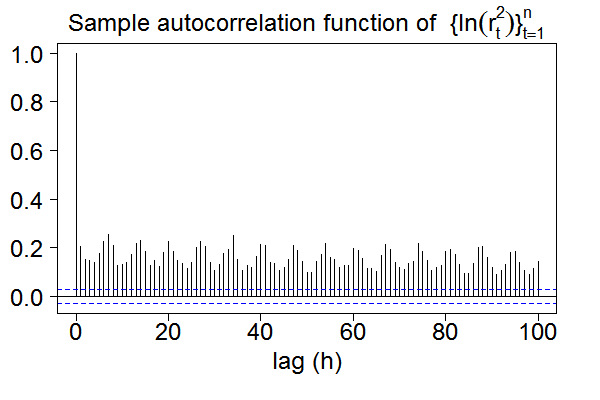}}
  \subfloat[]{\includegraphics[width = 0.45\textwidth]{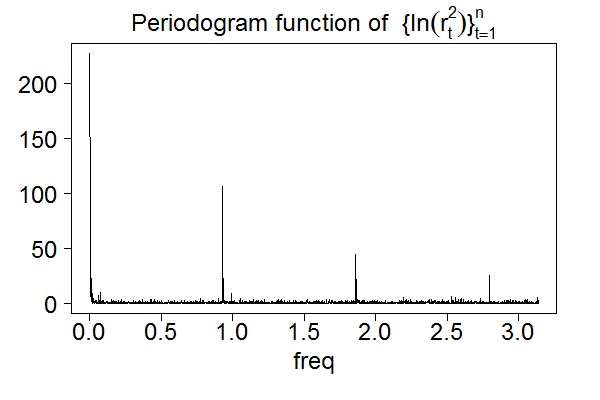}}
  \caption{ This figure considers the time series $\{\ln(r_t^2)\}_{t=1}^{4232}$, where $\{r_t\}_{t=1}^{4232}$
    is the S\&P500 US stock index log-return time series corresponding to the period from March 21, 2007 at
    09:30 am to August 13, 2009 at 03:00 pm.  Panel (a) shows the sample autocorrelation function of
    $\{\ln(r_t^2)\}_{t=1}^{4232}$. Panel (b) presents its periodogram function.}\label{fig:sp500acf2}
\end{figure}

Figure \ref{fig:sp500acf2} (a) and (b) present, respectively, the sample autocorrelation and the
periodogram functions for the time series $\{\ln(r_t ^2)\}_{t=1}^{4232}$.  We observe that
both functions  indicate long-memory and cyclical behavior, with seasonal parameter $s=7$ (one day cycle).
 To account for the long-memory cyclical behavior in the volatility,  we shall consider 
an SFIEGARCH$(p_2,d,q_2)_s$, described in Section \ref{sec:sfiegarch}. 
 To confirm the importance of including the seasonal effect associated
to long-memory\footnote{ The  FIEGARCH  model  captures non-seasonal long-memory and
   short-memory cyclical behaviors (if $p_2$ and $q_2$ are large enough).  The EGARCH model  is only
   able to describe short-memory cyclical behavior. }
we also consider a FIEGARCH$(p_2,d,q_2)$ \citep{BM1996},
and an EGARCH$(p_2,q_2)$ \citep{N1991} model. 

Recall that,  for all models just mentioned,  $\{X_t\}_{t\in\Z}$ is written as
$X_t = \sigma_tZ_t$, for any $t\in\Z$, where $\{Z_t\}_{t \in \mathds{Z}}$ is an i.i.d. sequence of random
variables with $\mathds{E}(Z_0) = 0$ and $\mbox{Var}(Z_0) = 1$. 
  For the SFIEGARCH$(p_2,d,q_2)_s$ model one has
\[
\ln(\sigma_t^2) = \omega +  \frac{\alpha(\B)}{\beta(\B)}(1-\B^s)^{-d}g(Z_{t-1}), \quad \mbox{for all } t\in\Z,
\]
with  $\alpha(\cdot)$, $\beta(\cdot)$ and $(1-\B^s)^{-d}$  be defined as in  \eqref{eq:alphabeta} and
\eqref{eq:pik},  and   $g(Z_t) = \theta Z_t +  \gamma\big[|Z_t|-\E\big(|Z_t|\big)\big]$.  The FIEGARCH and EGARCH models are particular cases of the SFIEGARCH model obtained from
  \eqref{eq:sigmat}, respectively,  when  $s=1$  and $d = 0$.
%





\begin{remark}
  By comparing the sample kurtosis values given in Table \ref{table:statistics} with the theoretical kurtosis
  values of an SFIEGARCH$(p,d,q)_s$ process (see Figures \ref{fig:kurtosisfig} and \ref{fig:kurtosisfig2}) we
  conclude that the best SFIEGARCH fit for the data more likely will have $p,q >0$.
\end{remark}

\begin{remark}
  Our first intention was to compare the performance of the SFIEGARCH model with the PLM-GARCH \citep{BEA2007,
    BEA2009}, since both models are able to describe long-memory cyclical behavior.  Analogously to the
  SFIEGARCH case, upon considering the PLM-GARCH we would also include  a HYGARCH \citep{D2004} and a GARCH
  model \citep{B1986}.  It turns out that we were not able to fit any PLM-GARCH model for which the squared  residuals
  time series $\{\hat z_t^2\}_{t=1}^n$ shows no correlation and at the same time the  positivity criteria for
  $\{\sigma_t^2\}_{t=1}^n$ would be satisfied.  The number of cases for which $\sigma_t^2 < 0$ were always
  too high to be replaced by a constant or by $|\sigma_t^2|$.
 \end{remark}

\subsection{Model Selection Procedure}

Parameter estimation is carried out by applying the so-called quasi-likelihood method.  In this method, the
estimate $\boldsymbol{\hat \eta}$ for the vector of unknown parameters $\boldsymbol{\eta }$ is the value that
maximizes
  \begin{equation*}
    \mathcal{L}_n(\boldsymbol\eta) :=-\frac{n}{2}\ln(2\pi)-\frac{1}{2}
    \sum_{t=1}^{n}\left[\ln(\sigma_t^2)+\frac{(r_t - \mu_t)^2}{\sigma_t^2}\right],
  \end{equation*}
  where $\mu_t := \mu + \sum_{k=1}^{p_1}\phi_k (r_{t-k} -\mu) + \sum_{j=1}^{q_1}\varphi_jX_{t-j}$, for all
  $t\in\Z$. The recursion starts by setting $r_t = \bar{r}$, where $\bar{r} = \frac{1}{n}\sum_{t = 1}^n r_t$
  is the sample mean of the log-return time series, and $X_t = g(Z_t) = 0$ and $X_t^2 = \sigma_X^2$, whenever
  $t\leq 0$, where $\sigma_X^2$ is the sample variance of $\{X_t\}_{t=1}^n$.

   Since we intent to compare the performance of the SFIEGARCH model with other ARCH-type models which
  incorporate or not seasonal long-memory cyclical behavior in the volatility, we shall consider a two step
  estimation procedure. First an ARMA$(p_1,q_1)$ model for the one-hour log-return time series is selected.
  The second step consists on fitting an SFIEGARCH$(p_2,d,q_2)_s$ (or any other ARCH-type) model to the
  residuals of the ARMA model.

  The information obtained through the analysis of the autocorrelation function of $\{r_t\}_{t=1}^{4232}$
  is applied to select the orders $p_1$ and $q_1$.  Given the large number of possible combinations of $p_1$ and $q_1$ we
  restrict our attention to  ARMA$(263,0)$ and ARMA$(0, 263)$ models with $\phi_i = \theta_j = 0$,
  whenever $i,j \notin \{7,8,13,18,53,60,109,119,139,190,191,203,218,263\}$.

  In order to select an SFIEGARCH$(p_2,d,q_2)_s$ 
  model for the residuals we   fix $s=7$ (estimated from the periodogram and sample autocorrelation functions) and consider different
  combinations of $p_2$ and $q_2$.  Once the right combination of $p_2$ and $q_2$ is found,  
  the same values are consider to select the FIEGARCH and the  EGARCH models. 

  The following criteria applies to both  estimation steps.

\begin{enumerate}[{\bf 1.}]
\item For any combination of $p_1,q_1$ or $p_2,q_2$, we start with the full model and remove the
  non-significant parameters (one at a time) until all p-values are smaller than 0.05.

\item  The standard deviations for the model parameters were obtained by considering the robust covariance
  matrix given by $n^{-1}H^{-1}(\boldsymbol{\hat\eta})B(\boldsymbol{\hat\eta})H^{-1}(\boldsymbol{\hat\eta})$,
  where $n^{-1}H(\boldsymbol{\hat\eta})$ and $n^{-1}B(\boldsymbol{\hat\eta})$ are, respectively, the Hessian
  and the outer product of the gradients \citep[see][]{BW1992}.

\item A model is considered to fit the data well if $\{\hat x_t\}_{t=1}^n$ (the residual of the ARMA model),
  $\{\hat z_t\}_{t=1}^n$ (the residual of the SFIEGARCH model) and $\{\hat z_t^2\}_{t=1}^n$ show no
  significant correlation. To test for correlation we consider both the Box-Pierce and Ljung-Box hypothesis
  tests (see Remark \ref{remark:lb}).

\item When more than one model satisfy the criteria in Step 3, model selection is performed based on the
  values of the log-likelihood, AIC, BIC and HQC criteria, obtained in Step 2.

\item In case two or more models present similar AIC, BIC, HQC and/or log-likelihood criteria values, we chose
  the more parsimonious one.
\end{enumerate}

\begin{remark}\label{remark:lb}
  When applying the Box-Pierce (or the Ljung-Box) hypothesis test, if the null hypothesis is rejected for
  $\{\hat x_t\}_{t=1}^n$ but it is not reject for both $\{\hat z_t\}_{t=1}^n$ and $\{\hat z_t^2\}_{t=1}^n$,
  the cumulative periodogram (also known as Kolmogorov-Smirnov hypothesis test) is  considered.  If this test does not
  reject the null hypothesis that $\{\hat x_t\}_{t=1}^n$ is a white noise process, the model is not discarded.
\end{remark}

Further residuals analysis is performed by following the same approach as \cite{HEA2004}.  The procedure
consists on employing a density transformation, as presented in \cite{DEA1998}, to test the assumption $\hat
x_t|\mathcal{F}_{t-1} \sim F_t$, for some given target distribution $F_t$.

\subsection{ Forecasting Procedure}

Once the parameters of the ARMA$(p_1,q_1)$-SFIEGARCH$(p_2,d,q_2)_s$ model are estimated, out-of-sample
forecasting is performed.  To obtain the predicted values $\hat r_{n+h}$, $\hat \sigma_{n+h}^2$ and $\hat
x_{n+h}^2$, given $\{r_t\}_{t=1}^n$, with $n = 4232$ and $h\in \{1,\cdots,270\}$, we proceed as described in
steps {\bf 1} - {\bf 9} below.  We shall denote by $\boldsymbol \eta$ the true parameters, namely,
\[
\boldsymbol \eta = (d, \theta, \gamma, \omega, \mu, \phi_1, \cdots, \phi_{p_1},\varphi_1, \cdots, \varphi_{q_1}, \alpha_1, \cdots, \alpha_{p_2}, \beta_{q_1}, \cdots, \beta_{q_2} )',
\]
and by $\boldsymbol{\hat \eta}$ the estimated values.  With  obvious identifications, the forecasting considering
the other ARCH-type models is analogous. 

\begin{enumerate}[{\bf 1.}]
\item The true parameters values $(d,\alpha_1,\cdots,\alpha_{p_2},\beta_1,\cdots,\beta_{q_2})'$ are replaced by
  the estimated ones, namely, $(\hat
  d,\hat\alpha_1,\cdots,\hat\alpha_{p_2},\hat\beta_1,\cdots,\hat\beta_{q_2})'$, and the recurrence formula
  given in Proposition \ref{prop:coefs} is used to calculate the corresponding coefficients $\{\hat
  \lambda_{d,k}\}_{k = 0}^{n+270}$.  Notice that, with obvious identifications, this recurrence formula can
  be also used to calculate the coefficients $\{\hat \psi_{d,k}\}_{k = 0}^{n+270}$ associated to $\hat\phi_1$,$
  \cdots$, $\hat\phi_{p_1}$, $\hat\varphi_{1}$,  $\cdots$, $\hat\varphi_{q_1}$.

\item Upon setting $ r_t = \hat \mu$ and $\hat x_t = g(\hat z_t) = 0$, whenever $t<0$, the time series
  $\{\hat x_t\}_{t=1}^{n}$, $\{\hat\sigma_t\}_{t=1}^{n}$ and $\{\hat z_t\}_{t=1}^{n}$ are obtained recursively
  as follows\footnote{If the pseudo-likelihood is used instead of the quasi-likelihood, the value
    $\sqrt{2/\pi}$ must be replace by the value of $\E(|Z_0|)$ associated to the distribution considered in the estimation procedure.}
  \begin{align*}
    \hat x_t & = \sum_{k=0}^{p_1}\phi_k (r_{t-k} -\hat\mu) + \sum_{j=1}^{q_1}\varphi_j\hat x_{t-j},\\
    \hat \sigma_t & =\exp\bigg\{\frac{\hat \omega}{2} +
    \frac{1}{2}\sum_{k=0}^{n-1}\hat\lambda_{d,k}\Big[\hat\theta z_{t-1-k} + \hat\gamma \Big(|\hat z_{t-1-k}| -
      \sqrt{2/\pi}\Big)\Big]\bigg\} \quad \mbox{and} \quad \hat z_t = \frac{\hat x_t}{\sigma_t},
  \end{align*}
  for all $ t\in\{1,\cdots, n\}$.  Note that, in particular, $\hat x_1 = r_1 - \hat \mu$, $\hat\sigma_1 =
  e^{\hat \omega 0.5}$ and $\hat z_1 = \hat x_1\hat \sigma_1^{-1}$.

\item Since the $h$-step ahead predictor for $X_{n+h}$ given $\mathcal{F}_n$ is zero, it is set $\hat x_{n+h}
  = 0$, for all $h>0$.

\item The $h$-step ahead forecast $\hat r_{n+h}$ for $r_{n+h}$ is given by \citep[see, for
  instance][]{BD1991}
  \begin{equation}\label{eq:arma_forecast}
    \hat r_{n+h} =  \hat \mu + \sum_{k=1}^{p_1}\phi_k (\hat r_{n+h-k} -\hat\mu) +
    \sum_{j=1}^{q_1}\varphi_j\hat x_{n+h-j},  \quad \mbox{ with } \hat r_{t} = r_t,
    \mbox{ if } t \leq n.
  \end{equation}

\item An estimate $\hat\sigma_g^2$ for $\sigma_g^2$ is obtained by replacing $\mathds{E}(|Z_0|)$ and
  $\mathds{E}(Z_0|Z_0|)$, in expression \eqref{eq:sigma2gzt}, by their respective sample estimates, that is,
  \[
  \hat \sigma^2_g = \hat\theta^2 + \hat\gamma^2 - (\hat\gamma \hat{\mu}_{|z|} )^2 + 2 \,
  \hat\theta\,\hat\gamma \bigg[\frac{1}{n}\sum_{t=1}^{n}\hat z_t|\hat z_t|\bigg], \quad \mbox{with} \quad
  \hat{\mu}_{|z|} := \frac{1}{n}\sum_{t=1}^n|\hat z_t|.
  \]

\item An estimative for $E_{\boldsymbol \eta}(\ell) := \E(\exp\{\lambda_{d,\ell}g(Z_{0})\})$, given in
  \eqref{eq:predictor}, is obtained by considering the respective sample estimator
  \[
  \hat E_{\boldsymbol{\hat \eta}}(\ell) = \frac{1}{n}\sum_{t = 1 }^n
  \exp\Big\{\hat\lambda_{d,\ell}\left[\hat\theta z_t + \hat\gamma (|\hat z_t| - \hat{\mu}_{|z|})\right]\Big\},
  \quad \mbox{for any} \quad \ell \in \{0, \cdots, h-2\}.
  \]

\item  Since $\sigma_{n+1}^2$ is $\mathcal{F}_{n}$-measurable, $\hat \sigma_{n+1}^2 = \sigma_{n+1}^2$ and it
  is  computed as in step {\bf 2}.

\item The predictor $\check \sigma_{n+h}^2 $ is obtained upon replacing the true parameter values by the
  estimated ones in \eqref{eq:forecast_log}, with the additional assumption $g(\hat z_t ) = 0$, if $t<0$.
  Then, from \eqref{eq:relation3}, $ \hat \sigma_{n+h}^2 $ and $\tilde{\sigma}_{n+h}^2$ are obtained by setting,
  respectively,
  \[
  \hat \sigma_{n+h}^2 = \check{\sigma}_{n+h}^2 \prod_{\ell=0}^{h-2}\hat E_{\boldsymbol{\hat \eta}}(\ell) \quad
  \mbox{and} \quad \tilde{\sigma}_{n+h}^2 = \tilde{\sigma}_{n+h}^2 \bigg[ 1 +\frac{1}{2}\hat
  \sigma^2_g\sum_{k=0}^{h-2}\hat\lambda_{d,k}^2\bigg]^{-1}\prod_{\ell=0}^{h-2}\hat E_{\boldsymbol{\hat
      \eta}}(\ell), \quad \mbox{for any} \quad h > 2.
  \]

\item The predictor $\hat r_{n+h}^2 $ is obtained through \eqref{eq:rt2forecast}, with the additional
  assumption $\hat x_{t} = 0$, if $t<0$.
\end{enumerate}

\begin{remark}
  In the literature, the time series $\{\hat \mu_t\}_{t = 1}^n$, given by,
  \[
  \hat\mu_t = \sum_{k=1}^{p_1}\phi_k (r_{t-k} -\hat\mu) + \sum_{j=1}^{q_1}\varphi_j\hat x_{t-j}, \quad
  \mbox{for } t\in \{1, \cdots, n\},
  \]
  is called \emph{fitted values} or \emph{in-sample forecasts}\footnote{ From \eqref{eq:arma_forecast} it is
    clear that $\hat\mu_t $ is the 1-step ahead forecast for $r_t$,  given $\mathcal{F}_{t-1}$, for any $t>
    0$. }.  Consequently, the residuals time series $\{\hat x_t\}_{t = 1}^n$ is also denoted \emph{in-sample
    errors of forecast}.  Furthermore, since $\hat z_t = \hat x_t\hat \sigma_t^{-1}$, for all
  $t\in\{1,\cdots,n\}$, the time series $\{\hat z_t\}_{t=1}^n$ is often called \emph{standardized residuals}.
\end{remark}


\subsection{ Forecasting Performance and Models Comparison}\label{sec:forecastComp}

Without loss of generality, let $\{\hat y_t\}_{t = t_0}^{t_0 -1+ n_p}$ denote either the in-sample ($t_0 = 1$
and $n_p = 4232$) or the  out-of-sample ($t_0 = 4233$ and $n_p = 270$) forecast values corresponding to the time
series $\{y_t\}_{t = t_0}^{t_0 -1+ n_p}$, where $y_t$ is either $r_t$ or $r_t^2$ and $n_p$ is the number of
predicted values.  Denote by $\mathcal{M}$ any model used to obtain $\{\hat y_t\}_{t = t_0}^{t_0 -1+ n_p}$.
The forecasting performance of model $\mathcal{M}$ is evaluated by computing the mean absolute error ($mae$),
the mean percentage error ($mpe$) and the maximum absolute error ($max_{ae}$) measures, namely,
\[
mae(\mathcal{M}) = \frac{1}{n_p}\sum_{t=1}^{n_p}|e_t|, \quad  mpe(\mathcal{M}) :=
\frac{1}{n_p}\sum_{t=1}^{n_p}\frac{|e_t|}{|y_t|} \quad \mbox{and} \quad
max_{ae}(\mathcal{M}) :=  \max_{t\in\{1,\cdots,n_p\}}\{|e_t|\}
\]
where  $e_t := y_t - \hat{y}_t$ denotes the forecasting error at step
$t$. The statistical significance of the out-of-sample forecasting performance is evaluated by using the so called
Diebold and Mariano hypothesis test \citep[see][]{DM1995}.

\begin{remark}
  The $mpe$ is an interesting measure since it considers not only the magnitude of the error (as does the
  $mae$) but also the proportion between the error and the true values so it is easier to decide whether the
  error is small or not. A drawback of the $mpe$ is that this measure is highly affected when observations are
  too close to zero.
\end{remark}

The predictive performance of model $\mathcal{M}$ is also  evaluated by measuring the quality of
the one-step ahead density forecasts \citep[see, for instance,][]{P2013}.  The measure used for this analysis
is the normalized sum of the realized predictive log-likelihood, given by
\[
S(\mathcal{M}) = \frac{1}{n_p} \sum_{t = n+1}^{n+n_p}\ln\!\big(f_{t|\mathcal{F}_{t-1}}^{\mathcal{M}}(y_t; \boldsymbol{\hat\eta})\big),
\]
where $\{y_t\}_{t=1}^{n+n_p}$ is a sample from $\{Y_t\}_{t\in\Z}$, $n$ is the size of the sample used to
estimate the parameters for model $\mathcal{M}$, $n_p$ is the number of predicted values,
$f_{t|\mathcal{F}_{t-1}}^{\mathcal{M}}(\cdot; \boldsymbol{\hat\eta})$ denotes the conditional probability
density function of $Y_t$ given $\mathcal{F}_{t-1}$ and  $\boldsymbol{\eta}$ is the parameter vector for model
$\mathcal{M}$.

Forecast efficiency regressions \citep[see][]{MZ1969} are also used to compare the quality of the volatility
forecasts among the different models fitted to the data.  The standard Mincer-Zarnowitz regression
for forecast efficiency is given by
\begin{equation}\label{eq:regression}
y_{t+h} = \gamma_0 + \gamma_1\hat y_{t+h} + \varepsilon_{t},
\end{equation}
where $y_{t+h}$ is the variable of interest and $\hat y_{t+h}$ is the $h$-step ahead forecast for $y_{t+h}$ given
$\mathcal{F}_t$.  Under the null hypothesis of forecast efficiency $\gamma_0 = 0$ and $\gamma_1 = 1$.  The
coefficients $\gamma_0$ and $\gamma_1$ in \eqref{eq:regression} are obtained by ordinary least square (OLS)
estimation.  The standard errors of the estimates are corrected for heteroskedasticity and autocorrelation by
using the HAC estimator \citep[see][]{NW1987}.  Since the forecasts $\hat y_{n+h}$ are obtained from a model
$\mathcal{M}$ for which the true parameter values are unknown, the uncertainty concerning parameter estimation
is corrected by multiplying the Newey-West standard errors by $\lambda = \sqrt{1 + n_p/n}$
\citep[see][]{WM1998}, where $n$ is the sample size used to fit the model and $n_p$ is the number of predicted
values.

Since the true volatility cannot be directly measured, the forecast efficiency regression is performed by
considering the realized volatility instead.  The ideas for this analysis were adapted from \cite{K2002}, where
 a slightly different definition\footnote{\cite{K2002} considers daily log-returns and the
  models fitted to the data do not include the ARMA regression.  So,  $r_t$ was written as $r_t = \mu + X_t$,
  where $X_t$ follows an ARCH-type model.  In this case, the author replaced $\bar y_t$ by $\hat \mu$ in the
  definition of $v_t$.  In our case, $\mu$ is replaced by $\mu_t$, which may vary according to each ARCH-type
  model associated to $X_t$.  Therefore we shall use the traditional definition of realized volatility to let
  $v_t$ be model free. } for the  ``observed volatility'' was considered.  The
definitions adopted here are given below.

\begin{defn}\label{def:realized}
Let $r_{(t-1)M + k}$ be the log-return value corresponding to the $k$-th period of day $t$, for
  $k\in\{1,\cdots,M\}$ and $t = 1,\cdots, N$, where $M$ is the number of intraday periods and $N$ is the number of
  observed days.
  \begin{enumerate}[{\bf a)}]
\item   The daily  log-return, denoted by  $r_t^{(d)}$,  is defined through $r_t^{(d)} = \sum_{k =
    1}^Mr_{(t-1)M +k}$, for all $t \in\{1, \cdots, N\}$.

    \item   The daily realized volatility,  denoted  by $v_t$,  is given by
  \[
v_t = \sum_{k = 1}^M(r_{(t-1)M + k} -  \bar r_t)^2,  \quad \mbox{where} \quad  \bar r_t := \frac{1}{M} \sum_{j =
  1}^Mr_{(t-1)M + j}, \quad \mbox{for all } t \in\{1, \cdots, N\}.
\]

  \item The log-return over the period of $h$ days,  denoted by $r_{t-1}^{(d)}[h]$, is given by
  \[
  r_{t-1}^{(d)}[h] =  \sum_{j = 0}^{h-1}r_{t+j}^{(d)} = \sum_{j = 0}^{h-1}\Bigg[\sum_{k = 1}^M r_{(t+j-1)M+k}\Bigg] =
  \sum_{k = 1}^{hM} r_{(t-1)M + k}, \quad
  \mbox{for any }  h\geq 1  \mbox{ and }  t \in\{1, \cdots, N\}.
  \]
  In particular, $r_{t-1}^{(d)}[1] = r_{t}^{(d)}$ so the log-return of period 1-day  is simply the daily log-return.

\item The realized volatility over the  period  of $h$ days,  denoted by  $v_{t-1}[h]$, is
given by
\[
v_{t-1}[h] =  \sum_{j=0}^{h-1}\Bigg[\sum_{k = 1}^M(r_{(t+j-1)M + k} -  \bar r_{t+j})^2\Bigg]  =  \sum_{j =
  0}^{h-1}v_{t+j}, \quad \mbox{for any }  h\geq 1  \mbox{ and }  t \in\{1, \cdots, N\}.
\]
In particular,  $v_{t-1}[h] = v_t$ so the realized volatility over the period of 1-day  is simply the daily
volatility.
\end{enumerate}
\end{defn}

By following the same steps as in the proof of proposition 2.2 in \cite{PL2013}, one can show that, if
$\{r_t\}_{t\in\Z}$ follows an ARMA$(p_1,q_1)$-SFIEGARCH$(p_2,d, q_2)_s$ model then the forecast for the
log-return over the period of $h$ days, given the information up to day\footnote{Note that, since we are
  considering intraday log-returns, the information up to day $t-1$ corresponds to $\mathcal{F}_{(t-1)M}$,
  where $M$ is number of the intraday periods. } $t-1$, is given by
 \begin{equation}\label{eq:return_h}
 \hat r_{t-1}^{(d)}[h] =  \sum_{j = 0}^{h-1}\hat r_{t+j}^{(d)} =   \sum_{k = 1}^{hM} \hat r_{(t-1)M + k}, \quad
  \mbox{for any }  h\geq 1  \mbox{ and }  t \in\{1, \cdots, N\}.
  \end{equation}
  where $\hat r_{(t-1)M+k} = \E(r_{(t-1)M+k}|\mathcal{F}_{(t-1)M})$ is the $k$-step ahead forecast for
  $r_{(t-1)M+k}$, obtained from the ARMA model. Moreover, the forecast for the conditional variance of
  log-return over the period of $h$, given the information up to day $t-1$, namely $\sigma_{t-1}^{2\, (d)}[h]
  = \var(r^{(d)}_{t-1}[h]|\mathcal{F}_{(t-1)M})$, is given by
  \begin{equation}\label{eq:volatility_h}
\hat  \sigma_{t-1}^{2 \, (d)}[h] = \sum_{j = 1}^{hM}\Psi_j\hat\sigma_{(t-1)M  + j}^2, \quad \mbox{with }
\Psi_j := \Bigg[\sum_{k = 0}^{hM -j}\psi_k\Bigg]^2,
\end{equation}
where $\hat\sigma_{(t-1)M + k}^2 = \E(X_{(t-1)M+k}^2|\mathcal{F}_{(t-1)M})$ is the $k$-step ahead forecast for
$X_{(t-1)M+k}^2$, obtained from the SFIEGARCH model.  Equivalent result is derived upon replacing the
SFIEGARCH by any other ARCH-type model.

\begin{remark}
  Notice that Definition \ref{def:realized} and expressions \eqref{eq:return_h} and \eqref{eq:volatility_h}
  assume $M$ constant.  When $M$ varies over time, similar results are derived upon making the following
  adjustments: replace $(t-1)M$ by the number of intraday log-returns available up to day $(t-1)$ (included);
  replace $hM$ by the number of intraday log-returns over the period from $t$ to $t+h$; replace $\sum_{k=1}^M$
  and $1/M$, respectively, by $\sum_{k=1}^{M_t}$ (in Definition \ref{def:realized} c) and d), by  $\sum_{k=1}^{M_{t+j}}$) and $1/M_t$,
  where $M_t$ is the number of intraday log-returns available for day $t$.
\end{remark}

 \begin{remark}
 For any $h\geq 1$ fixed, the forecast efficiency regression  is obtained upon replacing, in
 \eqref{eq:regression},  $y_{t+1}$ and $\hat y_{t+1}$, respectively,  by
 $v_t[h]$ and $\hat \sigma_{t}^{2 (d)}[h]$,  for $t = n,\cdots, N - h$, where $N$ is the sample size of the
 log-return time series and $n$ is the size of the sample used to fit the model.
 \end{remark}

\subsection{ Results}

Constrained ARMA$(263,0)$ and  ARMA$(0, 263)$ models,  with $\phi_i = \theta_j = 0$,
  whenever $i,j \notin \{7$, $8$, $13$, $18$, $53$, $60$, $109$, $119$, $139$, $190$, $191$, $203$, $218$, $263\}$, where analyzed.  It turns out that
  several initial considered  parameters where not significant and were removed from the models.  The most
  parsimonious found model  is given by (the number in parenthesis is the robust standard error)
  \begin{align*}
    r_t & = X_t - \underset{(0.0128)}{0.0391}X_{t-7} - \underset{(0.0169)}{0.0834}X_{t-13} +
    \underset{(0.0204)}{0.0803}X_{t-53} + \underset{(0.0215)}{0.0748}X_{t-109} -
    \underset{(0.0093)}{0.0433}X_{t-190} \\
    & \quad \quad + \underset{(0.0221)}{0.0742}X_{t-203}  + \underset{(0.0167)}{0.0680}X_{t-218}  -
    \underset{(0.0111)}{0.0421}X_{t-263}\\
    & := X_t + \mu_t,
  \end{align*}
  for all $t \in \{1, \cdots, 4232\}$, with $X_t = 0$, if $t\leq 0$.

\begin{figure}[!ht]
  \centering
  \subfloat[]{\includegraphics[width = 0.49\textwidth]{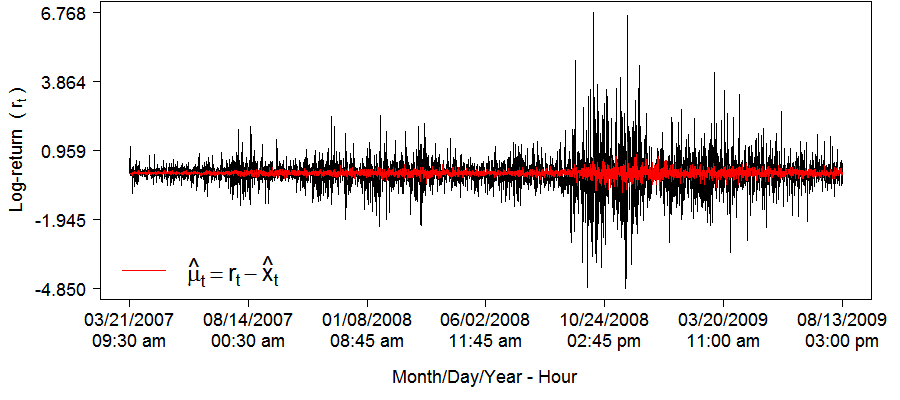}}
  \subfloat[]{\includegraphics[width = 0.49\textwidth]{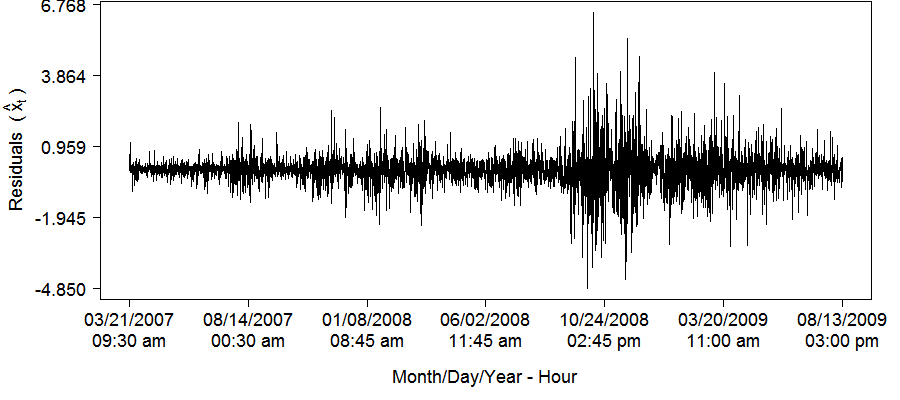}}\\
  \caption{This figure shows the S\&P500 US stock index log-return time series observed in the period from
    from  March 21, 2007 at 9:30 am      to  August 13, 2009 at 03:00 pm and the corresponding fitted values
    and residuals obtained form the constrained ARMA model.  Panel (a) gives the observed time series
    $\{r_t\}_{t=1}^{4232}$ (in black) and the fitted values $\{\hat{\mu}_t\}_{t=1}^{4232}$ (in red).
    Panel  (b) shows  the residuals  time series $\{\hat x_t\}_{t=1}^{4232}$.  }\label{fig:finalModelARMA}
\end{figure}

  The observed time series
  $\{r_t\}_{t=1}^{4232}$ and the corresponding fitted values $\{\hat \mu_t\}_{t=1}^n$, obtained from the ARMA
  model,  are given in Figure \ref{fig:finalModelARMA} (a). Figure \ref{fig:finalModelARMA} (b)
  shows the residuals time series $\{\hat x_t\}_{t=1}^{4232}$.   The p-values for the Box-Pierce and
  Ljung-Box test statistic for $\{\hat x_t\}_{t=1}^n$, that is, the residuals of the ARMA model, were smaller
  than 0.05 for any lag higher than 15.  On the other hand, upon applying the cumulative periodogram test, the
  null hypothesis that $\{\hat x_t\}_{t=1}^n$ is white noise process, was not rejected (the cumulative
  periodogram figure was omitted to save space and may be obtained from the authors upon request).  As
  expected, both the Box-Pierce (or Ljung-box) and the cumulative periodogram tests reject the null hypothesis
  that $\{\hat x_t^2\}_{t=1}^n$ is a white noise process.

 \begin{table}[!ht]
   \centering
   \caption{Estimated values and the corresponding robust standard errors  (in parenthesis) for the parameters of
     the SFIEGARCH, FIEGARCH and  EGARCH 
     models fitted to the
     S\&P500 US stock index log-return time series    observed in  the period from  March 21, 2007 at 9:30 am
     to  August 13, 2009 at 03:00 pm.  This table also presents the corresponding log-likelihood, AIC, BIC and
     HQC    criteria values.  } \label{table:finalModel}
 \begin{tabular}{ccrrr}
   \hline
   Parameter && \multicolumn{1}{c}{SFIEGARCH} & \multicolumn{1}{c}{FIEGARCH} &  \multicolumn{1}{c}{EGARCH} \\ 
   \cline{1-1}   \cline{3-5}
      $d$ && 0.4532 (0.0104) &  0.4529 (0.0057) & \multicolumn{1}{c}{-} \\
   $\omega$ && -1.1810 (0.0129) &   -1.6972 (0.0049) &  -0.8601 (0.0058) \\
   $\theta$ && -0.0820 (0.0115) &    -0.1100 (0.0017) &  -0.0954 (0.0016) \\
   $\gamma$ && 0.2127 (0.0213) & 0.2393 (0.0053)  &  0.2197 (0.0025) \\
   $\alpha_1$ &&  0.1655 (0.0040)  &  0.0657 (0.0065) &  -0.2718  (0.0015)\\
   $\alpha_2$ && 0.1963 (0.0263)  & 0.3021 (0.0088)&   0.3761  (0.0019)\\
   $\alpha_3$ && 0.1821   (0.0045) & 0.0151 (0.0095)&   -0.2503  (0.0021)\\
   $\alpha_4$ &&-0.3095   (0.0048) & -0.4206 (0.0079)&   -0.4603 (0.0039)\\
   $\alpha_5$ && -0.4115  (0.0094)  & -0.5962 (0.0136)&   -0.5355 (0.0038)\\
   $\alpha_6$ && -0.3139  (0.0112) & 0.4580 (0.0023)& 0.9040  (0.0019) \\
   $\beta_1$ && 0.3151   (0.0095)& 0.1378   (0.0039)& 0.3271   (0.0013)\\
   $\beta_2$ && -0.0668 (0.0174)&  0.1231  (0.0025)&  0.0779  (0.0016)\\
   $\beta_3$ && 0.2092 (0.0198)& -0.0081   (0.0030)& 0.1545   (0.0025)\\
   $\beta_4$ && -0.3621  (0.0079)& -0.2854   (0.0042)&   -0.4983 (0.0014)\\
   $\beta_5$ && 0.0785  (0.0124)& -0.1842  (0.0020)& 0.1895  (0.0023) \\
   $\beta_6$ &&0.6534  (0.0115)&  0.8896  (0.0029)& 0.7389   (0.0026)\\

\hline
   log-likelihood && -3053.9066 & -2936.1820& -2878.2630\\
   AIC && \phantom{-}6139.8131 &\phantom{-}5904.3641& \phantom{-}5786.5260\\
   BIC && \phantom{-}6241.4201 &\phantom{-}6005.9709& \phantom{-}5881.7824\\
   HQC && \phantom{-}6175.7272 &\phantom{-}5940.2781& \phantom{-}5820.1954\\
   \hline
   \end{tabular}
 \end{table}

  The estimated values and the corresponding standard errors for the parameters of considered ARCH-type models
  are given in Table \ref{table:finalModel}.    For any model in Table \ref{table:finalModel},  the p-values
  for both  the  Box-Pierce and Ljung-Box test statistics  corresponding to the time series
  $\{\hat{z}_t\}_{t=1}^{4232}$,  that is, the residual of the ARCH-type model,   were always  higher than
  0.05,  for any lag $h > 0$.   The same applies to the time series  $\{\hat{z}_t^2\}_{t=1}^{4232}$.    The
  histogram (or the kernel density) and the QQ-Plot (both omitted here to save space but available from the
  authors upon request) indicate that, although symmetric,  the distribution of   $\{\hat{z}_t\}_{t=1}^{4232}$
  is not Gaussian.

  \begin{remark}
    The SFIEGARCH model given in Table \ref{table:finalModel} is the most parsimonious model such that
    $\beta(\cdot)$ has no roots in the closed disk $\{z : |z| \leq 1\}$.  The p-value for $\alpha_3$, in the
    FIEGARCH model, is 0.1141.  On the other hand, the model fitted without this coefficient does not lead to
    uncorrelated residuals.  The polynomial $\beta(\cdot)$ associated to the FIEGARCH model has two roots with
    absolute value 1.0023 and two other roots with absolute value 1.0024.  Therefore, very close to the unit
    circle.  Analogously, the polynomial $\beta(\cdot)$ associated to the EGARCH model has two roots with
    absolute value 1.0024 and another one with absolute value 1.0001.  Despite this fact, the EGARCH model
    presents sligltly better performance in terms of  log-likelihood, AIC, BIC and HQC criteria.
  \end{remark}

 To apply the density transform procedure \citep[for details,  see][]{HEA2004, DEA1998} the  GED distribution,
 with different values for the tail-thickness parameter $\nu$.   Under  this  scenario,  the null hypothesis  to  be tested is
     \[
  H_0 :   \hat{x}_t| \mathcal{F}_{t-1} \sim  \mbox{GED}(\nu, 0,  \sigma_t) \quad \mbox{or, equivalently,
  }\quad
  H_0 :   F_t^{-1}(\hat{x}_t) \sim  \mathcal{U}(0,1),
     \]
     where GED$(\nu, 0, \sigma_t)$ denotes the generalized error distribution with tail-thickness parameter
     $\nu$, mean zero and standard deviation $\sigma_t$,  and $F_t(\cdot)$ is the corresponding cumulative
     distribution function.   In  particular, when $\nu = 2$, we have the Gaussin distribution. The time series $\{x_t\}_{t=1}^{4232}$
     corresponds to the residuals of the ARMA model fitted to the one-hour log-return time series
     $\{r_t\}_{t=1}^{4232}$ and $\{\hat\sigma_t^2\}_{t=1}^{4232}$ denotes the conditional variance of the
     log-returns, obtained from the SFIEGARCH, FIEGARCH or from the  EGARCH model. 
     Table \ref{table:dtl} reports the results for the Kolmogorov-Smirnov (K-S) hypothesis test used to compare
     the sample $\{F^{-1}(\hat{x}_t)\}_{t=1}^{4232}$ with the uniform distribution.

\begin{table}[!ht]
\centering
   \caption{Results  for the Kolmogorov-Smirnov (K-S) hypothesis test used to compare  the  sample
     $\{F_t^{-1}(\hat{x}_t)\}_{t=1}^{4232}$ with the uniform distribution.   The values reported
     are the p-value for the K-S test statistic.   The  null hypothesis  considered is
     $H_0 : \hat{x}_t| \mathcal{F}_{t-1} \sim  \mbox{GED}(\nu, 0,  \sigma_t)$, for different values of $\nu$.
     The time series  $\{x_t\}_{t=1}^{4232}$ corresponds to the residuals of the ARMA model fitted to the one-hour log-return
     time series $\{r_t\}_{t=1}^{4232}$ and    $\{\hat\sigma_t^2\}_{t=1}^{4232}$ denotes the conditional variance
     of the log-returns, obtained  from the   SFIEGARCH,  FIEGARCH or from the  EGARCH model. } \label{table:dtl}

 \begin{tabular}{ccccc}
   \hline
   $\nu$ && SFIEGARCH & FIEGARCH &  EGARCH\\ 
      \cline{1-1}   \cline{3-5}
   1.45 && 0.56 & 0.03 & 0.25 \\
   1.50 && 0.45 & 0.09 & 0.19 \\
   1.55 && 0.36 & 0.21 & 0.14 \\
   1.70 && 0.11  & 0.28 & 0.04\\
   2.00 && 0.00 & 0.00 & 0.00\\
   \hline
   \end{tabular}
  \end{table}

  Table \ref{table:dtl} confirms the results obtained with the QQ-Plot, that is, $\{\hat z_t\}_{t=1}^n$ does
  have Gaussian distribution.  This table also indicates that the assumption that $\{\hat z_t\}_{t=1}^n$
  follows a GED($\nu$) distribution holds for more than one value of $\nu$.  The next step in this analysis
  would be to replace the QMLE by the log-likelihood estimation procedure, using the GED distribution, and
  estimate $\nu$ alongside with the other parameters.  The information on the  $\nu$ parameter could then be
  incorporated in the forecasting equation to see whether forecast efficiency improves or not.  This analysis
  shall be performed in a future work.

  The mean absolute error ($mae$), the mean percentage error ($mpe$) and the maximum absolute error
  ($max_{ae}$) measures for the selected models are reported in Tables \ref{table:mae} and
  \ref{table:mae2}.   For the in-sample analysis, the $mae$, $mpe$ and $max_{ae}$ values were obtained by
  letting $e_t = r_t - \hat \mu_t$, for all $t\in \{1, \cdots, 4232\}$ (see Section \ref{sec:forecastComp}).
  For the out-of-sample comparison, we consider not only the forecasts for $r_{t+h}$ but also for
  $r_{t+h}^2$.  The  out-of-sample  $mae$, $mpe$ and $max_{ae}$ values were obtained by letting
  $e_{t+h} = r_{t+h} - \hat r_{t+h}$, for $t = 2432$ (fixed) and $h = 1, \cdots, 270$;   by setting $h = 1$
  fixed and letting $e_{t+1} = r_{t+1} - \hat r_{t+1}$ for $t = 4232, \cdots, 4501$; by letting $e_{t+h} = r_{t+h}^2
  - \hat r_{t+h}^2$, for $t = 2432$ (fixed) and $h = 1, \cdots, 270$; and also by letting $h = 1$ fixed and
  considering $e_{t+1} = r_{t+1}^2 - \hat r_{t+1}^2$, for $t = 4232, \cdots, 4501$.  For any ARCH-type model,
  $\hat r_{t+h}^2$, for any $h \geq 1$ and $t\in \Z$, was obtained according to Theorem
  \ref{thm:arma_sfiegarch}.   While   Table \ref{table:mae} reports the $mae$, $mpe$ and $max_{ae}$ associated to the in-sample and out-of-sample
forecasts for $r_{t+h}$, Table \ref{table:mae2} gives the values associated to $r_{t+h}^2$.
Notice that, since the ARMA model was selected independently of the ARCH-type   models, all  values in
Table \ref{table:mae}  do not depend   on the model for the conditional variance
$\sigma_t^2$.

\begin{table}[!ht]
  \centering
  \caption{The mean absolute error ($mae$),  the mean percentage error ($mpe$) and the maximum absolute error
    ($max_{ae}$) measures for the  out-of-sample forecasts for $r_{t+h}$.  Case 1 are the values corresponding
    to the in-sample forecasts, with  $e_t = r_t  - \hat \mu_t $, for  $t\in \{1, \cdots, 4232\}$.  Case 2 and
    3 correspond to the out-of-sample forecasts. For Case 2,  $e_{t+h} = r_{t+h}  - \hat r_{t+h} $, for  $t =  4232$ and
     $h\in\{1, \cdots, 270\}$. Case 3 assumes   $e_{t+1} = r_{t+1}  - \hat r_{t+1} $, for  $t\in \{4232,\cdots, 4501\}$.}\label{table:mae}
  \begin{tabular}{ccrrr}
      \hline
Measure && \multicolumn{1}{c}{Case 1} &  \multicolumn{1}{c}{Case 2} &  \multicolumn{1}{c}{Case 3} \\
  \cline{1-1} \cline{3-5}
  $mae$ && 0.4309 & 0.2611 & 0.1525 \\
  $mape$ && 1.6291  &  1.3480 & 1.7712 \\
  $max_{ae}$ && 6.4438  &  1.7892 & 0.9703 \\
  \hline
\end{tabular}
\end{table}

\begin{table}[!ht]
  \caption{The mean absolute error ($mae$),  the mean percentage error ($mpe$) and the maximum absolute error
    ($max_{ae}$) measures for the  out-of-sample forecasts for $r_{t+h}^2$. For each model the left column
    considers  $e_t = r_{t+h}^2  - \hat r_{t+h}^2 $, for  $t =  4232$ and
$h\in\{1, \cdots, 270\}$,  and the right column assumes  $e_{t+1} = r_{t+1}^2  - \hat r_{t+1}^2 $, for  $t\in \{4232,\cdots, 4501\}$.}\label{table:mae2}
  \begin{tabular*}{1\textwidth}{@{\extracolsep{\fill}}ccrrcrrcrr}
      \hline
Measure && \multicolumn{2}{c}{SFIEGARCH} &&  \multicolumn{2}{c}{FIEGARCH} &&   \multicolumn{2}{c}{EGARCH} \\
\cline{1-1} \cline{3-4} \cline{6-7} \cline{9-10}
  $mae$ && 0.2845 & 0.1865  && 0.2829 & 0.1724 &&  0.2272 & 0.1687 \\
  $mape$ &&  970.1824 & 531.1481 && 684.6831 & 409.1600  && 417.9443 & 372.8434\\
  $max_{ae}$ && 2.9271 & 2.7457 && 2.9426 & 2.9844 &&  2.9597 & 3.0164  \\
   \hline
   \end{tabular*}
  \end{table}

\begin{figure}[!ht]
    \centering
      \subfloat[]{\includegraphics[width = 0.49\textwidth]{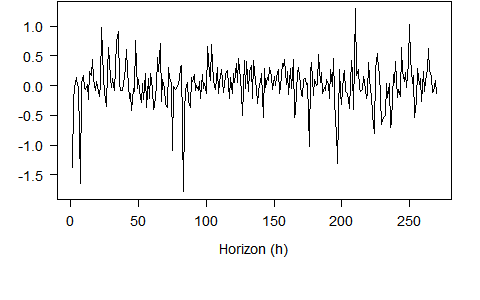}}
      \subfloat[]{\includegraphics[width = 0.49\textwidth]{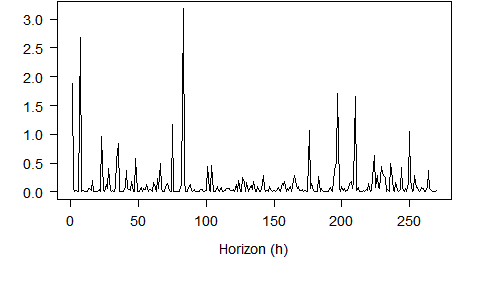}}
    \caption{Sample from the one-hour log-return time series corresponding to the period from August 14, 2009
      at 09:15 am to October 8, 2009 at 03:00 pm.   Panel (a) gives the log-returns time series
      $\{r_{t+h}\}_{h = 1}^{270}$. Panel (b)       gives the squared log-returns $\{r_{t+h}^2\}_{h =
        1}^{270}$.  The time  index  $t= 4232$ corresponds to August 13, 2009 at 03:00 pm.   }\label{fig:observed}
  \end{figure}

\begin{figure}[!ht]
        \centering
        \subfloat[]{\includegraphics[width = 0.33\textwidth]{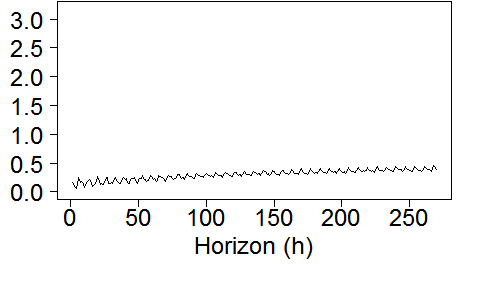}}
        \subfloat[]{\includegraphics[width = 0.33\textwidth]{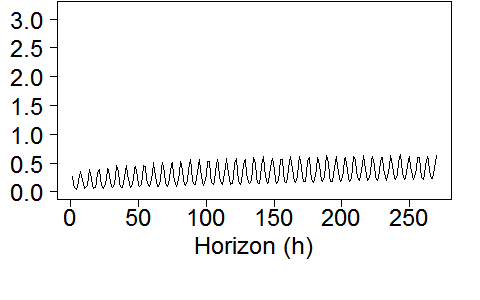}}
        \subfloat[]{\includegraphics[width = 0.33\textwidth]{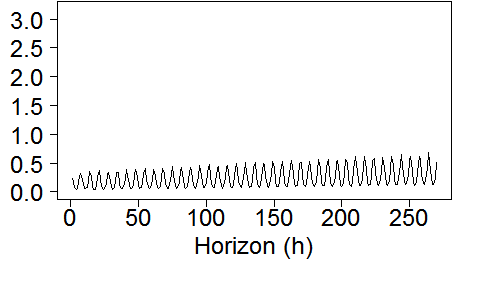}}   \\
        \subfloat[]{ \includegraphics[width = 0.33\textwidth]{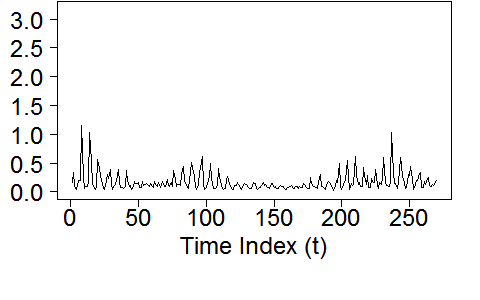}}
        \subfloat[]{\includegraphics[width = 0.33\textwidth]{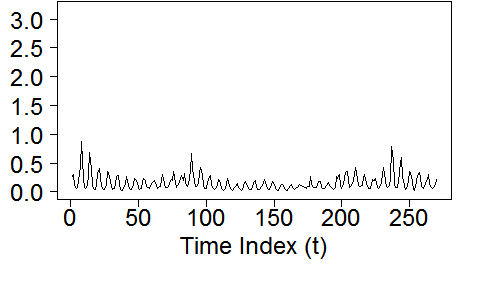}}
        \subfloat[]{\includegraphics[width = 0.33\textwidth]{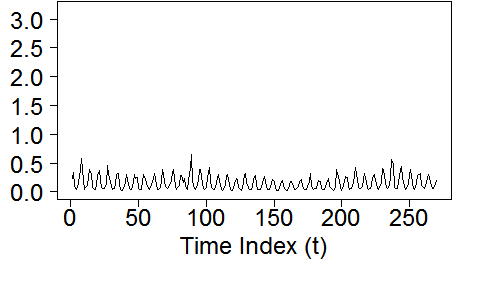}}

    \caption{Out-of-sample forecasts for the sample from the one-hour log-return time series corresponding to the
      period from August 14, 2009 at 09:15 am to October 8, 2009 at 03:00 pm.  Panels (a), (b) and (c) give
      the $h$-step ahead forecasts $\{\hat r_{t+h}^2\}_{h = 1}^{270}$, with forecasting origin $t = 2432$, for
      the SFIEGARCH, FIEGARCH and EGARCH models, respectively.  Panels (d), (e) and (b) give the one-step
      ahead forecast $\{\hat r_{t+1}^2\}_{t = 4232}^{4501}$, obtained by considering, respectively, the
      SFIEGARCH, FIEGARCH and EGARCH models fitted to the sample $\{r_t\}_{t=1}^{4232}$.
    }\label{fig:forecasts}
  \end{figure}

Table \ref{table:mae2} indicates that the SFIEGARCH model presents the best performance, among all models,
only in terms of $max_{ae}$.  Figures \ref{fig:observed} and \ref{fig:forecasts} show, respectively, the
observed values $r_{t+h}$, $r_{t+h}^2$ and the fitted ones, obtained from the ARCH-type models, for $h = 1,
\cdots, 270$.  These figures help to explain the reason why the SFIEGARCH model have higher $mae$ and
$mape$ than the other ARCH-type models.

From Figure \ref{fig:observed} it is clear that $r_{t+h}^2$ is very close to zero, for several values of
$h\in\{1, \cdots, 270\}$.  From Figure \ref{fig:forecasts} (first row, from top to bottom) it is clear that the $h$-step ahead forecasts
obtained from the SFIEGARCH model converge to a fixed value, which is expected since $\hat \sigma_{t+h}^2$
converges to the unconditional variance as $h$ goes to infinity.  From Figure \ref{fig:forecasts} (first row, from top to bottom)
is also evident that the forecasts for the FIEGARCH and EGARCH models do not converge at all.  On the
contrary, the forecast values for these two models seem to oscillate around a curve which converges to the
same value as the forecasts from the SFIEGARCH model.  Moreover, the amplitude of these oscillations increases
over time. This should be expected since both the FIEGARCH and the EGARCH models are such that the polynomial
$\beta(\cdot)$ has at least one unit root (or a root close enough to the unit root).

Figure \ref{fig:forecasts} (second row, from top to bottom)  also shows that the one-step ahead forecasts for the SFIEGARCH model
were able to capture the peaks in the square log-returns much better than the other two models.  Although this
is a good feature of the SFIEGARCH model it also makes the $mae$ and the $mape$ values increase since  the forecast
values tend to increase in the region around the peaks while the observed time series shows several values
close to zero in the same region.  Also, notice that,  since  the forecasts for $r_{t+h}^2$ are always positive, it is evident
that a model which oscilates as the FIEGARCH and the EGARCH do will provide $h$-step ahead forecasts close to
zero more often than a model for which the forecast value converges to a non-zero constant.

\begin{remark}
  We apply the Diebold and Mariano hypothesis test \citep[see][]{DM1995} to verify the statistical significance of
  the out-of-sample forecasting performance.  We consider the absolute error as loss function so the
  loss-differential series is given by $\{\mathrm{d}_t\}$, with $\mathrm{d}_t = |y_{t_0+t+1} - \hat y_{t_0+t+1}|
  - |y_{t_o+t+1}|$, for  $t_0 = 2432$ and  $t = 0, \cdots, 269$. The variable  $y_{t_0+t+1}$ denotes either the log-returns $r_{t_0+t+1}$ or
  the squared log-returns $r_{t_0+t+1}^2$ and $\hat y_{t_0+t+1}$ are the corresponding one-step ahead forecasts
  obtained from model $\mathcal{M}$.  For all models the p-value for the test statistic was smaller than 0.0002.
  Therefore, the null hypothesis that $\hat y_{t_0+t+1} = 0$, for all $t > 0$, was always rejected.
  \end{remark}

  Table \ref{table:predloglike} reports the values of the normalized sum of the realized predictive likelihood
  for each ARCH-type models.  The statistic was obtained by considering the GED($\nu$) probability density
  function, for different values of $\nu$. In particular, for $\nu = 2$ we have the Gaussian case.  The
  results in \ref{table:predloglike} show that $\nu = 1.45$ provides slightly better density forecasts for all
  three models, compared to the other values of $\nu < 2$. Moreover, density forecasts from the EGARCH model
  are slightly better than for the other two models.

  \begin{table}[!ht]
    \centering
    \caption{Normalized sum values of the realized predictive likelihood  for each ARCH-type models considering the
      GED($\nu$) probability density function for different values of $\nu$.}\label{table:predloglike}
 \begin{tabular}{ccccc}
   \hline
   $\nu$ && SFIEGARCH & FIEGARCH &  EGARCH \\ 
   \cline{1-1} \cline{3-5}
 1.45 &&1.0050   & 1.0959   & 1.1331 \\
 1.50 &&0.9973   & 1.0893   & 1.1269  \\
 1.55 &&0.9899   & 1.0829   & 1.1209  \\
 1.70 &&0.9688   & 1.0649   & 1.1040  \\
 2.00 &&0.9318   & 1.0338   & 1.0749  \\
 \hline
\end{tabular}
\end{table}

Table \ref{table:fore_eff} shows the estimated values of the coefficients $\gamma_0$ and $\gamma_1$ in
\eqref{eq:regression}.  From Table \ref{table:fore_eff} one concludes that, in all cases, the null hypothesis
of forecast efficiency ($\gamma_0 = 0$ and $\gamma_1 = 1$) is rejected.  Table \ref{table:fore_eff} also
indicates that the three models have a very similar performance.  On the other hand, as we mentioned earlier,
the log-return time series has several values very close to zero. Moreover, the volatility forecasts obtained from the
SFIEGARCH model converge to a constant  while the forecasts obtained from the FIEGARCH and EGARCH models
show an oscillating behavior leading to forecasts close to zero more often than the SFIEGARCH model.

\begin{table}[!ht]
  \centering
  \renewcommand{\tabcolsep}{2pt}
  \caption{Estimated values of the coefficients $\gamma_0$ and $\gamma_1$ defined in   the  forecast efficiency
    regression. The number in parenthesis are the standard errors  corrected for heteroskedasticity and
    autocorrelation (by using the HAC estimator)  and uncertainty concerning parameter estimation (upon
    multiplying by  $\lambda = \sqrt{1 + n_p /n}$, where $n$ is the sample size used to fit the model and $n_p$ is the number of predicted values).} \label{table:fore_eff}
  {\footnotesize
    \begin{tabular*}{1\textwidth}{@{\extracolsep{\fill}}ccrrcrrcrr}
      \hline
      \multirow{2}{*}{$h$} && \multicolumn{2}{c}{SFIEGARCH} && \multicolumn{2}{c}{FIEGARCH}  && \multicolumn{2}{c}{EGARCH} \\
      \cline{3-4} \cline{6-7} \cline{9-10}
      && \multicolumn{1}{c}{$\gamma_0$} & \multicolumn{1}{c}{$\gamma_1$} && \multicolumn{1}{c}{$\gamma_0$} &
      \multicolumn{1}{c}{$\gamma_1$}  && \multicolumn{1}{c}{$\gamma_0$} & \multicolumn{1}{c}{$\gamma_1$} \\
      \cline{1-1} \cline{3-4} \cline{6-7} \cline{9-10}
1 && 0.8447 (0.2668)  & -0.1087 (0.2043)  && 0.8018 (0.2453)  & -0.0789 (0.1796)  && 0.7636 (0.4685)  & -0.0458 (0.4522) \\
2 && 1.5941 (0.6158)  & -0.0942 (0.1874)  && 1.5559 (0.5474)  & -0.0840 (0.1799)  && 1.8090 (0.9887)  & -0.2285 (0.4327) \\
3 && 2.2665 (1.1055)  & -0.0657 (0.1933)  && 2.2181 (0.9726)  & -0.0538 (0.1758)  && 2.6404 (1.8851)  & -0.2181 (0.5425) \\
4 && 3.1568 (1.3799)  & -0.0906 (0.2194)  && 3.1464 (1.1611)  & -0.0918 (0.1865)  && 3.8928 (2.4312)  & -0.3186 (0.5390) \\
5 && 4.0411 (2.3326)  & -0.1074 (0.2713)  && 4.0263 (2.2802)  & -0.1077 (0.2753)  && 5.0712 (4.0225)  & -0.3666 (0.7009) \\
6 && 5.0353 (2.6067)  & -0.1345 (0.2473)  && 5.0360 (2.3323)  & -0.1374 (0.2297)  && 6.3321 (3.4384)  & -0.4120 (0.4943) \\
7 && 6.0834 (3.4448)  & -0.1598 (0.2711)  && 6.1792 (3.4706)  & -0.1743 (0.2811)  && 7.7567 (4.0976)  & -0.4684 (0.5108) \\
8 && 7.2811 (3.2498)  & -0.1938 (0.2080)  && 7.3698 (3.2783)  & -0.2050 (0.2247)  && 9.1453 (3.1329)  & -0.5013 (0.3597) \\
9 && 8.7566 (4.5220)  & -0.2508 (0.2554)  && 8.7339 (3.9764)  & -0.2485 (0.2409)  && 10.4351 (3.0866)  & -0.5120 (0.3238) \\
10 && 9.8906 (2.9644)  & -0.2636 (0.2124)  && 9.7976 (2.6865)  & -0.2548 (0.2050)  && 11.4789 (3.4867)  &-0.4921 (0.4433) \\
\hline
\end{tabular*}}
\end{table}

\section{Conclusions}\label{sec:conclusions}

In this work we presented several theoretical results regarding seasonal FIEGARCH (SFIEGARCH) processes.  The
similarities/differences between this model and the PLM-EGARCH model, introduced by \cite{BEA2009}, were also
discussed.

We proved here that $\{\ln(\sigma_t^2)\}_{t\in\Z}$ is a SARFIMA$(0,d,0)\times(p,0,q)_s$ process.  With this
result we provided a complete description of the process $\{\ln(\sigma_t^2)\}_{t\in\Z}$ given that necessary
and sufficient conditions for the existence, stationarity and ergodicity, as well as the autocovariance
structure and spectral representation of the SARFIMA processes are well known.  These results were used to
establish the conditions for the existence, stationarity and ergodicity of the process $\{X_t\}_{t \in
  \mathds{Z}}$ itself.  We also provided conditions for the existence of the $r$-th moment of the random
variables $\{\sigma_t^2\}_{t\in\Z}$ and $\{X_t^2\}_{t\in\Z}$ when the underlying distribution is GED.
Expressions for the asymmetry and kurtosis measures of any stationary SFIEGARCH process were also derived.

In this paper we also contributed to the theory of SARFIMA$(0,d,0)\times(p,0,q)_s$ processes by extending the
range of the parameter $d$ for the invertibility property and by providing an alternative asymptotic
expression for the autocovariance function $\gamma_{\ln(\sigma_t^2)}(\cdot)$ of the process
$\{\ln(\sigma_t^2)\}_{t\in\Z}$.  We also derived the exact and the asymptotic expressions for the
autocovariance and spectral density functions of the process $\{\ln(X_t^2)\}_{t\in\Z}$.

As an illustration, we analyzed the behavior of the intraday volatility of the S\&P$500$ US stock index
log-returns in the period from December 13, 2004 to October 10, 2009.    To account for serial correlation in
the log-return time series we considered a constrained  ARMA  model.  An SFIEGARCH model was used to
account for both the long memory and seasonal behavior for the volatility.      FIEGARCH and EGARCH models
were also considered  in order to analyze the influence of including or not the seasonal parameter in the
volatility equation.   We conclude that,  for this particular time series, not including the seasonal
parameter (FIEGARCH model) or  ignoring the  long-memory behavior (EGARCH model) lead to models which are
close enough to the non-stationary region.

\section*{Acknowledgments}

S.R.C. Lopes was partially supported by CNPq-Brazil, by CAPES-Brazil, by INCT em Matem\'atica and by Pronex
{\it Probabilidade e Processos Estoc\'asticos} - E-26/170.008/2008 -APQ1.  T.S. Prass was supported by
CNPq-Brazil. The authors are grateful to the (Brazilian) National Center of Super Computing (CESUP-UFRGS) for
the computational resources, to the Editor and two anonymous referees whose comments lead to considerable
improvement of this paper.

\appendix

\section*{Appendix A: Proofs}
\renewcommand{\theequation}{A.\arabic{equation}}

In this section we provide the proofs of all propositions, lemmas, corollaries and theorems stated in Sections
\ref{sec:sfiegarch} - \ref{sec:spectralrep}, in the same order as they appear in the text.

\subsubsection*{Proof of Proposition \ref{prop:gzt}:}
Straightforward.
\qed

\subsubsection*{Proof of Theorem \ref{thm:convergenceOrder}:}

From \eqref{eq:pik} and \eqref{eq:polilambda}, one has
\[
\sum_{k=0}^\infty\lambda_{d,k}z^k = \lambda(z)=\frac{\alpha(z)}{\beta(z)}(1-z^s)^{-d}
=\Bigg(\sum_{k=0}^{\infty}f_{k}z^k\Bigg)\Bigg(\sum_{j=0}^{\infty}\pi_{d,j}z^j\Bigg)
=\sum_{k=0}^{\infty}\Bigg(\sum_{j=0}^{k}f_{j}\pi_{d,k-j}\Bigg)z^k.
\]
\noindent It follows that, $\lambda_{d,k} = \sum_{j=0}^{k}f_{j}\pi_{d,k-j}$
or, equivalently,
\begin{equation}\label{eq:coefslambda}
  \lambda_{d,sk+r} = \sum_{j=0}^{sk+r}f_{j}\pi_{d,sk+r-j} =
  \sum_{j=0}^{k}f_{sj+r}\pi_{d,sk-sj}, \quad \mbox{for all $k \geq 0$  and $r\in
    \{0, \cdots, s-1\}$.}
\end{equation}
Let $\mathscr{K}(\cdot)$ be defined as
\[
\mathscr{K}(sk+r) = \sum_{j=0}^{\lfloor \sqrt{k}\rfloor}f_{sj+r}\frac{\pi_{d,sk-sj}}{\pi_{sk}}, \quad
\mbox{for all } k\in \N \quad \mbox{and all } r\in \{0,\cdots, s-1\}.
\]
From expression \eqref{eq:stirling}, $\pi_{d,sk-sj} \sim \pi_{d,sk}$, uniformly, for all $0 \leq j \leq
\lfloor \sqrt{k}\rfloor$.  Also, since $\beta(\cdot)$ has no roots in the closed disk $\{z: |z|\leq 1\}$,
there exist constants $a>0$ and $B>0$ such that $|f_k|< Be^{-ak}$, for all $k\in\N$
\citep[see][]{KT1994}. Hence, $\sum_{j>k_0}|f_j| \leq \frac{B\,e^{-a(k_0 + 1)}}{1-e^{-a}}$, for all
$k_0\in\N$. It follows that,
\[
\displaystyle\lim_{k\to\infty}\sum_{r=0}^{s-1}\mathscr{K}(sk+r) = \sum_{r=0}^{s-1}\bigg[\lim_{k\to\infty}
\sum_{j=0}^{\lfloor \sqrt{k}\rfloor}f_{sj+r}\frac{\pi_{d,sk-sj}}{\pi_{sk}}\bigg] = \sum_{j=0}^\infty f_j =
\frac{\alpha(1)}{\beta(1)}.
\]
Moreover,
\[
\bigg|\sum_{j= \lfloor \sqrt{k}\rfloor + 1}^{k}\hspace{-5pt}f_{sj+r}\pi_{d,sk-sj}\bigg| \leq \hspace{5pt}
\sup_{k\in\N}\big\{|\pi_{d,k}|\big\}\hspace{-5pt} \sum_{j > \lfloor \sqrt{k}\rfloor}\hspace{-5pt}|f_{sj+r}|
= o(k^{-\nu}), \quad \mbox{for any } \nu > 0.
\]
\noindent Thus, \vspace{-1\baselineskip}
\[
\lambda_{d,sk+r} = \pi_{d,sk} \sum_{j=0}^{\lfloor \sqrt{k}\rfloor}f_{sj+r}\frac{\pi_{d,sk-sj}}{\pi_{sk}}
+\hspace{-5pt} \sum_{j=\lfloor \sqrt{k}\rfloor+1}^{k}\hspace{-7pt}f_{sj+r}\pi_{d,sk-sj} =
\pi_{d,sk}\,\mathscr{K}(sk+r) + o(h ^{-\nu}),
\]
for any $\nu>0$, and
\[
\displaystyle \lim_{k\to \infty}\bigg[\bigg(\sum_{r=0}^{s-1}\lambda_{d,sk+r}\bigg)\bigg(
\pi_{sk}\frac{\alpha(1)}{\beta(1)}\bigg)^{-1}\bigg] =
\frac{\beta(1)}{\alpha(1)}\lim_{k\to\infty}\bigg[\sum_{r=0}^{s-1}\mathscr{K}(sk+r) +
\frac{o(h^{-\nu})}{\pi_{sk}}\bigg] = 1.
\]
Therefore, the result holds. \qed

\subsubsection*{Proof of Theorem \ref{thm:convergenceOrder2}:}

From expression \eqref{eq:coefslambda}, for all $r\in\{0,\cdots,s-1\}$,
\begin{align*}
  \lambda_{d,sk+r} - \frac{1}{\Gamma(d)k^{1-d}}\frac{\alpha(1)}{\beta(1)} & =
  \sum_{j=0}^{sk+r}f_{j}\left[\pi_{d,sk+r-j} -
    \frac{1}{\Gamma(d)k^{1-d}}\right] - \frac{1}{\Gamma(d)k^{1-d}}\sum_{j> sk+r} f_j\\
  & = \sum_{j=0}^{k}f_{sj+r}\left[\pi_{d,sk-sj} - \frac{1}{\Gamma(d)k^{1-d}}\right] -
  \frac{1}{\Gamma(d)k^{1-d}}\sum_{j=0}^{sk+r}f_{j}\Ind{\R\backslash
    \N}{\frac{|j-r|}{s}}\\
  & \quad - \frac{1}{\Gamma(d)k^{1-d}}\sum_{j> sk+r} f_j.
\end{align*}
Since $\sum_{j=0}^\infty |f_j| \leq \infty$, there exist $a>0$ and $B>0$ such
that $|f_j| < Be^{-aj}$, for all $j\in\N$. Consequently,
\[
\frac{1}{\Gamma(d)k^{1-d}}\hspace{-5pt}\sum_{j>sk+r}\hspace{-5pt}f_j =
o(k^{-\nu}) \quad \mbox{and} \quad
\frac{1}{\Gamma(d)k^{1-d}}\sum_{j=0}^{sk+r}f_{j}\,\Ind{\R\backslash\N}{\frac{|j-r|}{s}}
= O(k^{d-1})\Ind{\N\backslash\{0,1\}}{s},
\]
for any $\nu>0$, as $k$ goes to infinity.  From Theorem \ref{thm:convergenceOrder}, one
concludes that
\begin{align*}
  \sum_{j=0}^{k}\hspace{-1pt}f_{sj+r}\hspace{-1pt}\left[\pi_{d,sk-sj} -
    \frac{1}{\Gamma(d)k^{1-d}}\right] & =
  \sum_{j=0}^{\lfloor\sqrt{k}\rfloor}\hspace{-1pt}f_{sj+r}\hspace{-1pt}\left[\pi_{d,sk-sj} -
    \frac{1}{\Gamma(d)k^{1-d}}\right] + o(k^{-\nu}) -    \frac{1}{\Gamma(d)k^{1-d}}\hspace{-10pt}\sum_{j=\lfloor\sqrt{k}\rfloor+1}^{k}\hspace{-10pt}f_{sj+r}\\
  & = \sum_{j=0}^{\lfloor\sqrt{k}\rfloor}\hspace{-1pt}f_{sj+r}\hspace{-1pt}\left[\pi_{d,sk-sj} -
    \frac{1}{\Gamma(d)k^{1-d}}\right] + O(k^{2-d}),
\end{align*}
for any $r\in\{0,\cdots, s-1\}$.  Since
\begin{align*}
  \int_{1} ^{\sqrt{k}}e^{-a(sx+r)}\left[\Big(\frac{k-x}{k}\Big)^{d-1} -
    1\right]\mathrm{d}x &= e^{-ar}\int_{1}
  ^{\sqrt{k}}e^{-asx}\left[\Big(1-\frac{x}{k}\Big)^{d-1} -
    1\right]\!\mathrm{d}x\quad \mbox{\scriptsize (setting $x = ky$)}\\
  & = e^{-ar}k\int_{1/k}^{\sqrt{k}/k} e^{-asky} \big[(1-y)^{d-1} -
  1\big]\mathrm{d}y\\
  & \leq  (1-d)2^{1-d}e^{-ar}k\int_{0}^{\sqrt{k}/k} e^{-asky}y\, \mathrm{d}y \quad \mbox{\scriptsize (setting $u = asky$)}\\
  & = (1-d)2^{1-d}\frac{e^{-ar}}{ask}\int_{0}^{as\sqrt{k}} ue^{-u}\mathrm{d}u
  = O(k^{-1}), \quad \mbox{as } k\to \infty,
\end{align*}
by using equality \eqref{eq:pik}, one concludes that
\begin{align*}
  \sum_{j = 0} ^{ \lfloor\sqrt{k}\rfloor}\hspace{-4pt}f_{sj+r}
  \left[\pi_{d,sk-sj} - \frac{1}{\Gamma(d)k^{1-d}}\right] & =
  O(k^{d-2})\sum_{j = 0} ^{ \lfloor\sqrt{k}\rfloor}\hspace{-4pt}f_{sj+r} +
  \frac{1}{\Gamma(d)k^{1-d}}\sum_{j = 0}
  ^{\lfloor\sqrt{k}\rfloor}\hspace{-4pt}f_{sj+r}
  \left[\Big(\frac{k-j}{k}\Big)^{d-1} - 1\right]\\
  &= O(k^{d-2}), \quad \quad \quad \mbox{as } k\to \infty.
\end{align*}
Therefore, equation \eqref{eq:lambdaapprox} holds.
\qed

\subsubsection*{Proof of Proposition \ref{prop:coefs}:}

By definition,
\begin{equation}
  \frac{\alpha(z)}{\beta(z)}(1-z^s)^{-d} = \sum_{k=0}^{\infty}\lambda_{d,k}z^k \quad \Longrightarrow
  \quad \alpha(z)  = \beta(z)(1-z^s)^d\Big(\sum_{k=0}^{\infty}\lambda_{d,k}z^k\Big).\label{eq:equality}
\end{equation}
\noindent Set
\begin{equation}
  \alpha_m^* = \left\{
    \begin{array}{ccc}
      \alpha_m, & \mbox{if} &0 \leq m \leq p,\vspace{.3cm}\\
      0,       & \mbox{if}  &  m > p,
    \end{array}
  \right. \quad \mbox{and}\quad
  \beta_m^* = \left\{
    \begin{array}{ccc}
      \beta_m,  & \mbox{if} & 	 0\leq m \leq q,\vspace{.3cm}\\
      0,   & \mbox{if} &  m > q.
    \end{array}
  \right.		 \
\end{equation}
\noindent Then,
\begin{equation}
  \beta(z)(1-z^s)^d =
  \Big(\sum_{k=0}^{\infty}-\beta_k^*z^k\Big)\Big(\sum_{k=0}^{\infty}\delta_{d,k}\B^{sk}\Big) =
  \sum_{i=0}^{\infty}\sum_{j=0}^{\infty}-\beta_i^*\delta_{d,j}\B^{sj + i} = \sum_{k=0}^{\infty} \tau_k \B^{k}, \label{eq:coefstau}
\end{equation}
where $\tau_k = \sum_{i = 0}^k-\beta_i^*\delta_{d,\frac{k-i}{s}}^*$, for all
$k \in\N$, and $\delta^*_m$ is defined in \eqref{eq:newcoefs}.

Thus, \eqref{eq:equality} can be rewritten as
\begin{align}
  \beta(\B)(1-\B^s)^d\Big(\sum_{k=0}^{\infty}\lambda_{d,k}\B^k\Big)\!\!\!
  &= \Big(\sum_{k=0}^{\infty}\tau_k\B^k\Big)\Big(\sum_{k=0}^{\infty}\lambda_{d,k}\B^k\Big)
  = \sum_{k=0}^{\infty}\Bigg(\sum_{i=0}^{k}\lambda_{d,i}\tau_{k-i}\Bigg)\B^k\nonumber\\
  \!\!\!&=\sum_{k=0}^{\infty}\Bigg[\sum_{i=0}^{k}\lambda_{d,i}\Big(-\sum_{j=0}^{k-i}\beta_j^*\delta_{d,\frac{k-i-j}{s}}^*\Big)\Bigg]\B^k\nonumber\\
  \!\!\!&=\lambda_{d,0} + \sum_{k=1}^{\infty}\Bigg[\lambda_{d,k}-\sum_{i=0}^{k-1}\lambda_{d,i}\Big(\sum_{j=0}^{k-i}\beta_j^*\delta_{d,\frac{k-i-j}{s}}^*\Big)\Bigg]\B^k.\label{eq:part1}
\end{align}

Therefore, from \eqref{eq:newcoefs} and \eqref{eq:part1}, expression
\eqref{eq:equality} holds if and only if
\begin{equation*}
  -\alpha_0^* = \lambda_{d,0} \quad \mbox{and  } \quad -\alpha_k^* = \lambda_{d,k}-\sum_{i=0}^{k-1}\lambda_{d,i}\Big(\sum_{j=0}^{k-i}\beta_j^*\delta_{d,\frac{k-i-j}{s}}^*\Big),\quad \mbox{for all} \quad k\geq 1,
\end{equation*}
\noindent and the result holds.
\qed

\subsubsection*{Proof of Lemma \ref{lemma:welldefined}:}
From Theorem \ref{thm:convergenceOrder}, $\sum_{k=0}^\infty \lambda_{d,k}^2 <\infty$ if and only if $d <
0.5$. Therefore, if \eqref{eq:sigma_sum} is a.s. convergent, by applying the three series theorem
\citep[see][]{BI1995}, one concludes that, necessarily, $d < 0.5$.  On the other hand, if $d < 0.5$, by
applying Kolmogorov's convergence criteria \citep[theorem 22.6]{BI1995}, one concludes that
\eqref{eq:sigma_sum} is a.s. convergent.  Finally, from Theorem \ref{thm:convergenceOrder}, if $d<0$,
$\sum_{k=0}^\infty|\lambda_{d,k}| <\infty$ and, from proposition 3.1.1 in Brocwel and Davis (1991), the series
\eqref{eq:sigma_sum} is absolutely a.s. convergent.  \qed

\subsubsection*{Proof of Corollary \ref{lemma:SARFIMA}:}

The result follows immediately from Proposition \ref{prop:gzt}, Lemma \ref{lemma:welldefined} and the
definition of SARFIMA processes in \cite{BL2009}. In particular, if $s = 1$, it is an ARFIMA$(q,d,p)$
process \citep[see][]{BD1991}.  \qed

\subsubsection*{Proof of Corollary \ref{cor:xfinite}:}
From Lemma \ref{lemma:welldefined}, if $d<0.5$, the random variable $\sigma_t^2$ is finite with probability
one, for all $t\in \Z$. Since, by assumption $\E(Z_t^2) = \E(Z_0^2) = 1$, the random variable $Z_t$ is finite
with probability one, for all $t\in\Z$.  Therefore, from expression \eqref{eq:xt}, the result follows.  \qed

\subsubsection*{Proof of Theorem \ref{thm:invertibility}:}

By H\"older's inequality, it suffices to prove the result for $\p=2$.  From \cite{BL2009}, the spectral
density function of $\{Y_t\}_{t\in\Z}$ is given by
\begin{equation}
  f_{\ln(\sigma_t^2)}(\lambda)=\frac{\sigma_g^2}{2\pi}
  \frac{|\alpha(e^{-i\lambda})|^2}{|\beta(e^{-i\lambda})|^2}|1-e^{-is\lambda}|^{-2d}
  = \frac{\sigma_g^2}{2\pi}   \frac{|\alpha(e^{-i\lambda})|^2}{|\beta(e^{-i\lambda})|^2}\left[2\sin\left(\dfrac{s\lambda}{2}\right)\right]^{-2d}, \label{eq:fln}
\end{equation}
\noindent for all $\lambda \in [0,\pi]$. The authors also show that, for all
$k = 0, \cdots, \left\lfloor \frac{s}{2} \right\rfloor$,
\begin{equation}
  f_Y(\lambda) \sim  \frac{\sigma_\varepsilon^2}{2\pi}
  \left|\frac{\alpha(e^{-i\lambda_k})}{\beta(e^{-i\lambda_k})}\right|^2
  s^{-2d}\left|\lambda - \lambda_k\right|^{-2d}, \quad \mbox{as} \quad
  \lambda \rightarrow \lambda_k := \frac{2\pi k }{s}.\label{eq:fln2}
\end{equation}

Suppose first that the $L^2$ convergence holds. Notice that, there exists a real constant $c>0$ such that
$f_{Y}(\lambda)\sim c\lambda^{-2d}$, as $\lambda \rightarrow 0^+$.  Consequently, from proposition 1.5.8 in
\cite{BEA1987}, if $d \geq 0.5$, the function $f_{Y}(\cdot)\notin L^1$ and hence cannot be a spectral
density function.  From Theorem \ref{thm:convergenceOrder}, if $d \leq -1$, then $\tilde\lambda_{d,k}
\nrightarrow 0$, when $k\rightarrow \infty$.  Consequently, $\tilde\lambda_{d,k}Y_{t-k} \nrightarrow 0$ in
$L^\p$-norm, when $k\rightarrow \infty$, and the series representation $\sum_{k =
  0}^{\infty}\tilde\lambda_{d,k}Y_{t-k}$ cannot converge in $L^\p$-norm, for any $0<\p\leq 2$. Therefore,
necessarily, $d\in (-1,0.5)$.

Suppose now that $d\in(-1, 0.5)$ and assume $s>1$ \citep[the case $s=1$ can be found in][]{BP2007}. Notice
that $d = d_1 + d_2$, where $d_1\in(-0.5,0)$ and $d_2 \in (-0.5,0.5)$.  Define the functions $f_1(\cdot)$
and $f_2(\cdot)$ as follows,
\[ f_{1}(\lambda) := \frac{f_{Y}(\lambda)}{f_2(\lambda)} \quad \mbox{and}
\quad f_{2}(\lambda):= |1-e^{\im s\lambda}|^{-2d_2}, \quad \mbox{for all }
\lambda \in [0,\pi].
\]
\noindent Since $d_1\in(-0.5,0)$, from expressions \eqref{eq:fln} and \eqref{eq:fln2}, it is obvious that,
$f_1(\cdot) \in L^{\infty}$ and
\[
f_1(\lambda) \sim l\Big(\frac{1}{\lambda - 2\pi k/s}\Big)\Big|\lambda -
\frac{2\pi k}{s}\Big|^{-2d_1}, \quad \mbox{as } \lambda \rightarrow
\frac{2\pi k}{s}, \quad \mbox{for all } 0 \leq k \leq \left\lfloor
  \frac{s}{2}\right\rfloor,
\]
\noindent where $l(\cdot)$ is a slowly varying function at infinity \citep[for details, see][]{BEA1987}.
Since $d_2 \in (-0.5, 0.5)$, from expression \eqref{eq:fln}, one easily concludes that, for any $\epsilon
\in (0, 2\pi/s)$, the function $f_{1}^{-1}(\cdot)$ is bounded in the interval $[2\pi k/s + \epsilon, 2\pi
(k+1)/s - \epsilon]$, for all $0 \leq k \leq \lfloor s/2 \rfloor$. Moreover, from expression
\eqref{eq:fln2}, there exists an $\epsilon \in (0, 2\pi/s)$, such that
\[
f_{1}^{-1}(\lambda) \leq 2l_1\Big(\frac{1}{\lambda - 2\pi k/s}\Big)\Big|\lambda
- \frac{2\pi k}{s}\Big|^{2d_1}, \quad \mbox{ for all} \quad \lambda \in
\Big(\frac{2\pi k}{s} - \epsilon, \frac{2\pi k}{s}\Big)\cup
\Big(\frac{2\pi k}{s}, \frac{2\pi k}{s} + \epsilon\Big)
\]
\noindent and $k\in\{0, \cdots ,\lfloor s/2 \rfloor\},$ where $l_1(\cdot) = l^{-1}(\cdot)$, is also a slowly
varying function.  It follows that $\int_0^\pi f_1^{-1}(\lambda)d\lambda$ can be written as
\[
\int_0^\pi f_1^{-1}(\lambda)d\lambda   = \left\{
  \begin{array}{cc}
    S_1, & \mbox{if $s$ is even;}\vspace{0.2cm}\\
    S_1 + S_2, & \mbox{if $s$ is odd,}
  \end{array}
\right.
\]
\noindent where
\begin{align*}
  S_1 & := \sum_{k = 0}^{\lfloor s/2\rfloor -1}\left[
    \int_{2\pi k/s}^{2\pi k/s+\epsilon}f_1^{-1}(\lambda)d\lambda +
    \int_{2\pi k/s+\epsilon}^{2\pi(k+1)/s-\epsilon}f_1^{-1}(\lambda)d\lambda
    + \int_{2\pi(k+1)/s - \epsilon}^{2\pi(k+1)/s}f_1^{-1}(\lambda)d\lambda
  \right], \\
  S_2 &:= \int_{ 2\pi\lfloor s/2\rfloor/s}^{2\pi\lfloor s/2 \rfloor /s +
    \epsilon}f_1^{-1}(\lambda)d\lambda + \int_{2\pi\lfloor s/2\rfloor/s +
    \epsilon}^{\pi}f_1^{-1}(\lambda)d\lambda,\\
\end{align*}

\vspace{-1\baselineskip}
\noindent and it satisfies
\[
\quad \int_0^\pi f_1^{-1}(\lambda)d\lambda   \leq \left\{
  \begin{array}{cc}
    S_1^*, & \mbox{if $s$ is even;}\vspace{0.2cm}\\
    S_1^* + S_2^*, & \mbox{if $s$ is odd,}
  \end{array}
\right.
\]
\noindent where
\begin{align*}
  S_1^* & := \sum_{k = 0}^{\lfloor s/2 \rfloor -1}\left[
    2\int_{2\pi k/s}^{2\pi k/s+\epsilon}l_1\Big(\frac{1}{\lambda -
      2\pi k/s}\Big)\Big|\lambda - \frac{2\pi k}{s}\Big|^{2d_1}d\lambda +
    \int_{2\pi k/s+\epsilon}^{2\pi(k+1)/s-\epsilon}f_1^{-1}(\lambda)d\lambda
  \right.\\
  & \left. \hspace{2cm} + \ 2\int_{2\pi(k+1)/s -
      \epsilon}^{2\pi(k+1)/s}l_1\Big(\frac{1}{\lambda -
      2\pi k/s}\Big)\Big|\lambda - \frac{2\pi k}{s}\Big|^{2d_1}d\lambda
  \right],\\
  S_2^* &:= 2\int_{2\pi\lfloor s/2\rfloor/s}^{\pi\lfloor s/2\rfloor2/s +
    \epsilon}l_1\Big(\frac{1}{\lambda - 2\pi k/s}\Big)\Big|\lambda -
  \frac{2\pi k}{s}\Big|^{2d_1}d\lambda + \int_{ 2\pi\lfloor s/2\rfloor/s +
    \epsilon}^{\pi}f_1^{-1}(\lambda)d\lambda.
\end{align*}
\noindent Thus, since $d_1 > -0.5$, by proposition 1.5.10 in \cite{BEA1987} we conclude that
$f_1^{-1}(\cdot) \in L^1$.  By comparing the function $f_2(\cdot)$ with the corresponding one in
\cite{BP2007}, we conclude that $f_2(\lambda) = |1-e^{\im s\lambda}|$ satisfies condition $A_{p}$, with
$p=2$, from theorem 3 in \cite{B1985}.  Therefore, taking $p =1$ in theorem 4 from \cite{B1985} we conclude
that $\tilde\lambda(\cdot)$ has a Fourier series that converges in $L^2(f)$, where $F(\cdot)$ is the
spectral distribution function of $Y_t$, for all $t\in\Z$, and the result follows.  \qed

\subsubsection*{Proof of Lemma \ref{lemma:lsstationary}:}
See \cite{BL2009}.
\qed

\subsubsection*{Proof of Corollary \ref{cor:sigmastationary}:}
It follows immediately from Lemma \ref{lemma:lsstationary}.
\qed

\subsubsection*{Proof of Theorem \ref{thm:est_ergod}}
Suppose $d<0.5$. (i) By hypothesis, $Z_t$ is finite with probability one, for all $t\in\Z$. From Corollary
\ref{lemma:welldefined}, the random variables $\ln(\sigma_t^2)-\omega$ is finite with probability one, for all
$t\in\Z$, so it is $\sigma_t$.  Thus, from theorem 3.5.8 in \cite{S1974}, $\{X_t\}_{t \in \mathds{Z}}$ is a
strictly stationary and ergodic process.\vspace{0.1cm}

\noindent (ii) Assume that $\E(|\ln(Z_0^2)|) < \infty$. It follows that the random variable $|\ln(Z_t^2)\big|$
is finite with probability one, for all $t\in\Z$. From Corollary \ref{lemma:welldefined},
$\E(|\ln(\sigma_t^2)|) < \infty$, for all $t\in\Z$.  By expression \eqref{eq:sigmat}, $\ln(X_t^2) =
\ln(\sigma_t^2) + \ln(Z_t^2)$, for all $t\in\Z$. It follows that the random variable $\ln(X_t^2)$ is finite
with probability one, for all $t\in\Z$, and hence the stochastic process $\{\ln(X_t^2)\}_{t\in\Z}$ is well
defined.  From theorem 3.5.8 in \cite{S1974}, the stochastic process $\{\ln(X_t^2)\}_{t\in\Z}$ is strictly
stationary and ergodic.  Moreover, if $\E([\ln(Z_0^2)]^2) < \infty$ then $\var\big(\ln(X_t^2)\big) =
\var\big(\ln(X_0^2)\big) = \var(\ln(\sigma_0^2)) + \var(\ln(Z_0 ^2)) <\infty$, for all $t\in\Z$. Therefore,
$\{\ln(X_t^2)\}_{t\in\Z}$ is weakly stationary.  \qed

\subsubsection*{Proof of Theorem \ref{thm:estacionar1}:}

From Corollary \ref{cor:sigmastationary} and Theorem \ref{thm:est_ergod} both processes $\{X_t\}_{t \in
  \mathds{Z}}$ and $\{\sigma_t\}_{t\in\Z}$ are strictly stationary and hence, any existing moments are time
invariant. Let $r>2$ be any real number such that $\E(|Z_t|^{\r})<\infty$. By the independence hypothesis
$\E(|X_t|^{\r}) = \E(|\sigma_t|^{\r})\E(|Z_t|^{\r})$, for all $t\in\Z$.  Since $\E(|Z_t|^{\r})<\infty$, one
has to show that $\E(|\sigma_t|^{\r})<\infty$.  From expression \eqref{eq:sigmat}, and from the i.i.d.\!
hypothesis on the random variables $Z_t$, for all $t\in\Z$, one has
\begin{equation}\label{eq:moment2}
  \E(|\sigma_t|^{\r}) =  \E\Big(\Big|\exp\Big\{\frac{\omega}{2} +
  \frac{1}{2}\sum_{k=0}^\infty\lambda_{d,k}g(Z_{t-1-k})\Big\}\Big|^{\r}\Big)  =     e^{\frac{\r\omega}{2}}\prod_{k=0}^\infty\E\Big(\exp\Big\{\frac{\r}{2}\lambda_{d,k}g(Z_{0})\Big\}\Big).
\end{equation}
\noindent From \eqref{eq:convergencecriteria}, expression \eqref{eq:moment2} converges to a non-zero
constant \citep[see, section 0.25 in][]{GR2000}. By H\"older's inequality the result follows for all
$0<\m<\r$.  \qed

\subsubsection*{Proof of Corollary \ref{cor:estacionary}:}
By hypothesis, $d<0.5$ and $\beta(z)$ has no roots in the closed disk $\{z:|z|\leq 1\}$.  From Theorem
\ref{thm:convergenceOrder}, it follows that $\sum_{k=0}^{\infty}\lambda_{d,k}^2<\infty$. Moreover, from
expression A2.4 in \cite{N1991}, for all $r>0$, $ \E(\exp\{\r\lambda_{d,k}g(Z_{0})\}) = 1 +
O(\lambda_{d,k}^2)$, as $k\to \infty$.  Thus, the result follows immediately from Theorem
\ref{thm:estacionar1}.  \qed

\subsubsection*{Proof of Proposition \ref{prop:skewnessKurtosys}:}

Let $\{X_t\}_{t \in \mathds{Z}}$ be any stationary SFIEGARCH process.  Let $\lambda(\cdot)$ be the polynomial
defined by \eqref{eq:polilambda}.  Since $\E(X_t^\r)=\E(\sigma_t^\r Z_t^\r)$, for all $\r>0$, $\E(X_t) = 0$
and $\E(Z_t^2) =1$, for all $t\in\Z$, the asymmetry and kurtosis measures of $\{X_t\}_{t \in \mathds{Z}}$ are
given, respectively, by
\begin{equation}
  A_X =       \frac{\E(X_t^3)}{\big(\E(X_t^2)\big)^{3/2}}=\frac{\E(\sigma_t^3)\E(Z_t^3)}{\big(\E(\sigma_t^2)\big)^{3/2}}
  \quad \mbox{and} \quad K_X =       \frac{\E(X_t^4)}{\big(\E(X_t^2)\big)^2}=\frac{\E(\sigma_t^4)\E(Z_t^4)}{\big(\E(\sigma_t^2)\big)^2}. \label{eq:ak}
\end{equation}
\noindent Upon replacing expression \eqref{eq:moment2} in \eqref{eq:ak} the result follows.  \qed

\subsubsection*{Proof of Lemma \ref{lemma:expressioncovln}:}
See \cite{BL2009}.
\qed

\subsubsection*{Proof of Corollary \ref{cor:rhoconv}:}

Notice that, for any $\eta \in \N$ one can write
\[
\gamma_{\ln(\sigma_t^2)}(sh +r) = \sum_{|k|\leq
  \eta} \gamma_{\mbox{\tiny$A$}}(sk+r)\gamma_{\mbox{\tiny$V$}}(sh - sk) \ + \sum_{|k|>
  \eta } \gamma_{\mbox{\tiny$A$}}(sk+r)\gamma_{\mbox{\tiny$V$}}(sh - sk).
\]
Since $\sum_{h\in\Z}|\gamma_{\mbox{\tiny$A$}}(h)| < \infty$, there exist $a>0$ and $B>0$ such that
$|\gamma_{\mbox{\tiny$A$}}(h)|< Be^{-ah}$, for all $h\in\N$ \citep[see][]{KT1994}.  Also,
$\gamma_{\mbox{\tiny$A$}}(sk+r) = \gamma_{\mbox{\tiny$A$}}(-s|k|+r) = \gamma_{\mbox{\tiny$A$}}(s|k|-r)$, for
all $k< -\eta$. Thus,
\begin{align*}
  \Big|\sum_{|k|> \eta }
  \gamma_{\mbox{\tiny$A$}}(sk+r)\gamma_{\mbox{\tiny$V$}}(sh - sk) \Big|& \leq
  \gamma_{\mbox{\tiny$V$}}(0)\sum_{k > \eta
  }\Big(|\gamma_{\mbox{\tiny$A$}}(sk+r)| + |\gamma_{\mbox{\tiny$A$}}(sk-r) | \Big)\\
  &\leq
  \frac{B\gamma_{\mbox{\tiny$V$}}(0)
    (e^{-ar}+e^{ar})}{1-e^{-a}}e^{-as(\sqrt{h}+1)} = o(h^{-\nu}), \quad \mbox{for
    all } \nu >0.
\end{align*}
\noindent On the other hand,
\begin{equation*}
  \sum_{|k|\leq \eta }
  \gamma_{\mbox{\tiny$A$}}(sk+r)\gamma_{\mbox{\tiny$V$}}(sh - sk) =
  \gamma_{\mbox{\tiny$V$}}(sh)\sum_{|k|\leq \eta }     \gamma_{\mbox{\tiny$A$}}(sk+r)\frac{\gamma_{\mbox{\tiny$V$}}(sh-sk)}{\gamma_{\mbox{\tiny$V$}}(sh)}
  := \gamma_{\mbox{\tiny$V$}}(sh)\mathscr{G}(sh + r).
\end{equation*}
\noindent Let $\eta =\lfloor\sqrt{h}\rfloor$. Then, for any $|k|\leq \eta = \lfloor\sqrt{h}\rfloor$, one has
$(h-k)^{2d-1} \sim h^{2d-1}$, uniformly, as $h \to \infty$.  From expression \eqref{eq:covln},
\begin{align}
  \gamma_{\mbox{\tiny$V$}}(sh) &= \sigma_g^2
  \frac{(-1)^{h}\Gamma(1-2d)}{\Gamma(1-d+h)\Gamma(1-d-h)} = \sigma^2_g
  \frac{\Gamma(1-2d)}{\Gamma(1-d)\Gamma(d)} \frac{\Gamma(h+d)}{\Gamma(h+1-d)} \sim \sigma^2_g \frac{\Gamma(1-2d)}{\Gamma(1-d)\Gamma(d)}h^{2d-1}, \quad
  \mbox{as} \quad h\to \infty, \label{eq:covlninf}
\end{align}
where $\sigma_g^2$ is given by \eqref{eq:sigma2gzt}. Hence, $\gamma_{\mbox{\tiny$V$}}(sh- sk) =
\gamma_{\mbox{\tiny$V$}}(s(h-k)) \sim \gamma_{\mbox{\tiny$V$}}(sh)$, uniformly, for all $|k|\leq \eta =
\lfloor\sqrt{h}\rfloor$. Since $\sum_{h\in\Z}|\gamma_{\mbox{\tiny$A$}}(h)| < \infty$, it follows that
\[
\lim_{h\to \infty} \sum_{r=0}^{s-1}\mathscr{G}(sh + r) = \lim_{h\to \infty}
\sum_{r=0}^{s-1}\sum_{|k|\leq \eta }
\gamma_{\mbox{\tiny$A$}}(sk+r)\frac{\gamma_{\mbox{\tiny$V$}}(sh-sk)}{\gamma_{\mbox{\tiny$V$}}(sh)}
= \sum_{k\in\Z }\gamma_{\mbox{\tiny$A$}}(k),
\]
\noindent and expression \eqref{eq:covlsinfty} holds.
\qed

\subsubsection*{Proof of Theorem \ref{thm:resultsXt}:}

From Theorem \ref{thm:est_ergod}, $\{\ln(X_t^2)\}_{t\in\Z}$ is a stationary (weakly and strictly)
process. Thus,
\[
\gamma_{\ln(\X_t^2)}(h) =
\cov(\ln(X_t^2),\ln(X_{t+h}^2)) = \cov(\ln(X_0^2),\ln(X_{h}^2)), \quad \mbox{for all } t,h\in\Z
\]
From expression \eqref{eq:xt}, one concludes that
\begin{align*}
  \gamma_{\ln(\X_t^2)}(h)  & =  \cov(\ln(\sigma_0^2),\ln(\sigma_{h}^2))  +
  \cov(\ln(\sigma_0^2),\ln(Z_{h}^2)) + \cov(\ln(Z_0^2),\ln(\sigma_{h}^2)) +\cov(\ln(Z_0^2),
  \ln(Z_{h}^2)),
\end{align*}
for all $h\in\Z.$ Notice that $\cov(\ln(\sigma_0^2),\ln(\sigma_{h}^2)) = \gamma_{\ln(\sigma_t^2)}(h)$ is the
autocovariance function of $\{\ln(\sigma_t^2)\}_{t\in\Z}$, given in Lemma \ref{lemma:expressioncovln}, and
$\cov(\ln(Z_0^2), \ln(Z_{h}^2)) = \sigma^2_{\mbox{\tiny$\ln(Z_t^2)$}}\Ind{\{0\}}{h}$. Moreover, from
expression \eqref{eq:sigmat},
\begin{equation}
  \cov(\ln(\sigma_{h}^2), \ln(Z_{0}^2))  =  \cov\Big(\sum_{k=0}^\infty
  \lambda_{d,k}g(Z_{h-1-k}),\ln(Z_{0}^2)\Big)
  = \left\{
    \begin{array}{ccc}
      0, & \mbox{if}& h\leq 0;\\
      C_1\lambda_{d,h-1},&\mbox{if}& h > 0,
    \end{array}
  \right.\label{eq:covln02}
\end{equation}
\noindent where $ C_1 = \cov(\ln(Z_0^2),g(Z_0)) = \theta\,\E( Z_0\ln(Z_0^2)) +\gamma \E( |Z_0|\ln(Z_0^2)) -
\gamma \E( |Z_0|)\E(\ln(Z_0^2))$. Thus, from expression \eqref{eq:covln02},
$\cov(\ln(\sigma_0^2),\ln(Z_{h}^2)) + \cov(\ln(Z_0^2),\ln(\sigma_{h}^2)) =
C_1\lambda_{d,|h|-1}\Ind{\Z^*}{h}$, for all $h\in\Z$, and expression \eqref{eq:covlx} holds.  Expression
\eqref{eq:gammainfty} follows directly from Corollary \ref{cor:rhoconv} and Theorem
\ref{thm:convergenceOrder}.  \qed

\subsubsection*{Proof of Corollary \ref{cor:asymptcov}:}

From Corollary \ref{cor:rhoconv} and equation
\eqref{eq:covlninf},
\[
\sum_{r=0}^{s-1}\gamma_{\ln(\sigma_t^2)}(sh + r) =  \gamma_{\mbox{\tiny$V$}}(sh)\sum_{r=0}^{s-1}\mathscr{G}(sh+r) + o(h^{-\nu})  \sim \mathscr{C}_1 h^{2d-1},\quad \mbox{as } h\to \infty.
\]
From Theorem \ref{thm:convergenceOrder} and equation \eqref{eq:stirling} one concludes that $\mathscr{K}(sh
- 1) = \mathscr{K}(s(h-1) + (s-1))$, for all $h\geq1$, and
\[
\sum_{r=0}^{s-1}\lambda_{d,sh+r-1}   = \pi_{d,s(h-1)}  \mathscr{K}(sh-1) +
\pi_{d,sh}\sum_{r=1}^{s-1}\mathscr{K}(sk+r-1)  + o(h^{-\nu}) \sim
\mathscr{C}_2k^{d-1},  \mbox{ as } h\to\infty,
\]
for any $\nu>0$, where $\mathscr{C}_1(\cdot)$ and $\mathscr{C}_2(\cdot)$ are given in equation
\eqref{eq:c1c2}. Since
\begin{align*}
  \lim_{h\rightarrow  \infty} \dfrac{h^{2d-1}}{h^{d-1}} =  \lim_{h\rightarrow
    \infty}h^d = 0, \quad \mbox{if } d<0, & \qquad
  \lim_{h\rightarrow  \infty} \frac{h^{d-1}}{h^{2d-1}} =  \lim_{h\rightarrow
    \infty}h^{-d} = 0, \quad \mbox{if } d>0,\\
  \lim_{h\to \infty} \mathscr{C}_1(h) = \sigma^2_g
  \frac{\Gamma(1-2d)}{\Gamma(1-d)\Gamma(d)}\sum_{k\in\Z}\gamma_{\mbox{\tiny$A$}}(k)
  & \qquad \mbox{and} \qquad
  \lim_{h\to \infty} \mathscr{C}_2(h) = \frac{1}{\Gamma(d)}\frac{\alpha(1)}{\beta(1)},
\end{align*}
expression \eqref{eq:asymptcov}  holds.
\qed

\subsubsection*{Proof of Theorem \ref{thm:resultsXt2}:}

From Theorem \ref{thm:est_ergod}, if $d<0.5$, the stochastic process $\{\ln(X_t^2)\}_{t\in\Z}$ is strictly
stationary and ergodic. Moreover, if $\var(\ln(Z_0^2)) < \infty$, then it is weakly stationary and hence, it
has a spectral distribution function.  Thus, from Herglotz's theorem \citep[see][]{BD1991}, it suffices to
show that $f_{\ln(\X_t^2)}(\cdot)$, given by \eqref{eq:spectraldensity}, is a continuous, non-negative
function and it satisfies
\begin{equation*}
  \gamma_{\ln(\X_t^2)}(h)=\int_{(-\pi,\pi\mbox{]}}e^{\im h\lambda}f_{\ln(\X_t^2)}(\lambda)d\lambda, \quad \mbox{for all } h\in \Z,
\end{equation*}
\noindent with $\gamma_{\ln(\X_t^2)}(\cdot)$ given in Theorem \ref{thm:resultsXt}.

The continuity of $f_{\ln(\X_t^2)}(\cdot)$ follows immediately from its definition.  To prove non-negativity
notice that, from expression \eqref{eq:fln} and from the i.i.d.\! hypothesis on the $\{Z_t\}_{t \in
  \mathds{Z}}$ process, one has
\[
f_{\ln(\sigma_t^2)}(\lambda)  = \frac{\sigma^2_g}{2\pi}|\Lambda(\lambda)|^2
\quad \mbox{and} \quad  f_{\ln(Z_t^2)}(\lambda)  =
\frac{\sigma^2_{\ln(Z_t^2)}}{2\pi}, \quad \mbox{for all } \lambda \in [0,\pi],
\]
\noindent where $\Lambda(z) := \lambda(e^{-\im z})$ and $\lambda(\cdot)$ is
defined in \eqref{eq:polilambda}. Moreover,  $|C_1| \leq \sigma_g
\sigma_{\ln(Z_t^2)}$ and $|\Re(z)| \leq |z|$,  for all   $z\in\mathds{C}$. Thus,
\begin{align*}
  f_{\ln(\sigma_t^2)}(\lambda) +  \frac{C_1}{\pi}\Re\big(e^{-\im \lambda}\Lambda(\lambda)\big) +
  f_{\ln(Z_t^2)}(\lambda)  &\geq  f_{\ln(\sigma_t^2)}(\lambda) -
  \frac{\sigma_g \sigma_{\ln(Z_t^2)}}{\pi}\big|e^{-\im \lambda}\Lambda(\lambda)\big| +
  f_{\ln(Z_t^2)}(\lambda) \phantom{\Bigg|}\\
  & =  \left(\frac{\sigma_g}{\sqrt{2\pi}}|\Lambda(\lambda)| -
    \frac{\sigma_{\ln(Z_t^2)}}{\sqrt{2\pi}}\right)^2 \phantom{\Bigg|} \geq 0,
  \quad \mbox{for all $\lambda \in [0,\pi]$.}
\end{align*}
\noindent To complete the proof, observe that
\begin{equation*}
  \gamma_{\ln(\sigma_t^2)}(h)=\int_{(-\pi,\pi\mbox{]}}e^{\im
    h\lambda}f_{\ln(\sigma_t^2)}(\lambda)d\lambda,  \quad \quad \gamma_{\ln(Z_t^2)}(h)=\int_{(-\pi,\pi\mbox{]}}e^{\im h\lambda}f_{\ln(Z_t^2)}(\lambda)d\lambda\nonumber
\end{equation*}
and
\begin{align*}
  \int_{(-\pi,\pi\mbox{]}} e^{\im h \lambda} \frac{C_1}{\pi}\Re\big(e^{-\im
    \lambda}\Lambda(\lambda)\big)  d\lambda &=
  \frac{C_1}{\pi}\sum_{k=0}^\infty\lambda_{d,k}\int_{(-\pi,\pi\mbox{]}} e^{\im h
    \lambda}  \cos((k+1)\lambda) d\lambda\\
  &= \frac{C_1}{\pi}\sum_{k=0}^\infty\lambda_{d,k}\pi \Ind{\{0\}}{k+1-|h|}
  = C_1\lambda_{d,|h|-1}\Ind{\N^*}{|h|},
\end{align*}
for all  $h \in \Z$. Therefore, the result holds.
\qed

\subsubsection*{Proof of Theorem \ref{thm:predictor}}
Let $\mathscr{S}_1 := \prod_{k = 0}^{h-2}\exp\big\{\lambda_{d,k}g(Z_{n+h-1-k})\big\}$ and $\mathscr{S}_2 :=
\prod_{k = h-1}^\infty\exp\big\{\lambda_{d,k}g(Z_{n+h-1-k})\big\}$, for any $n\in \N$ and $h>1$ fixed.  Notice
that, from expression \eqref{eq:sigmat}, one can write
\begin{equation}\label{eq:cond_moment}
  \sigma_{n+h}^2 =  e^\omega \prod_{k = 0}^\infty\exp\big\{\lambda_{d,k}g(Z_{n+h-1-k})\big\} := e^\omega\mathscr{S}_1\mathscr{S}_2, \quad \mbox{for all } n\in \N  \mbox{ and }  h>1.
\end{equation}
Also, observe that  the hypothesis $\E(\sigma_t^2) < \infty$ implies  $0 < \mathscr{S}_1, \mathscr{S}_2 < \infty$  with probability one.

The  $\mathcal{F}_{n}$-measurability of $Z_{n+h-1-k}$,  when $k \geq h-1$, and the the i.i.d. property of  $\{Z_t\}_{t\in\mathds{Z}}$ imply that
\[
\E( \sigma_{n+h}^2|\mathcal{F}_n) = e^\omega\mathscr{S}_2\E(\mathscr{S}_1) =
e^\omega\prod_{k=h-1}^{\infty}\!\!\exp\big\{\lambda_{d,k}g(Z_{n+h-1-k})\big\}
\prod_{k=0}^{h-2}\E\big(\exp\big\{\lambda_{d,k}g(Z_{0})\big\}\big),
\]
and  expression \eqref{eq:predictor} holds.

Now,  the independence of $\mathscr{S}_1$ and $\mathscr{S}_2$ implies that
\[
\E\big(\big[\sigma_{n+h}^2 - \hat \sigma_{n+h}^2 \big]^2\big) = \E\big(\big[e^\omega\mathscr{S}_1\mathscr{S}_2 - e^\omega\mathscr{S}_2\E(\mathscr{S}_1)\big]^2\big) =  e^{2\omega}\E(\mathscr{S}_2^2)\E\big(\big[\mathscr{S}_1 - \E(\mathscr{S}_1)\big]^2\big), \quad \mbox{for any } h>1.
\]
Since $\E\big(\big[\mathscr{S}_1 - \E(\mathscr{S}_1)\big]^2\big) = \E(\mathscr{S}_1^2) -
\big[\E(\mathscr{S}_1)\big]^2$, expression \eqref{eq:mse_sigma} holds.

To conclude the proof observe that
\[
\E\big(\big[X_{n+h}^2 - \hat X_{n+h}^2 \big]^2\big) = \E\big(\big[X_{n+h}^2 - \hat \sigma_{n+h}^2 \big]^2\big)
= \E\big(\big[e^\omega\mathscr{S}_1\mathscr{S}_2Z_{n+h}^2 -
e^\omega\mathscr{S}_2\E(\mathscr{S}_1)\big]^2\big), \quad \mbox{for any } h>1.
\]
Then, by adding and subtracting $\E\big(\big[e^{\omega}\mathscr{S}_1\mathscr{S}_2\big]^2\big)$ to the right
hand side of the above equation, one can rewrite
\begin{align*}
  \E\big(\big[X_{n+h}^2 - \hat X_{n+h}^2 \big]^2\big) &= \E\big(\big[e^{\omega}\mathscr{S}_1\mathscr{S}_2\big]^2\big)\big[\E(Z_{n+h}^4) - 1\big] +  e^{2\omega}\E(\mathscr{S}_2^2)\big( \E(\mathscr{S}_1^2) - \big[\E(\mathscr{S}_1)\big]^2\big)\\
  & = \E(\sigma_{0}^4)\big[\E(Z_0^4) - 1\big] + mse(\sigma_{n+h}^2), \quad \mbox{for any } h>1,
\end{align*}
and the result follows.  \qed

\subsubsection*{Proof of Theorem \ref{thm:forecast_log}}
Expression \eqref{eq:forecast_log} and the first equation in \eqref{eq:mse_log} follow immediately by
mimicking the proof of proposition 4 in \cite{LP2013}.  The second equation in \eqref{eq:mse_log} is obtained
by replacing $\check \ln(X_{n+h}^2)$ by $\hat \ln(X_{n+h}^2)$ and mimicking the proof of proposition 5 in
\cite{LP2013}.  \qed

\subsubsection*{Proof of Theorem \ref{thm:arma_sfiegarch}}
From   \eqref{eq:rt}, for any fixed $n\in \Z$,
\[
r_{n+h}^2 = \mu^2 + 2\mu\sum_{i = 0}^\infty\psi_iX_{n+h-i} + \sum_{k = 0}^\infty\sum_{j =
  0}^\infty\psi_k\psi_jX_{n+h-k}X_{n+h-j}, \quad \mbox{for all } h>0.
\]
Thus the result follows by observing that $\E(X_{n+h-i}|\mathcal{F}_n) = X_{n+h-i}$, if $i\geq h$, and 0
otherwise, and that
\[
\E(X_{n+h-k}X_{n+h-j} | \mathcal{F}_n) =
\left\{
  \begin{array}{cc}
    \E(X_{n+h-k}^2| \mathcal{F}_n) = \hat \sigma_{n+h-k}^2 ,  & \mbox{for any } k,j < h, \mbox{ with } k = j;\\
    0, & \mbox{for any } k,j < h, \mbox{ with } k\neq j;\\
    X_{n+h-k}X_{n+h-j}, & \mbox{for any } k,j \geq h, \mbox{ with } k\neq j.
  \end{array}
\right.
\]
for any $k,j \in \N$ and $h>0$. \qed

\end{document}